%% file: Highcorr.tex
\documentclass{amsart}
\newtheorem{theorem}{Theorem}[section]
\newtheorem{lemma}[theorem]{Lemma}
\newtheorem{corollary}[theorem]{Corollary}
\theoremstyle{definition}

\theoremstyle{remark}

\numberwithin{equation}{section}

\newcommand{\scr}{\scriptstyle}
\def\sumprime_#1{\setbox0=\hbox{$\scriptstyle{#1}$}
\setbox2=\hbox{$\displaystyle{\sum}$}
\setbox4=\hbox{${}'\mathsurround=0pt$}
\dimen0=.5\wd0 \advance\dimen0 by-.5\wd2
\ifdim\dimen0>0pt
\ifdim\dimen0>\wd4 \kern\wd4 \else\kern\dimen0\fi\fi
\mathop{{\sum}'}_{\kern-\wd4 #1}}
\font\ger=eufm10
\newcommand{\gs}{\hbox{\ger S}}
\def\ndiv{\not \hskip .03in \mid}

\begin{document}
\include{Main}

\end{document}

%% file: Main.tex
\title[HIGHER CORRELATIONS OF DIVISOR SUMS I]{Higher Correlations of Divisor Sums Related to Primes I: Triple Correlations}
\author{D. A. Goldston}
\address{Department of Mathematics and Computer Science, San Jose
State University, San Jose, CA 95192, USA}
\email{goldston@mathcs.sjsu.edu}
\thanks{The first author was supported in part by an NSF Grant}
\author{C. Y. Y{\i}ld{\i}r{\i}m}
\address{ Department of Mathematics, Bilkent University, Ankara 06533, Turkey}
\email{yalcin@fen.bilkent.edu.tr}

\date{\today}


\begin{abstract}
We obtain the triple correlations for a truncated divisor sum related to primes. We also obtain the mixed correlations for this divisor sum when it is summed over the primes, and give some applications to primes in short intervals.
\end{abstract}

\maketitle

\section{  Introduction}

This is the first in a series of papers concerned with the calculation of higher correlations of short divisor sums that are approximations for the von Mangoldt function $\Lambda(n)$, where $\Lambda(n)$ is defined to be $\log p$ if $n=p^m$, $p$ a prime, $m$ a positive integer, and to be zero otherwise. These higher correlations have applications to the theory of primes which is our motivation for their study.
In this first paper we will calculate the pair and triple  correlations for
\begin{equation} \Lambda_R(n)= \sum_{\stackrel{\scr d |n}{\scr d\le R}}\mu(d) \log (R/d), \qquad \mathrm{ for} \ n\ge 1, \label{1.10} \end{equation}
and $\Lambda_R(n) = 0$ if $n\le 0$.
In later papers in this series we will examine quadruple and higher correlations, and also examine the more delicate divisor sum
\begin{equation} \lambda_R(n) = \sum_{r\leq R} \frac{\mu^2(r)}{\phi(r)} \sum_{\stackrel{\scr d|r}{ \scr d|n}} d\mu(d), \qquad \mathrm{ for} \ n\ge 1, \label{1.20} \end{equation}
 and  $\lambda_R(n)=0$ if $n\le 0$.
The correlations we are interested in evaluating are
\begin{equation} \mathcal{ S}_k(N,\mbox{\boldmath$j$}, \mbox{\boldmath $a$}) =\sum_{n=1}^N \Lambda_R(n+j_1)^{a_1}\Lambda_R(n+j_2)^{a_2}\cdots \Lambda_R(n+j_r)^{a_r}\label{1.30}\end{equation}
and
\begin{equation}\tilde{\mathcal{S}}_k(N,\mbox{\boldmath$j$}, \mbox{\boldmath$a$}) =\sum_{n=1}^N \Lambda_R (n+j_1)^{a_1}\Lambda_R (n+j_2)^{a_2}\cdots \Lambda_R(n+j_{r-1})^{a_{r-1}}\Lambda(n+j_r) \label{1.40}\end{equation}
where $\mbox{\boldmath$j$} = (j_1,j_2, \ldots , j_r)$ and $\mbox{\boldmath$a$} = (a_1,a_2, \ldots a_r)$, the $j_i$'s are distinct integers, $a_i\geq 1$ and $\sum_{i=1}^r a_i = k$. In \eqref{1.40} we assume that $r\ge 2$ and take $a_r=1$. For later convenience we define
\begin{equation} \tilde{\mathcal{ S}}_1(N, \mbox{\boldmath$j$}, \mbox{\boldmath$a$}) =\sum_{n=1}^N\Lambda(n+j_1) \sim N \label{1.50}\end{equation}
with $|j_1|\le N$ by the prime number theorem.
For $k=1$ and $k=2$ these correlations have been evaluated before \cite{GO2} (and for $\lambda_Q(n)$ they have been evaluated in \cite{GO3}); the results show that $\Lambda_R$ and $\lambda_R$ mimic the behavior of $\Lambda$, and this is also the case in arithmetic progressions, see \cite{HB}, \cite{Hoo}, \cite{GY}.

When $k\ge 3$ the procedure for evaluating these correlations is complicated, and it is easy to make mistakes in the calculations. Therefore we have chosen to first treat the triple correlations in detail. The main terms in the theorems can often be obtained in an easier way by evaluating the multiple sums in a different order or with a different decomposition of the initial summands; the method used here was chosen to control the error terms and generalize to higher correlations. Recently we have found a somewhat different method which is preferable for higher values of $k$. This method will be used in the third paper in this series.   We can not compute correlations which contain a factor  $\Lambda(n)\Lambda(n+k)$, $k\neq 0$,  without knowledge about prime twins. This limits our applications, and further the mixed correlations \eqref{1.40}  can only be calculated for shorter divisor sums than the pure correlations \eqref{1.30} of $\Lambda_R(n)$, which degrades to some extent the results we obtain. When we assume the Elliott-Halberstam conjecture we can eliminate this latter problem and obtain stronger results.

One motivation for the study of the correlations of $\Lambda_R(n)$ or $\lambda_R(n)$ is to provide further information on the moments
\begin{equation} M_k(N,h,\psi) = \sum_{n=1}^N(\psi(n+h)-\psi(n))^k \label{1.60} \end{equation}
where $\psi(x) = \sum_{n\le x}\Lambda(n)$.  We always take $N\to \infty$, and let
\begin{equation}h \sim \lambda \log N ,\label{1.70}\end{equation}
where we will usually be considering the case $\lambda \ll 1$. When $h$ is larger we need to subtract the expected value  $h$ in the moments above, which leads to more delicate questions which we will not consider in this paper (see \cite{MS}). Gallagher \cite{GA} proved that the moments in (\ref{1.60}) may be computed from the Hardy-Littlewood prime $r$-tuple conjecture \cite{HL}. This conjecture states that for $\mbox{\boldmath$j$} = (j_1,j_2, \ldots , j_r)$ with the $j_i$'s distinct integers,
\begin{equation}\psi_{\mbox{\boldmath$j$}}(N) = \sum_{n=1}^N \Lambda(n+j_1)\Lambda(n+j_2)\cdots \Lambda(n+j_r) \sim \gs(\mbox{\boldmath$j$}) N \label{1.80}\end{equation}
when $\gs(\mbox{\boldmath$j$})\neq 0$,
where
\begin{equation} \gs(\mbox{\boldmath$j$}) = \prod_p\left( 1- \frac{1}{p}\right)^{-r}\left(1-\frac{\nu_p(\mbox{\boldmath$j$})} {p}\right) \label{1.90}\end{equation}
and $\nu_p(\mbox{\boldmath$j$})$ is the number of distinct residue classes modulo $p$ that the $j_i$'s occupy. If $r=1$ we see $\gs(\mbox{\boldmath$j$})=1$, and for $|j_1|\le N$ equation (\ref{1.80}) reduces to \eqref{1.50}, which is the only case where (\ref{1.80}) has been proved. To compute the moments in (\ref{1.60}) we have
\begin{eqnarray*} M_k(N,h,\psi) &=& \sum_{n=1}^N \left(\sum_{1\le m\le h}\Lambda(n+m)\right)^k \\
&=& \sum_{\stackrel{\scr 1\le m_i\le h }{\scr 1\le i\le k}}\sum_{n=1}^N \Lambda(n+m_1)\Lambda(n+m_2)\cdots \Lambda(n+m_k).
\end{eqnarray*}
Now suppose that the $k$ numbers $m_1, m_2, \ldots , m_k$ take on $r$ distinct values
$j_1, j_2, \ldots , j_r$ with $j_i$ having multiplicity $a_i$, so that $\sum_{1\le i \le r}a_i =k$. Grouping the terms above, we have that
\begin{equation} M_k(N,h,\psi) =  \sum_{r=1}^k \sum_{\stackrel{\scr
a_1 , a_2, \ldots, a_r }{ \scr a_i \ge 1, \sum a_i =k}}
    \left(
        \begin{array}{c} k \\  a_1 , a_2 , \ldots , a_r         \end{array}
    \right)
\sum_{1\le j_1 < j_2 < \cdots < j_r \le h }\psi_k(N,\mbox{\boldmath$j$}, \mbox{\boldmath$a$}),\label{1.100}
\end{equation}
where
\begin{equation}\psi_k(N,\mbox{\boldmath$j$}, \mbox{\boldmath$a$})=\sum_{n=1}^N \Lambda(n+j_1)^{a_1}\Lambda(n+j_2)^{a_2} \cdots \Lambda(n+j_r)^{a_r}\label{1.110} \end{equation}
and the multinomial coefficient
counts the number of different innermost sums that occur. If $n + j_i$ is a prime then $\Lambda(n+j_i)^{a_i} = \Lambda(n+j_i)(\log (n+j_i))^{a_i -1},$
and we easily see that
\begin{equation} \begin{split}\psi_k(N,\mbox{\boldmath$j$}, \mbox{\boldmath$a$}) & =(1+o(1))  (\log N)^{k-r} \sum_{n=1}^N \Lambda(n+j_1)\Lambda(n+j_2) \cdots \Lambda(n+j_r) + O(N^{\frac{1}{2} +\epsilon})\\ &= (1+o(1))(\log N)^{k-r}\psi_{\mbox{\boldmath$j$}}(N) + O(N^{\frac{1}{2} +\epsilon}).\label{1.120}
\end{split}
\end{equation}
Hence we may apply the conjecture (\ref{1.80}) assuming it is valid uniformly for $\max_i|j_i| \le h$ and obtain
\[ M_k(N,h,\psi) \sim  N\sum_{r=1}^k (\log N)^{k-r}\sum_{\stackrel{\scr a_1 , a_2, \ldots, a_r}{ \scr a_i \ge 1, \sum a_i =k}}
    \left(
        \begin{array}{c} k \\  a_1 , a_2 , \ldots , a_r         \end{array}
    \right)
\sum_{ 1\le j_1 < j_2 < \cdots < j_r \le h}\gs(\mbox{\boldmath$j$}).\]
Gallagher \cite{GA} proved that, as $h\to \infty$,
 \begin{equation}\sum_{\stackrel{\scr 1\le j_1 , j_2 , \cdots ,j_r \le h}  {\scr  \mathrm{ distinct}}}\gs(\mbox{\boldmath$j$})\sim h^r, \label{1.130}
\end{equation}
and since this sum includes $r!$ permutations of the specified vector $\mbox{\boldmath$j$}$ when the components are ordered, we have
\[M_k(N,h,\psi) \sim N(\log N)^k \sum_{r=1}^k \frac{1}{r!}(\frac{h}{\log N})^{r}\sum_{\stackrel{\scr  a_1 , a_2, \ldots, a_r} { \scr a_i \ge 1, \sum a_i =k}}
    \left(
        \begin{array}{c} k \\  a_1 , a_2 , \ldots , a_r         \end{array}
    \right)
.\]
 Letting
\( \left\{
        \begin{array}{c} k \\ r
        \end{array}
    \right\}\)
denote the Stirling numbers of the second type, then it may be easily verified (see \cite{GKP}) that
\begin{equation} \sum_{\stackrel{\scr a_1 , a_2, \ldots, a_r}{\scr a_i \ge 1, \sum a_i =k}}
    \left(
        \begin{array}{c} k \\  a_1 , a_2 , \ldots , a_r         \end{array}
    \right)
 = r! \left\{
        \begin{array}{c} k \\ r
        \end{array}
    \right\}
.\label{1.140}\end{equation} We conclude that for $h\sim \lambda \log N$,
\begin{equation} M_k(N,h,\psi) \sim N(\log N)^k \sum_{r=1}^k
 \left\{
        \begin{array}{c} k \\ r
        \end{array}
    \right\}
\lambda^r ,\label{1.150}\end{equation}
which are the moments of a Poisson distribution with mean $\lambda$.
The first 4 moments are, for $\lambda \ll 1$,
\begin{eqnarray*}  M_1(N,h,\psi)& \sim &\lambda N\log N, \\
  M_2(N,h,\psi)& \sim &(\lambda +\lambda^2)N\log^2 N, \\
M_3(N,h,\psi)& \sim &(\lambda + 3\lambda^2 +\lambda^3)N\log^3 N, \\
M_4(N,h,\psi)& \sim &(\lambda +7 \lambda^2 +6\lambda^3 + \lambda^4)N\log^4N .\end{eqnarray*}
The asymptotic formula for the first moment is known to be true as a simple consequence of the prime number theorem. The other moment formulas have never been proved. It is known that the asymptotic formula for the second moment follows from the assumption of the Riemann Hypothesis and the pair correlation conjecture for zeros of the Riemann zeta-function \cite{GM}.

Turning to our approximation $\Lambda_R(n)$, we define
\begin{equation} \psi_R(x) = \sum_{n\le x}\Lambda_R(n)\label{1.160} \end{equation}
and first wish to examine the moments $M(N,h,\psi_R)$ defined as in
(\ref{1.60}). The same computation used for $\psi$ to obtain (\ref{1.100}) clearly applies and therefore we obtain
\begin{equation} M_k(N,h,\psi_R) =  \sum_{r=1}^k \sum_{\stackrel{\scr a_1 , a_2, \ldots, a_r}{\scr a_i \ge 1, \sum a_i =k}}
    \left(
        \begin{array}{c} k \\  a_1 , a_2 , \ldots , a_r         \end{array}
    \right)
\sum_{1\le j_1 < j_2 < \cdots < j_r \le h } \mathcal{S}_k(N,\mbox{\boldmath$j$}, \mbox{\boldmath$a$}), \label{1.170}\end{equation}
where $\mathcal{S}_k(N,\mbox{\boldmath$j$}, \mbox{\boldmath$a$})$ is the correlation given in (\ref{1.30}).  Since $\Lambda_R(n)$ is not supported on the primes and prime powers as $\Lambda(n)$ is, we can not use (\ref{1.120}) to reduce the problem to correlations without powers, and as we will see these powers sometimes effect the correlations for $\Lambda_R(n)$. Our main result on these correlations is contained in the following theorem.

\begin{theorem} \label{Theorem1} Given $1\le k\le 3$,  let $\mbox{\boldmath$j$} = (j_1,j_2, \ldots , j_r)$ and $\mbox{\boldmath$a$} = (a_1,a_2, \ldots a_r)$, where the $j_i$'s are distinct integers, and $a_i\geq 1$ with $\sum_{i=1}^r a_i = k$. Assume $\max_i |j_i|\ll R^\epsilon$ and $R\gg N^\epsilon$. Then
we have
\begin{equation} \mathcal{S}_k(N,\mbox{\boldmath$j$},\mbox{\boldmath$a$}) = \big(\mathcal{ C}_k(\mbox{\boldmath$a$})\gs(\mbox{\boldmath$j$})+o(1)\big)N(\log R)^{k-r} +O(R^k),\label{1.180}\end{equation}
where $\mathcal{ C}_k(\mbox{\boldmath$a$})$ has the values
\begin{eqnarray}  \mathcal{ C}_1(1)&=&1, \nonumber \\
\mathcal{ C}_2(2)&=&1, \quad  \mathcal{ C}_2(1,1)=1, \nonumber \\
 \mathcal{ C}_3(3)&=& \frac{3}{ 4}, \quad \mathcal{ C}_3(2,1)=1, \quad \mathcal{ C}_3(1,1,1)=1.\nonumber \end{eqnarray}
\end{theorem}

Here we have used the notational convention of dropping extra parentheses, so for example $\mathcal{ C}_2((1,1))=\mathcal{C}_2(1,1)$.
The method of proof used in this paper is not limited to $k\le 3$, but it does become extremely complicated even for $k=4$. Using the new method mentioned before we will prove Theorem \ref{Theorem1} holds for all $k$ in the third paper in this series. The computation of the constants $\mathcal{ C}_k(\mbox{\boldmath$a$})$ as $k$ gets larger becomes increasingly difficult. We also believe the error term
$O(R^k)$ can be improved. This has been done in the case $k=2$ (unpublished) where the error term $O(R^2)$ may be replaced by $O(R^{2-\delta})$ for a small constant $\delta$. In the special case of $S_2(N,(0),(2))$ Graham \cite{Gr} has removed the error term $O(R^2)$ entirely.

In proving Theorem \ref{Theorem1} we will assume $j_1=0$. This may be done without loss of generality since we may shift the sum over $n$ in $\mathcal{ S}_k$ to $m=n +j_1$ and then return to the original summation range with an error $O(N^\epsilon)$ since $|j_1|\ll R^\epsilon$ and $\Lambda_R(n)\ll n^\epsilon$. Further $\gs(\mbox{\boldmath$j$})= \gs(\mbox{\boldmath$j$}-\mbox{\boldmath$j_1$})$ where $\mbox{\boldmath$j_1$}$ is a vector with $j_1$ in every component, and so the singular series are unaffected by this shift.

We now apply Theorem \ref{Theorem1} in \eqref{1.170}, and obtain immediately using \eqref{1.130} that
\begin{equation} M_k(N,h,\psi_R) =  (1+o(1))P_k(\lambda, R) N(\log R)^k  +O(h^kR^k)\label{1.190}
\end{equation}
where
\begin{equation} P_k(\lambda, R)=
\sum_{r=1}^k \frac{1}{r!} \left(\frac{h}{ \log R}\right)^r
\sum_{\stackrel{\scr a_1 , a_2, \ldots, a_r }{\scr a_i \ge 1, \sum a_i =k}}
    \left(
        \begin{array}{c} k \\  a_1 , a_2 , \ldots , a_r         \end{array}
    \right)
\mathcal{ C}_k(\mbox{\boldmath$a$})  .\label{1.200} \end{equation}
Using the values of the constants $\mathcal{ C}(\mbox{\boldmath$a$})$ in Theorem \ref{Theorem1} we obtain the following result on moments of $\psi_R$.

\begin{corollary} For $h\sim \lambda \log N$, $\lambda \ll R^\epsilon$, and $R=N^{\theta_k}$,  where $\theta_k$ is fixed and $0<\theta_k<\frac{1}{k}$ for $ 1\le k \le 3$, we have
 \begin{eqnarray*} M_1(N,h,\psi_R)& \sim& \lambda N\log N,   \\
 M_2(N,h,\psi_R)& \sim &(\theta_2 \lambda +\lambda^2)N\log^2 N, \\
 M_3(N,h,\psi_R)& \sim &(\frac{3}{4}{\theta_3}^2 \lambda + 3\theta_3 \lambda^2 +\lambda^3)N\log^3 N.
\end{eqnarray*}
\label{Corollary1}
\end{corollary}

We next consider the mixed moments
\begin{equation} \tilde{M}_k(N,h,\psi_R) = \sum_{n=1}^N(\psi_R(n+h)-\psi_R(n))^{k-1}
(\psi(n+h)-\psi(n))\label{1.210}\end{equation}
for $k\ge 2$, while if $k=1$ we take
\begin{equation}\tilde{M}_1(N,h,\psi_R)= M_1(N,h,\psi) \sim \lambda N \log N \label{1.220}\end{equation}
for $1\le R\le N$ by the prime number theorem.
Now assume $k \ge 2$. On multiplying out and grouping as before for the terms involving $\Lambda_R$ we have
\begin{equation} \tilde{M}_k(N,h,\psi_R) =  \sum_{r=2}^{k} \sum_{\stackrel{\scr
a_1 , a_2, \ldots, a_{r-1}}{ \scr a_i \ge 1, \sum a_i =k-1}}
\frac{1}{(r-1)!}
    \left(
        \begin{array}{c} k-1 \\  a_1 , a_2 , \ldots , a_{r-1}       \end{array}
    \right)
\sum_{\stackrel{\scr 1\le j_1 ,j_2 ,\cdots , j_{r-1} \le h}{ \scr  {\rm distinct}}}V_{r-1}(N,\mbox{\boldmath$j$},\mbox{\boldmath$a$}),\label{1.230}\end{equation}
where
\[V_{r-1}(N,\mbox{\boldmath$j$}, \mbox{\boldmath$a$}) = \sum_{1\le m\le h}\sum_{n=1}^N \Lambda_R(n+j_1)^{a_1}\Lambda_R(n+j_2)^{a_2} \cdots \Lambda_R(n+j_{r-1})^{a_{r-1}}\Lambda(n+m)\]
Since, provided  $n +j>R$ and $n+j\neq p^m$ for some $m\ge 2$ and $p<R$, we have
\[ \Lambda_R(n+j)^a\Lambda(n+j) = (\log R)^a\Lambda(n+j), \]
we see
\begin{equation*}\begin{split} V_{r-1}&(N,\mbox{\boldmath$j$}, \mbox{\boldmath$a$})\\ &= \sum_{i=1}^{r-1}(\log R)^{a_i}\sum_{n=1}^N \Big(\prod_{\stackrel{\scr 1\le s\le r-1}{ \scr s\neq i}}\Lambda_R(n+j_s)^{a_s}\Big)\Lambda(n+j_i)\\
&+ \sum_{\stackrel{\stackrel{\scr 1\le j_r \le h}{ \scr  j_r\neq j_i}} { \scr 1\le i \le r-1}}
\sum_{n=1}^N \Lambda_R(n+j_1)^{a_1}\Lambda_R(n+j_2)^{a_2} \cdots \Lambda_R(n+j_{r-1})^{a_{r-1}}\Lambda(n+j_r) + O( RN^\epsilon)\\
&=\sum_{i=1}^{r-1}(\log R)^{a_i}\tilde{\mathcal{ S}}_{k-a_i}(N,\mbox{\boldmath$j_i$}, \mbox{\boldmath$a_i$})+\sum_{
\stackrel{
\stackrel{\scr 1\le j_l \le h}{\scr  j_l\neq j_i}}
{\scr 1\le i \le r-1}}
\tilde{\mathcal{ S}}_k(N,\mbox{\boldmath$j$}, \mbox{\boldmath$a$})
+ O( RN^\epsilon)\end{split}\end{equation*}
where
\begin{equation} \mbox{\boldmath$j_i$}= (j_1,j_2, \ldots ,j_{i-1},j_{i+1}, \ldots , j_{r-1}, j_i),     \quad \mbox{\boldmath$a_i$}= (a_1,a_2,\ldots , a_{i-1},a_{i+1}, \ldots , a_{r-1}, 1). \label{1.240} \end{equation}
We conclude for $k\ge 2$
\begin{equation} \tilde{M}_k(N,h,\psi_R) =  \sum_{r=2}^{k} \sum_{\stackrel{\scr
a_1 , a_2, \ldots, a_{r-1}}{\scr a_i \ge 1, \sum a_i =k-1}}
\frac{1}{(r-1)!}
    \left(
        \begin{array}{c} k-1 \\  a_1 , a_2 , \ldots , a_{r-1}       \end{array}
    \right)
W_r(N,\mbox{\boldmath$j$}, \mbox{\boldmath$a$}) + O( RN^\epsilon),\label{1.250}\end{equation}
where
\begin{equation} W_r(N,\mbox{\boldmath$j$}, \mbox{\boldmath$a$})= \sum_{i=1}^{r-1} (\log R)^{a_i} \sum_{\stackrel{\scr 1\le j_1 ,j_2 ,\cdots , j_{r-1} \le h}{ \scr  \mathrm{ distinct}}} \tilde{\mathcal{ S}}_{k-a_i}(N,\mbox{\boldmath$j_i$}, \mbox{\boldmath$a_i$}) +
\sum_{\stackrel{\scr 1\le j_1 ,j_2 ,\cdots , j_{r} \le h}{ \scr  \mathrm{ distinct}}}
\tilde{\mathcal{ S}}_{k}(N,\mbox{\boldmath$j$}, \mbox{\boldmath$a$}) .\label{1.260}\end{equation}

We have reduced the calculation of the mixed moments to mixed correlations. Our method for evaluating the mixed correlations will prove as a by-product that the mixed correlations are asymptotically equal to the corresponding pure correlations in a certain range of $R$. Our results depend on the uniform distribution of primes in arithmetic progressions. We let
\begin{equation} \psi(x;q,a) = \sum_{\stackrel{\scr n\le x}{ \scr n\equiv a (q)} }\Lambda(n),\label{1.270}\end{equation}
and
\begin{equation}    E_{a,b} = \left\{
        \begin{array}{ll}
        1, &\mbox{if $(a,b)=1$,}\\
                  0, & \mbox{if $(a,b)>1$}.
        \end{array}
    \right.
\label{1.280}\end{equation}
On taking
\begin{equation} E(x;q,a)= \psi(x;q,a) - E_{a,q}\frac{x}{ \phi(q)},\label{1.290}\end{equation}
the estimate we need is, for some fixed $0<\vartheta \le 1$,
\begin{equation} \sum_{1\le q\le x^{\vartheta - \epsilon}} \max_{\stackrel{\scr a }{ \scr (a,q)=1}}|E(x;q,a)|  \ll \frac{x}{\log ^\mathcal{A}x},\label{1.300} \end{equation}
for any $\epsilon >0$, any $\mathcal{A}=\mathcal{A}(\epsilon)>0$, and $x$ sufficiently large.  This is a weakened form of the Bombieri-Vinogradov theorem if $\vartheta=\frac{1}{2}$, and therefore \eqref{1.300} holds unconditionally if $\vartheta \le \frac{1}{2}$.  Elliott and Halberstam conjectured \eqref{1.300} is true with $\vartheta=1$. The range of $R$ where our results on mixed correlations hold depends on $\vartheta$ in (\ref{1.300}).  We first prove the following general result.
\begin{theorem} \label{Theorem2} Given $k\ge 2$,  let $\mbox{\boldmath$j$} = (j_1,j_2, \ldots , j_r)$ and $\mbox{\boldmath$a$} = (a_1,a_2, \ldots a_r)$, where the $j_i$'s are distinct integers, $a_i\geq 1$, $a_r=1$, and $\sum_{i=1}^r a_i = k$ . Assume $\max_i |j_i|\ll R^{\frac{1}{k}}$ and that $ N^\epsilon\ll R\ll N^{\min(\frac{\vartheta}{ k-1}, \frac{1}{k})-\epsilon}$. Then
we have, with $\mathcal{A}$ from \eqref{1.300},
\begin{equation}
\tilde{\mathcal{S}}_k(N,\mbox{\boldmath$j$},\mbox{\boldmath$a$}) = \mathcal{ S}_k(N,\mbox{\boldmath$j$},\mbox{\boldmath$a$})
+ O(R^{k}) +O\big(\frac{N}{(\log N)^{\frac{\mathcal{A}}{2} - 4^{k-3/2}+2^{k-1}-1 }}\big).
\label{1.310}
\end{equation}
\end{theorem}

The proof of Theorem \ref{Theorem2} involves proving that both $\tilde{\mathcal{S}}_k(N,\mbox{\boldmath$j$},\mbox{\boldmath$a$})$ and $\mathcal{ S}_k(N,\mbox{\boldmath$j$},\mbox{\boldmath$a$})$ are asymptotic to the same main term and therefore they are asymptotic to each other in the  range where both asymptotic formulas hold.
Using Theorems \ref{Theorem1} and \ref{Theorem2} we can now immediately evaluate the mixed moments. There is, however, an inefficiency in the use of Theorem \ref{Theorem2} which imposes the condition that $R \ll N^{\min(\frac{\vartheta}{k-1}, \frac{1}{k}) -\epsilon}$. The restriction $R\ll N^{\frac{1}{k}-\epsilon}$ in this condition arises from applying Theorem \ref{Theorem1}, but by directly evaluating  the main term that arises in the proof of Theorem \ref{Theorem2} we can remove this condition and prove the following result.

\begin{theorem} \label{Theorem3} Given $2\le k\le 3$,  let $\mbox{\boldmath$j$} = (j_1,j_2, \ldots , j_r)$ and $\mbox{\boldmath$a$} = (a_1,a_2, \ldots a_r)$, where $r\ge 2$, $a_r=1$, and where the $j_i$'s are distinct integers, and $a_i\geq 1$ with $\sum_{i=1}^r a_i = k$. Assume $\max_i |j_i|\ll R^\epsilon$. Then
we have, for $N^{\epsilon}\ll R \ll N^{\frac{\vartheta}{k-1} -\epsilon} $ where \eqref{1.300} holds with $\vartheta$,
\begin{equation} \tilde{\mathcal{S}}_k(N,\mbox{\boldmath$j$},\mbox{\boldmath$a$}) = \big(\gs(\mbox{\boldmath$j$}) +o(1)\big)N(\log R)^{k-r} .\label{1.320}\end{equation}
\end{theorem}
For larger $k$  the constants $\mathcal{ C}(\mbox{\boldmath$a$})$ will appear in this theorem as in Theorem \ref{Theorem1}, but for $k\le 3$ all these constants for the mixed correlations are equal to 1.

Next, using \eqref{1.250} we are able to evaluate the first three mixed moments.

\begin{corollary} \label{Corollary2} For $h\sim \lambda \log N$, $\lambda \ll R^\epsilon$, and $R=N^{\theta_k}$,  where $\theta_k$ is fixed, $0<\theta_1 \le 1$,  and  $0<\theta_k<\frac{\vartheta}{k-1}$ for $ 2\le k \le 3$, we have,
\begin{eqnarray*} \tilde{M}_1(N,h,\psi_R)& \sim &\lambda N\log N,  \\
  \tilde{M}_2(N,h,\psi_R)& \sim &(\theta_2 \lambda +\lambda^2)N\log^2 N, \\
\tilde{M}_3(N,h,\psi_R)& \sim& ({\theta_3}^2 \lambda + 3\theta_3 \lambda^2 +\lambda^3)N\log^3 N.
\end{eqnarray*}
\end{corollary}

The starting point of Bombieri and Davenport's \cite{BD} work on small gaps between primes  is essentially equivalent to the inequality
\begin{equation}
\sum_{n=N+1}^{2N}\Big(\big(\psi(n+h)-\psi(n)\big)-
\big(\psi_R(n+h)-\psi_R(n)\big)\Big)^2 \ge 0.\label{1.330}
\end{equation}
Letting
\begin{equation} M_k'(N,h,\psi) = M_k(2N,h,\psi)-M_k(N,h,\psi),
\label{1.340}
\end{equation}
with the corresponding definition for $M_k'(N,h,\psi_R)$ and $\tilde{M}_k'(N,h,\psi_R)$, we see that Corollary \ref{Corollary1} holds with $M_k(N,h,\psi_R)$ replaced by $M_k'(N,h,\psi_R)$ and Corollary
\ref{Corollary2} holds with $\tilde{M}_k(N,h,\psi_R)$ replaced with $\tilde{M}_k'(N,h,\psi_R)$.
On expanding \eqref{1.330} we have
\[ M_2'(N,h,\psi) \ge  2 \tilde{M}_2'(N,h,\psi_R) - M_2'(N,h,\psi_R)  \]
which implies on taking $\theta_2 = 1/2-\epsilon$ in Corollaries \ref{Corollary1} and  \ref{Corollary2} that
\begin{equation}
M_2'(N,h,\psi) \ge ( \frac{1}{2} \lambda + \lambda^2 -\epsilon)N\log^2N.
\label{1.350}
\end{equation}
 Let $p_j$ denote the $j$-th prime. If it is the case that $p_{j+1}-p_j > h=\lambda \log N$ for all $\frac{N}{2} < p_j \le 2N$, then each of the intervals $(n,n+h]$ for $N<n\le 2N$ contains at most one prime. Hence, since the prime powers may be removed with a negligible error, we have that $M_2'(N,h,\psi) \sim (\log N) M_1'(N,h,\psi)\sim \lambda N\log^2N$ so that \eqref{1.350} implies
\[  \lambda \ge \frac{1}{2} \lambda + \lambda^2 -\epsilon \]
which is false if $\lambda > \frac{1}{2}$.
We conclude that
\begin{equation} \liminf_{n\to \infty} \left( \frac{p_{n+1}-p_n}{\log p_n}\right) \le \frac{1}{2} .\label{1.360}\end{equation}
More generally, we define for $r$ any positive integer
\begin{equation}  \Xi_r =  \liminf_{n\to \infty} \left( \frac{p_{n+r}-p_n}{\log p_n}\right) \label{1.370}
\end{equation}
and see that
if $p_{n+r}-p_n > h =\lambda \log N$ for $N< p_n \le 2N$ then
\[ M_2'(N,h,\psi) \le (r+\epsilon)(\log N) M_1'(N,h,\psi) \le  (r+\epsilon)\lambda N\log ^2N \]
which then implies that
\[ r\lambda \ge \frac{1}{2} \lambda + \lambda^2 -\epsilon \]
and hence
\begin{equation} \Xi_r \le r -\frac{1}{2} .\label{1.380}\end{equation}
Bombieri and Davenport were also able to improve \eqref{1.360} by incorporating an earlier method of Erd\"os into their argument. This method depends on the sieve upper bound for primes differing by an even number $k$ given by
\begin{equation} \sum_{n\le N}\Lambda(n)\Lambda(n+k) \le (\mathcal{B}+ \epsilon)\gs(k) N \label{1.390}
\end{equation}
where $\gs(k) = \gs(\mbox{\boldmath$j$})$ with $\mbox{\boldmath$j$}=(0,k)$, and $\mathcal{B}$ is a constant. In \cite{BD} Bombieri and Davenport proved that \eqref{1.390} holds with $\mathcal{B}=4$, and using this value  they improved \eqref{1.360} and obtained
\begin{equation} \Xi_1 \le \frac{2 +\sqrt{3}}{8}= 0.46650\ldots  .\label{1.400}\end{equation}
While \eqref{1.350} has never been improved, the refinements based on the Erd\"os method together with the choice of certain weights in a more general form of \eqref{1.350} has led to further improvements. Huxley \cite{Hu1} \cite{Hu2} proved that,
letting $\theta_r$ be the smallest positive solution of
\begin{equation}\theta_r +\sin\theta_r = \frac{\pi}{\mathcal{B} r},\quad \sin\theta_r < (\pi + \theta_r) \cos \theta_r \ , \label{1.410} \end{equation}
then
\begin{equation}  \Xi_r \le \frac{2r-1}{4\mathcal{B}r}\left\{ \mathcal{B}r +(\mathcal{B}r-1)\frac{\theta_r}{\sin\theta_r}\right\}.\label{1.420}
\end{equation}
With the value $\mathcal{B}=4$ this gives
\[ \Xi_1 \le 0.44254\ldots , \quad \Xi_2 \le 1.41051\ldots, \quad \Xi_3 \le 2.39912\ldots , \quad \Xi_4 \le 3.39326\ldots. \]
We note that the expression on the right-hand side of \eqref{1.420} is equal to
\[  r-\frac{1+\frac{1}{\mathcal{B}}}{2} + O(\frac{1}{r}), \]
and thus for large $r$ this bound approaches $r- \frac{5}{8}$ with $\mathcal{B}=4$.

The best result known for $\mathcal{B}$ which holds uniformly for all $k$ is
$\mathcal{B}=3.9171\ldots$ due to Chen \cite{Ch}.
However, in the application to obtain \eqref{1.420} one only needs \eqref{1.390} to hold uniformly for a restricted range of $k$; the condition $0 < |k| \le \log^2N$ is more than sufficient. In this case there have been a string of improvements. For ease of comparison with the value $\mathcal{B}= 4$ used above, the value $\mathcal{B}= 3.5$  obtained by Bombieri, Friedlander, and Iwaniec \cite{BFI} gives the values
\[ \Xi_1 \le 0.43493\ldots , \quad \Xi_2 \le 1.39833\ldots, \quad \Xi_3 \le 2.38519\ldots , \quad \Xi_4 \le 3.37842\ldots. \]
All of these results above actually hold for a positive percentage of gaps.  Maier \cite{Ma} introduced a new method to prove  that
\begin{equation} \liminf_{n\to \infty} \left( \frac{p_{n+1}-p_n}{\log p_n}\right) \le e^{-\gamma} =  0.56145\ldots \ .\label{1.430}\end{equation}
This method, which applies to special sets of sparse intervals, may be combined with the earlier methods to include this factor of $e^{-\gamma}$ times the earlier results. The argument was carried out with $\mathcal{B}=4$ in \cite{Ma}, which then gives in particular
\[ \Xi_1 \le 0.24846\ldots , \quad \Xi_2 \le 0.79194\ldots, \quad \Xi_3 \le 1.34700\ldots , \quad \Xi_4 \le 1.90518\ldots. \]

Our approach for examining gaps between primes is to consider the mixed third moment
\begin{equation}
\tilde{M}_3'(N, h,\psi_R,C) =  \sum_{n=N+1}^{2N}\big(\psi(n+h)-\psi(n)\big)
\big(\psi_R(n+h)-\psi_R(n)- C\log N\big)^2 ,\label{1.440}
\end{equation}
Here $C$   may be chosen as a function of $h$ and $R$ to optimize the argument. The idea behind the use of  $\tilde{M}_3'(N, h,\psi_R,C)$ is that it will approximate and thus provide some of the same information as the third moment $M_3'(N,h,\psi)$.
If $p_{n+r}-p_n >h=\lambda \log N$ for all $N<p_n\le 2N$ then, removing prime powers as before, we have
\begin{equation}\tilde{M}_3'(N, h,\psi_R,C) \le (r +\epsilon) \log N  \sum_{n=N+1}^{2N}
\big(\psi_R(n+h)-\psi_R(n)- C\log N\big)^2 .\label{1.450}
\end{equation}
Corollaries \ref{Corollary1} and \ref{Corollary2} allow us to evaluate both sides of \eqref{1.450}, and on choosing $C$ appropriately we are able to prove the following result.

\begin{theorem} For $r\ge 1$, we have
\begin{equation} \label{1.460}
\Xi_r \le r -\frac{1}{2}\sqrt{r}.
\end{equation}
Further, assuming that  for $h\ll \log N$ and $R\le N^{\frac{1}{4}}$
\begin{equation}
M_4(N,h,\psi_R) \ll N\log^4N,
\label{1.470}
\end{equation}
then we have that, for $h = \lambda \log N$ and $\lambda > r -\frac{1}{2}\sqrt{r}$,
\begin{equation}
\sum_{\substack{N+1\le p_n \le 2N \\ p_{n+r}-p_{n}< h}}1 \gg_\lambda
\frac{N}{\log N} .\label{1.480}
\end{equation}
\label{Theorem5}
\end{theorem}

Thus we have
\[ \Xi_1 \le \frac{1}{2} , \quad \Xi_2 \le 1.29289\ldots, \quad \Xi_3 \le 2.13397\ldots , \quad \Xi_4 \le 3. \]
We see that our result improves on the results of Huxley when $r\ge 2$, although Maier's results are still better. Our theorem corresponds to \eqref{1.380} in that it does not use the Erd\"os method.
It is possible to incorporate the Erd\"os method into our method too, but this requires we first obtain an asymptotic formula for $M_4(N,h,\psi_R)$.  One should also be able to incorporate Maier's method as well, which would then give better results than are currently known for $r\ge2$.

The result in \eqref{1.480} shows that the small gaps produced in the theorem form a positive proportion of all the gaps assuming that \eqref{1.470} holds. We will prove \eqref{1.470} in a later paper in this series and thus show that \eqref{1.480} holds unconditionally.

We will actually prove
\begin{equation} \label{1.490}
\Xi_r \le r -\sqrt{\frac{\vartheta r}{2}},
\end{equation}
where $\vartheta$ is the number in \eqref{1.300}. Therefore, assuming the Elliott-Halberstam Conjecture in the form that one may take $\vartheta =1$ in \eqref{1.300}, we have
\[ \Xi_r \le r -\sqrt{\frac{r}{2}} \]
which in particular gives
\[ \Xi_1 \le 0.29289\ldots , \quad \Xi_2 \le 1, \quad \Xi_3 \le 1.77525\ldots , \quad \Xi_4 \le 2.58578\ldots. \]
These results are in contrast to the method of Bombieri and Davenport where the Elliott-Halberstam conjecture does not improve their results directly. (The Elliott-Halberstam conjecture does allow one to take $\mathcal{B}=2$ in \eqref{1.390}, and therefore leads to small improvements in Huxley's results, which for $r\ge 2$ are weaker than the result in Theorem \ref{Theorem5}.)
We can not extend these last results obtained assuming an Elliott-Halberstam conjecture to a positive proportion of gaps because we can not prove \eqref{1.470} for $\theta >\frac{1}{4}$.  Our proof gives that the number of gaps we produce in this case is $\gg N\log^{-B} N$ for some positive constant $B>1$.

Our method can also be used to examine larger than average gaps between primes. In this case much more is known than for small gaps; the latest result being that \cite{Pi}
\begin{equation}
\max_{p_n \le N}\frac{ p_{n+1}-p_n}{\log p_n} \ge (2e^\gamma - o(1)) \frac{\log\log N \log \log \log \log N}{(\log \log \log N)^2}.\label{1.500}
\end{equation}
If one were to ask however for a positive proportion of gaps larger than the average, then it is a remarkable fact that nothing non-trivial is known. \footnote{The first-named author of this paper learned of this from Carl Pomerance after a talk in which the author had claimed to have such a result.} What can be proved is that a positive proportion of the interval $(N,2N]$ is contained between consecutive primes whose difference is a little larger than average. To formalize this, we let $\Theta_r$ be the supremum over all $\lambda$ for which
\begin{equation}  \sum_{\substack{N<p_n\le 2N\\ p_{n+r}-p_n \ge \lambda \log N }}(p_{n+r}-p_n) \gg_\lambda  N \label{1.510}
\end{equation}
for all sufficiently large $N$.
Then using the Erd\"os method one finds that \cite{C-G}
\begin{equation}  \Theta_1 \ge 1 + \frac{1}{2\mathcal{B}}
\label{1.520}
\end{equation}
where $\mathcal{B}$ is the number in \eqref{1.390}.

To apply \eqref{1.440} to this problem, we assume that $p_{j+r}-p_j <h=\lambda \log N$ for all $N<p_j\le 2N$ in which case the interval $(n,n+h]$ always contains at least $r$ primes, and therefore \eqref{1.450} holds with the inequality reversed. On optimizing $C$ we obtain the following result.

\begin{theorem} \label{Theorem6} Assume that \eqref{1.470} holds. For $r\ge 1$
 we have that
 \begin{equation}\label{1.530}
\Theta_r \ge r +\frac{1}{2}\sqrt{r}.
\end{equation}
\end{theorem}

As mentioned above, we will prove \eqref{1.470} in a later paper in this series, which will show that  Theorem \ref{Theorem6}  holds  unconditionally.

The proof of Theorems \ref{Theorem5} and \ref{Theorem6} only require the asymptotic formula for the third mixed moment in Corollary \ref{Corollary2} and the second moment for $\psi_R$ in Corollary \ref{Corollary1}. Thus the results in sections 6--10 which are concerned with triple correlations of $\psi_R$ may be skipped by the reader who is only interested in our applications to primes.

\bigskip

\textit{Notation.}
In this paper $N$ will always be a large integer, $p$ denotes a prime number, and sums will start at 1 if a lower limit is unspecified. When a sum is denoted with a dash as $\sumprime_{} $ this always indicates we will sum over all variables expressed by inequalities in the conditions of summation and these variables will all be pairwise relatively prime with each other. We will always take the value of a void sum to be zero and the value of a void product to be 1. The letter $\epsilon$ will denote a small positive number which may change in each equation.
We will also use the Iverson notation \cite{GKP} that putting brackets around a true-false statement will replace the statement by 1 if it is true, and 0 if it is false:
\begin{equation}
    [P(x)]=
        \left\{
        \begin{array}{ll}
        1, &\mbox{if $P(x)$ is true,}\\
                  0, & \mbox{if $P(x)$ is false}.
        \end{array}
        \right. \label{1.540}
    \end{equation}
As usual, $(a,b)$ denotes the greatest common divisor of $a$ and $b$ and $[a_1,a_2, \cdots , a_n]$ denotes the least common multiple of $a_1,a_2, \ldots , a_n$.
\section{Lemmas}

For $j$ a non-negative integer, we define the arithmetic function $\phi_j(n)$  on the primes by
    \begin{equation}
\phi_j(p) = p-j, \label{2.10}
    \end{equation}
$\phi_j(1)=1$, and extend the definition to squarefree integers by multiplicativity.
 Thus for $n$ squarefree $\phi_0(n) = n$, and $\phi_1(n) = \phi(n)$. We will not need to extend the definition beyond the squarefree integers here.
Letting
    \begin{equation}
    p(j) =
        \left\{
        \begin{array}{ll}
        j, &\mbox{if $j$ is a prime,}\\
                  1, & \mbox{otherwise},
        \end{array}
        \right. \label{2.20}
    \end{equation}
we next define
    \begin{equation}
 H_j(n) = \prod_{\stackrel{\scr p|n}{ \scr p\neq j-1,\ p\neq j}}\left(1+\frac{1}{p-j}\right)=
\prod_{\stackrel{\scr p|n}{\scr p\neq j-1,\  p\neq j}}\left(1+\frac{1}{\phi_j(p)}\right)=\sum_{\stackrel{\scr d|n}{\scr (d,p(j-1)p(j))=1}}\frac{\mu^2(d)}{\phi_j(d)}.\label{2.30}
    \end{equation}
We see that for $n$ squarefree $H_0(n)=\sigma(n)/n$,  $H_1(n) = n/\phi(n)$, and in general for $j\ge 1$
    \begin{equation}
H_j(n) = \prod_{\stackrel{\scr p|n}{\scr p\neq j-1,\ p\neq j}}\left(\frac{p-j+1}{p-j}\right) =
       \frac{\phi_{j-1}(\frac{n}{(n,p(j-1)p(j))})}{ \phi_j(\frac{n}{(n,p(j-1)p(j))})},  \qquad (\mu^2(n)\neq 0). \label{2.40}
    \end{equation}
Next, we define the singular series for $j\ge 1$ and $n\neq 0$ by
    \begin{equation}
    \gs_j(n) =
        \left\{
        \begin{array}{ll}
        C_jG_j(n)H_j(n), &\mbox{if $p(j) |n$,}\\
                  0, & \mbox{otherwise}.
        \end{array}
        \right. \label{2.50}
    \end{equation}
where
\begin{equation} G_j(n) = \prod_{\stackrel{\scr p|n}{ \scr  p= j-1\ {\rm or}\ p=j}}\left(\frac{p}{p-1}\right), \label{2.60}
    \end{equation}
and
\begin{equation} C_j =  \prod_{\stackrel{\scr p}{ \scr  p\neq j-1,\ p\neq j}}\left( 1 - \frac{j-1}{(p-1)(p-j+1)}\right).\label{2.70}
    \end{equation}

The case of $j=3$ is special because it is the only case where $p(j-1)$ and $p(j)$ are both greater than one. We see that for $j= 1$ and $n\neq 0$
    \begin{equation}
\gs_1(n)=\prod_{p|n}\left(\frac{p}{p-1}\right)= \frac{n}{ \phi(n)}\label{2.80}
    \end{equation}
and for $j=2$ we have the familiar singular series for the Goldbach and prime twins conjectures
    \begin{equation}
 \gs_2(n) = \left\{ \begin{array}{ll}
      {\displaystyle 2C_2\prod_{\stackrel{\scr p \vert n}{\scr p>2}}\left(\frac{p-1}{p-2}\right),} &
        \mbox{if $n$ is  even,
    $n\neq 0$;}\\
0,   &\mbox{if $n$ is  odd;}\\
\end{array}
\right.\label{2.90}
\end{equation}
where
    \begin{equation}
 C_2 = \prod_{p>2}\left( 1 - \frac{1}{(p-1)^2}\right). \label{2.100}
    \end{equation}

    \begin{lemma}
For $R\ge 1$, $j \ge 0$, $p(j)|k$, and $0\le \log|k| \ll \log R$,   we have
    \begin{equation} \sum_{\stackrel{\scr d\le R} {\scr (d,k)=1}} \frac{\mu(d)}{\phi_j(d)} \log \frac{R}{d} = \gs_{j+1}(k) + r_j(R,k),\label{2.110}
    \end{equation}
where
    \begin{equation}
r_j(R,k) \ll_je^{-c_1\sqrt{\log R}},\label{2.120}
    \end{equation}
and  $c_1$ is an absolute positive constant.  Also,
\begin{equation} \sum_{\stackrel{\scr d\le R} {\scr (d,k)=1}} \frac{\mu(d)}{\phi_j(d)} \ll_j e^{-c_1\sqrt{\log R}}.\label{2.130}
    \end{equation}
\label{Lemma1}
\end{lemma}

Special cases of Lemma \ref{Lemma1} have been proved before. When $j=0$ this was used by Selberg \cite{S}, and also Graham \cite{Gr}, but we have made the error term stronger with regard to $k$ by an argument suggested in \cite{FG}. It is easy to make the $j$ dependence explicit in the error term, but in this paper we will assume $j$ is fixed (actually we only use $j\le 2$.) We will sometimes use Lemma \ref{Lemma1} in the weaker form
    \begin{equation*}
\sum_{\stackrel{\scr d\le R}{\scr (d,k)=1}} \frac{\mu(d)}{\phi_j(d)} \log \frac{R}{d} = \gs_{j+1}(k) + O_j(\frac{1}{(\log 2R)^A})
    \end{equation*}
for $A$ any positive number, and assuming the same conditions as in Lemma \ref{Lemma1}. Further, in handling error terms, we will need to remove the restriction $0
\le \log |k| \ll \log R$ in Lemma \ref{Lemma1}, in which case we have
the error estimate
\begin{equation} r_j(R,k) \ll_j m(k) e^{-c_1\sqrt{\log R}} ,\label{2.140} \end{equation}
which holds uniformly for $k\ge 1$ and $R\ge 1$, where $m(k)$ is defined below in equation \eqref{2.180}.  This estimate also holds in \eqref{2.130}.

The next lemma is a generalization of a result of Hildebrand \cite{Hi}.

\begin{lemma} For $R\ge 1$, $j\ge 1$ and $p(j)|k$, we have
    \begin{equation}
\sum_{\stackrel{\scr d\le R}{\scr (d,k)=1}} \frac{\mu^2(d)}{\phi_j(d)}
= \left\{ \begin{array}{ll}
        \frac{1}{\gs_j(k)}\left( \log R + D_j + h_j(k)\right) + O(\frac{m(k)}{\sqrt{R}} )
 & \mbox{if $p(j-1)|k$,} \\
       O(\frac{m(k)}{\sqrt{R}} )
,   &\mbox{if $p(j-1)\ndiv k$,}\\
\end{array}
\right.\label{2.150}
\end{equation}
where
\begin{equation}
 D_j = \gamma + \sum_{p\neq j-1}\frac{(2-j)\log p}{(p-j+1)(p-1)}, \label{2.160}
\end{equation}
    \begin{equation}\begin{split}
h_j(k) &= \sum_{p|k}\frac{\log p}{p-1} - \sum_{\stackrel{\scr p|k}{\scr p\neq j-1}}\frac{(2-j)\log p}{(p-j+1)(p-1)}\\
&= \sum_{\stackrel{\scr p|k}{\scr p\neq j-1}}\frac{\log p}{(p-j+1)}+ [(p(j-1),k)>1]\frac{\log(j-1)}{j-2}, \label{2.170}\end{split}
\end{equation}
and
\begin{equation} m(k) = \sum_{d|k}\frac{\mu^2(d)}{ \sqrt{d}}=\prod_{p|k}\left( 1 + \frac{1}{\sqrt{p}}\right).\label{2.180}
\end{equation}
\label{Lemma2}
\end{lemma}

The case $j=1$ of this lemma is Hilfssatz 2 of \cite{Hi}. The proof of this generalization only requires minor modifications in Hildebrand's proof which we sketch. In applying this lemma we will sometimes use the simple estimates (see \cite{GY})
\begin{equation} h_j(k) \ll_j \log\log 3k,  \qquad m(k) \ll \exp\left(\frac{c\sqrt{\log k}}{\log\log 3k}\right). \label{2.190}
\end{equation}
We will frequently use the estimate, for $p(j)|k$ and $\log |k| \ll \log R$,
\begin{equation}
\sum_{\substack{d\le R\\ (d,k)=1}}\frac{\mu^2(d)}{\phi_j(d)} \ll \log 2R \label{2.200}
\end{equation}
which follows immediately from Lemma \ref{Lemma2} or may be seen directly.
We also need the following result that is obtained by partial summation in Lemma \ref{Lemma2}.
For $j\ge 1$ and $p(j)|k$, we have
    \begin{equation}\begin{split}
\sum_{\stackrel{\scr d\le R}{\scr (d,k)=1}}& \frac{\mu^2(d)}{\phi_j(d)}\log \frac{R}{d}  \\
&= \left\{ \begin{array}{ll}
        \frac{1}{\gs_j(k)}\left( \frac{1}{2}\log^2 R + (D_j + h_j(k))\log R + E_j(k) \right) + O(\frac{m(k)}{\sqrt{R}} )
 & \mbox{if $p(j-1)|k$,} \\
      E_j(k)+ O(\frac{m(k)}{\sqrt{R}} )
,   &\mbox{if $p(j-1)\ndiv k$,}\\
\end{array}
\right. \label{2.210} \end{split}
\end{equation}
where $E_j(k)$ is given by
\begin{equation}
E_j(k) = \int_1^\infty \Big(\sum_{\stackrel{\scr d\le u}{\scr (d,k)=1}} \frac{\mu^2(d)}{\phi_j(d)} -
        \frac{1}{\gs_j(k)}\left( \log u + D_j + h_j(k)\right)\Big) \, \frac{du}{u} .\label{2.220}
\end{equation}

\begin{lemma} For $R\ge 1$, $j\ge 1$, $p(j)|k$, and $0\le \log|k|\ll \log R$, we have
\begin{eqnarray}\lefteqn{ \sum_{\stackrel{\scr d\le R}{\scr (d,k)=1}}\frac{\mu(d)}{\phi_{j}(d)}\gs_{j+1}(dk)\log \frac{R}{d}} \nonumber  \\  & & =  \mu(p(j+1))\mu((k,p(j+1)))\gs_{j+1}(kp(j+1))\Big(\gs_{j+2}(kp(j+1)) \label{2.230}
 \\ & & +  r_{j+1}(\frac{R(k,p(j+1))}{p(j+1)},kp(j+1))\Big), \nonumber \end{eqnarray}
where $r_j(R,k)$ is the error term in Lemma \ref{Lemma1}.
\label{Lemma3}
\end{lemma}
\bigskip
\begin{lemma}  For $R\ge 1$, $j\ge 1$ and $p(j)|k$, we have
\begin{equation} \begin{split} \sum_{\stackrel{\scr d\le R}{\scr (d,k)=1}}\frac{\mu^2(d)}{\phi_{j}(d)}\gs_{j+1}(dk)
&=  \log \left(\frac{R(k,p(j+1))}{p(j+1)}\right) + D_{j+1} + h_{j+1}(kp(j+1))\\& \qquad + O\left(\frac{\gs_{j+1}(kp(j+1))m(kp(j+1))}{\sqrt{\frac{R(k,p(j+1))}{p(j+1)}}}\right)
\label{2.240}\end{split}
\end{equation}
and
\begin{equation} \sum_{\stackrel{\scr d|r}{\scr (d,k)=1}}\frac{\mu^2(d)}{\phi_{j}(d)}\gs_{j+1}(dk)
= \gs_{j+1}(rk). \label{2.250}
\end{equation}
\label{Lemma4}
\end{lemma}
\bigskip
Our final lemma relates the singular series given in (\ref{1.90}) for $r$ equal to two and three to the singular series in (\ref{2.50}).
\begin{lemma} \label{Lemma5} For $\mbox{\boldmath{$k$}}=(0,k)$, with $k\neq 0$, we have
\begin{equation}
\gs(\mbox{\boldmath{$k$}}) = \gs_2(k),
\label{2.260}
\end{equation}
and for $\mbox{\boldmath{$k$}}=(0,k_1,k_2)$, $k_1\neq k_2 \neq 0$, $\kappa = (k_1,k_2)$, and $\Delta = k_1k_2(k_2-k_1)$, we have
\begin{equation} \gs(\mbox{\boldmath{$k$}}) = \gs_2(\kappa)\gs_3(\Delta).
\label{2.270}
\end{equation}
\end{lemma}

\section{ Proof of the Lemmas}

\emph{ Proof of Lemma \ref{Lemma1}.} We assume $p(j)|k$. Let $s=\sigma +it$, $\sigma >0$, and define
\begin{eqnarray} F(s) &=& \sum_{\stackrel{\scr n=1}{ \scr (n,k)=1}}^\infty \frac{\mu(n)} {\phi_j(n)n^s} = \prod_{p\ndiv k}
\left(1- \frac{1}{(p-j)p^s}\right) \nonumber \\
&=&\frac{1}{ \zeta(s+1)} \prod_{p|k}\left(1-\frac{1}{ p^{s+1}}\right)^{-1}\prod_{p \ndiv k}
\left(1- \frac{1}{(p-j)p^s}\right) \left(1-\frac{1}{ p^{s+1}}\right)^{-1}\nonumber \\
&=& \frac{1}{\zeta(s+1)}\prod_{p|k}\left(1-\frac{1}{ p^{s+1}}\right)^{-1}  \prod_{   p \ndiv k}\Big( 1 -\frac{j}{ (p-j)(p^{s+1}-1)}\Big) \nonumber \\
& =&\frac{1}{\zeta(s+1)}g_k(s)h_k(s).\label{3.1}
\end{eqnarray}
We see the product for $h_k(s)$ converges absolutely for $\mathrm{Re}(s) > -1$, and therefore $F(s)$ is an analytic function in this half-plane except possibly for poles at the zeros of $\zeta(s+1)$.

We now apply the formula, for $m \ge 2$ and $b>0$,
\begin{equation*} \frac{(m-1)!}{ 2\pi i}\int_{b-i\infty}^{b+i\infty} \frac{x^s}{ s^m}\, ds =
\left\{ \begin{array}{ll}
       0,
        &\mbox{if $0<x\le 1$,
    }  \\
(\log x)^{m-1},& \mbox{if $x\ge 1$,}\nonumber
\end{array}
\right.
\end{equation*}
which, in the case $m=2$, gives
\begin{equation}\sum_{\stackrel{\scr d\le R }{ \scr (d,k)=1}}\frac {\mu(d)}{ \phi_j(d)} \log \frac{R}{ d} =
\frac {1}{ 2\pi i}\int_{b-i\infty}^{b+i\infty} F(s)\frac{R^s}{ s^2}\, ds .\label{3.20} \end{equation}

By Theorem 3.8 and (3.11.8) of \cite{T}  there exists a small positive constant $c$ such that $\zeta(\sigma +it)\neq 0$ in the region  $\sigma \ge 1 -\frac{c}{ \log(|t|+2)} $ and all $t$, and further
\begin{equation}
\frac{1}{ \zeta(\sigma +it)} \ll \log(|t|+2)
\label{3.30}
\end{equation}
  in this region. (There are stronger results but this suffices for our needs.)
We now move the contour to the left  to the path $\mathcal{ L}$ given by $s= -\frac{c}{ \log(|t|+2)} +it$.  When $p(j+1)\ndiv k$, $h_k(s)$ has a simple zero at $s=0$ and hence $F(s)/s^2$ is analytic at $s=0$ and  no contribution occurs, but when $p(j+1)|k$ (including $p(j+1)=1$), $F(s)/s^2$ has a simple pole at $s=0$ which gives a contribution from the residue of
\begin{eqnarray}  g_k(0)h_k(0)&=& \prod_{ p | k}\left(1-\frac{1}{ p}\right)^{-1}\prod_{p \ndiv k }\left(1-\frac{j}{(p-1)(p-j) }\right)\nonumber  \\
&=&\prod_{\stackrel{\scr p|k}{ \scr  p= j \ {\rm or}\  p=j+1}}\left(\frac{p}{ p-1}\right)\prod_{\stackrel{\scr p | k}{ \scr p\neq j ,\ p\neq j+1}} \left(1-\frac{1}{ p}\right)^{-1}\left(1-\frac{j}{(p-1)(p-j) }\right)^{-1}\nonumber \\  & & \qquad \prod_{\stackrel{\scr p}{ \scr  p \neq j , \ p\neq j+1 }}\left(1-\frac{j}{(p-1)(p-j) }\right)\nonumber  \\
&=& G_{j+1}(k)\prod_{\stackrel{\scr p}{ \scr p\neq j,  p \neq j+1}}\left(1-\frac{j}{(p-1)(p-j) }\right) \prod_{\stackrel{\scr p | k }{ \scr p\neq j, \ p\neq j+1}}\left(\frac{p-j}{ p-j-1}\right)\nonumber  \\
& =& C_{j+1}G_{j+1}(k)H_{j+1}(k).\nonumber
\end{eqnarray}
Hence we have
\begin{equation}
\sum_{\stackrel{\scr d\le R }{ \scr (d,k)=1}}\frac{\mu(d)}{ \phi_j(d)} \log \frac{R}{ d} = \gs_{j+1}(k) + \frac{1}{ 2\pi i}\int_{\mathcal{ L}} F(s)\frac{R^s}{ s^2}\, ds =\gs_{j+1}(k) +  r_j(R,k).\label{3.40}
\end{equation}
It remains to estimate the integral in \eqref{3.40}.
On $\mathcal{ L}$ we have $-\frac{1}{ 4} \le \sigma <0$, and therefore
\[|h_k(s)| \ll \prod_{p}\big(1 + O_j(\frac{1}{ p^{2 +\sigma}})\Big) \ll_j 1.\]
For $g_k(s)$, we see that $g_1(s)=1$ since the product defining $g_1(s)$ is void, and for $k>1$
\begin{eqnarray*}
 \log |g_k(s)| &\le &   - \sum_{p|k}\log(1-\frac{1}{ p^{\sigma +1}})  \\
&\ll &\sum_{p|k} \frac{1}{ p^{\sigma +1}}\\ &\ll& \sum_{p<2\log 2k}\frac{1}{ p^{\sigma +1}}\ll (\log k)^{-\sigma}   \end{eqnarray*}
and hence
\begin{equation} |g_k(s)| \ll  e^{\sqrt[4]{\log k}}, \quad \mathrm{for}  \ k\ge 1 .\label{3.50}
\end{equation}
Thus  the integral in \eqref{3.40} is
\[ \ll_j  e^{\sqrt[4]{\log k}}\int_{-\infty}^\infty R^{-\frac{c}{ \log(|t|+2)}}\log(|t|+2)\frac{1}{\big(\frac{c}{ \log(|t|+2))}\big)^2 + t^2}\, dt \]
This last integral is, for any $w\ge 2$,
\[ \ll  \int_0^w R^{-\frac{c}{\log(|t|+2)}}\, dt + \int_w^\infty \frac{\log t}{ t^2}\, dt
 \ll w e^{\frac{-c\log R}{ \log w}} + \frac{\log w }{ w} \]
and hence, on choosing $\log w = \frac{1}{ 2}\sqrt{c\log R}$ we have, since $\log |k| \ll \log R$, that the error term is
\[\ll_j  e^{\root 4 \of {\log k}} \left( e^{-\sqrt{c\log R}} + \sqrt{c\log R}e^{-\frac{1}{ 2} \sqrt{c\log R}}\right) \ll_je^{-c_1\sqrt{\log R}},\]
which proves the first part of Lemma 1.

The bound in equation \eqref{2.140} follows from the previous argument when we replace the estimate for $g_k(s)$ used above by the bound, for $-\frac{1}{4}<\sigma <0$,
\[ |g_k(s)| \ll \prod_{p|k}\left(1+\frac{1}{p^{1+\sigma}}\right) \ll m(k), \]
which follows from
\begin{eqnarray*}
 \log |g_k(s)| &\le &   - \sum_{p|k}\log(1-\frac{1}{ p^{\sigma +1}})  \\
&= &\sum_{p|k}\left( \log(1 + \frac{1}{p^{\sigma +1}})+ O(\frac{1}{p^{2(\sigma +1)}})\right)\\
&=& \sum_{p|k} \log(1 + \frac{1}{p^{\sigma +1}}) + O(\sum_p \frac{1}{p^{\frac{3}{2}}})\\
&=& \sum_{p|k} \log(1 + \frac{1}{p^{\sigma +1}}) +O(1). \end{eqnarray*}

To prove \eqref{2.130}, we apply Perron's formula (see \cite{T}, Chapter 3) in the usual way to obtain, with $b= 1/\log R$,
\[ \sum_{\stackrel{\scr d\le R} {\scr (d,k)=1}} \frac{\mu(d)}{\phi_j(d)}
= \frac{1}{2\pi i}\int_{b - iT}^{b+iT} F(s) \frac{R^s}{s}\, ds + O_j(\frac {\log^2R}{T}). \]
Moving the contour to the left to $\mathcal{L}$ for $-T\le t\le T$ we have no residue, and we may estimate  the integral along $\mathcal{L}$ and the upper and lower horizontal paths by
\[ \ll_j  e^{\sqrt[4]{\log k}} \Big( \log^2T e^{-c \frac{\log R}{\log T}} + \frac{\log T}{T} \Big).\]
Now $\log|k| \ll \log R$, and taking $T= e^{\sqrt{\log R}}$  shows the above error is $\ll e^{-c_1\sqrt{ \log R}} $ for a sufficiently small constant $c_1$. This completes the proof of Lemma \ref{Lemma1}.

\emph{ Proof of Lemma \ref{Lemma2}.}  We follow Hildebrand's proof of the case $j=1$, indicating only the main steps. We assume $p(j)|k$, and $j\ge 1$. Letting

\begin{equation} f_k(n) = \left\{ \begin{array}{ll}
       \frac{\mu^2(n) n}{ \phi_j(n)} , &
        \mbox{if $ (n,k)=1$,
    }  \\
0,   & \mbox{otherwise, }
\end{array} \right. \label{3.60} \end{equation}
and defining $g_k(n)$ by
\[ f_k(n) = \sum_{d|n}g_k(d) \]
so that
\begin{equation} g_k(p^m) = f_k(p^m)-f_k(p^{m-1}) = \left\{ \begin{array}{ll}
       \frac{j}{ p-j}, &
        \mbox{$ m=1, p\not | k $,
    }  \\
-1,   &\mbox{$m=1, p|k $},  \\
\frac{-p}{ p-j}  & \mbox{ $ m=2, p\not | k $}, \\
0, & \mbox{ $m=2$ and $p|k$, or $m>2$,}
\end{array} \right.\label{3.70}
\end{equation}
then we have
\begin{eqnarray} \sum_{\stackrel{\scr d\le R}{ \scr (d,k)=1}} \frac{\mu^2(d)}{ \phi_j(d)}&=&
\sum_{n\le R} \frac{f_k(n)}{ n} \nonumber \\ &=& \sum_{n \le R} \frac{1}{ n}\sum_{d|n}g_k(d) \nonumber  \\ & =& \sum_{d\le R}\frac{g_k(d)}{ d}\sum_{m \le R/d} \frac{1}{ m}\nonumber \\
&=& \sum_{d\le R}\frac{g_k(d)}{ d}\left( \log\frac{R}{ d} + \gamma + O(\frac{d}{ R})\right) \nonumber  \\
&=& \sum_{d\le R}\frac{g_k(d)}{ d}\log \frac{R}{ d} + \gamma \sum_{d\le R}\frac{g_k(d)}{ d} + O( \frac{1}{ R}\sum_{d\le R}|g_k(d)|)\nonumber  \\
&=& S_1 + \gamma S_2 +O(\frac{1}{ R}S_3).
  \label{3.80} \end{eqnarray}
Using  \eqref{3.70} and the multiplicativity of $g_k(n)$ we easily verify that
\[ M_j(k) = \sum_{d=1}^\infty \frac{g_k(d)}{ d} = \left\{ \begin{array}{ll}
       \frac{1 }{ \gs_j(k)}
 , &
        \mbox{if $p(j-1)|k$,
    }  \\
0,   &\mbox{if $p(j-1)\not | k$. }
\end{array} \right. \]
 As in \cite{Hi}  we find
\[ S_3 \ll_j \sqrt{R}m(k)\]
which  by partial summation implies
\[\sum_{d >R } \frac{g_k(d)}{ d} \ll\frac {m(k)}{ \sqrt{R}}.\]
These results  now imply
\[ S_2 = M_j(k) + O(\frac{m(k)}{ \sqrt{R}})\]
and
\[ S_1 = M_j(k)\log R  -\sum_{d=1}^\infty \frac{g_k(d)\log d}{ d} + O(\frac{m(k)}{ \sqrt{R}}).\]
Finally, letting
\[G_k(s) = \sum_{n=1}^\infty \frac{g_k(n)}{ n^{s+1}},\]
we have
\[- \sum_{d=1}^\infty \frac{g_k(d)\log d}{ d} = {G_k}'(0)\]
which on using \eqref{2.160} and the Euler product for $G_k(s)$ gives
${G_k}'(0) = M_j(k)(D_j + h_j(k)-\gamma)$
where $h_j(k)$ is given by \eqref{2.170}.

\emph{ Proof of Lemma \ref{Lemma3}.} We have, assuming $p(j)|k$,
\begin{eqnarray*} \lefteqn{\sum_{\stackrel{\scr d\le R }{ \scr (d,k)=1}} \frac{\mu(d)}{ \phi_{j}(d)}\gs_{j+1}(dk)\log \frac{R}{ d} }\\ &= &C_{j+1}G_{j+1}(k)H_{j+1}(k) \sum_{\stackrel{\stackrel{\scr d\le R }{ \scr (d,k)=1}}{ \scr p(j+1)|dk}}\frac{\mu(d)}{ \phi_{j}(d)}G_{j+1}(d)H_{j+1}(d)\log \frac{R}{ d}  \\
& = & C_{j+1}G_{j+1}(k)H_{j+1}(k) \sum_{\stackrel{\scr d\le R }{\stackrel{\scr (d,k)=1 }{ \scr \frac{p(j+1)}{ (k, p(j+1))}|d}}}\frac{\mu(d)}{ \phi_j(d)}G_{j+1}(d)H_{j+1}(d)\log \frac{R}{ d}  \\
& =&\mu(p(j+1))\mu((k,p(j+1)))C_{j+1}G_{j+1}(kp(j+1))H_{j+1}(k)\\
& & \qquad \sum_{\stackrel{\scr m\le \frac{R(k,p(j+1))}{ p(j+1)} }{ \scr (m,kp(j+1))=1}}\frac{\mu(m)}{ \phi_{j+1}(m)}\log \Big(\frac{R(k,p(j+1))}{ mp(j+1)}\Big) \\
&=&\mu(p(j+1))\mu((k,p(j+1)))\gs_{j+1}(kp(j+1))\\
& & \qquad \Big(\gs_{j+2}(kp(j+1)) + r_{j+1}\big(\frac{R(k,p(j+1))}{ p(j+1)},kp(j+1)\big)\Big),
\end{eqnarray*}
by Lemma \ref{Lemma1}.

\emph{ Proof of Lemma \ref{Lemma4}.} The same argument used above to prove Lemma \ref{Lemma3} shows that the sum in (\ref{2.240}) is equal to
\[ \gs_{j+1}(kp(j+1))\sum_{\stackrel{\scr m\le \frac{R(k, p(j+1))}{ p(j+1)} }{ \scr (m,kp(j+1))=1}}\frac{\mu^2(m)}{ \phi_{j+1}(m)}.\]
Since $p(j)|k$, equation (\ref{2.240}) now follows from Lemma \ref{Lemma2}.

To prove \eqref{2.250}, we proceed as before in the proof of Lemma \ref{Lemma3} to obtain
\begin{equation} \sum_{\stackrel{\scr d|r }{ \scr (d,k)=1}}\frac{\mu^2(d)}{ \phi_{j}(d)}\gs_{j+1}(dk) = C_{j+1}G_{j+1}(k)H_{j+1}(k) \sum_{\stackrel{\stackrel{\scr d|r }{ \scr (d,k)=1}}{ \scr p(j+1)|dk}}\frac{\mu^2(d)}{ \phi_{j}(d)}G_{j+1}(d)H_{j+1}(d).\label{3.90}
\end{equation}
Suppose first that $p(j+1)|k$. Then since $p(j)|k$ also we have that $G_{j+1}(d)=1$ in the sum on the right, and hence our expression becomes
\begin{eqnarray} &=& C_{j+1}G_{j+1}(k)H_{j+1}(k) \sum_{\stackrel{\scr d|r }{ \scr (d,k)=1}}\frac{\mu^2(d)}{ \phi_{j}(d)}H_{j+1}(d)\nonumber  \\
&=& C_{j+1}G_{j+1}(k)H_{j+1}(k) \sum_{\stackrel{\scr d|r }{ \scr (d,k)=1}}\frac{\mu^2(d)}{ \phi_{j+1}(d)}\nonumber  \\
 &=& C_{j+1}G_{j+1}(k)H_{j+1}(k) \prod_{\stackrel{\scr p|r}{ \scr p\ndiv k}}\left( 1 + \frac{1}{ p-j-1}\right) \nonumber  \\
&=& C_{j+1}G_{j+1}(k)H_{j+1}(rk) = \gs_{j+1}(rk) \nonumber
\end{eqnarray}
since here $G_{j+1}(k) = G_{j+1}(p(j)p(j+1))= G_{j+1}(rk)$.

Now assume the alternative case that $p(j+1)\ndiv k$. Then the right-hand side of \eqref{3.90} is
\[=C_{j+1}G_{j+1}(k)H_{j+1}(k) \sum_{\stackrel{\scr p(j+1)e|r }{ \scr (e,kp(j+1))=1}}\frac{\mu
^2(e)}{ \phi_{j}(e)}G_{j+1}(ep(j+1))H_{j+1}(e).\]
If $p(j+1)\ndiv r$ this sum has no terms and is zero which proves \eqref{2.250} in this case. If $p(j+1)|r$
our expression becomes
\begin{eqnarray} &=& C_{j+1}G_{j+1}(p(j)p(j+1))H_{j+1}(k) \sum_{\stackrel{\scr e|r }{ \scr (e,kp(j+1))=1}}\frac{\mu^2(e)}{ \phi_{j+1}(e)}\nonumber  \\
&=& C_{j+1}G_{j+1}(p(j)p(j+1))H_{j+1}(k)\prod_{\stackrel{\scr p| r }{ \scr p\ndiv kp(j+1)}}\left( 1 + \frac{1}{ p-j-1}\right) \nonumber  \\
&=& \gs_{j+1}(rk), \nonumber
\end{eqnarray}
which completes the proof of Lemma \ref{Lemma4}.
\bigskip

\emph{ Proof of Lemma \ref{Lemma5}.} For the case $\mbox{\boldmath{$k$}}=(0,k)$ we have by (\ref{1.90}) that
\[ \gs(\mbox{\boldmath{$k$}}) =4\left(1 - \frac{\nu_2(\mbox{\boldmath{$k$}})}{2}\right)\prod_{p>2} \left( 1 -\frac{1}{p}\right)^{-2}\left(1-\frac {\nu_p(\mbox{\boldmath{$k$}})}{p}\right). \]
Now $\nu_p(\mbox{\boldmath{$k$}}) = 1$ if $p|k$ and $\nu_p(\mbox{\boldmath{$k$}})=2$ if $p\ndiv k$, and hence
\begin{eqnarray*}
\gs(\mbox{\boldmath{$k$}}) &=& 2[2|k]\prod_{\substack{p|k\\ p>2}}  \left( 1 -\frac{1}{p}\right)^{-1}\prod_{\substack{p\ndiv k\\ p>2}}\left(1-\frac {1}{p}\right)^{-2}\left(1 - \frac{2}{p}\right) \\
&=& 2[2|k]\prod_{\substack{p|k\\ p>2}}  \left( 1 -\frac{1}{p}\right)\left(1-\frac{2}{p}\right)^{-1}\prod_{p>2}\left(1-\frac {1}{p}\right)^{-2}\left(1 - \frac{2}{p}\right) \\
&=& 2[2|k]\prod_{\substack{p|k\\ p>2}}  \left( \frac{p-1}{p-2}\right)\prod_{p>2}\left(1-\frac{1}{(p-1)^2}\right)\\
&=& \gs_2(k).
\end{eqnarray*}
Next, in general,
\begin{eqnarray*}
\gs(\mbox{\boldmath{$k$}}) &=& \prod_p \left( 1 -\frac{1}{p}\right)^{-r}\left(1-\frac {\nu_p(\mbox{\boldmath{$k$}})}{p}\right) \\
&=& \left( 1 -\frac{1}{2}\right)^{-r}\left( 1 -\frac{1}{3}\right)^{-r}\left(1-\frac {\nu_2(\mbox{\boldmath{$k$}})}{2}\right)\left(1-\frac {\nu_3(\mbox{\boldmath{$k$}})}{3}\right)\prod_{p>3}\left(\frac{p}{p-1}\right)^{r-1}\left(\frac{p-\nu_p(\mbox{\boldmath{$k$}})}{p-1}\right).
\end{eqnarray*}
With $r=3$, $\mbox{\boldmath{$k$}}=(0,k_1,k_2)$, $\kappa = (k_1,k_2)$, and $\Delta =k_1k_2(k_2-k_1)$ we have
    \begin{equation*}
    \nu_2(\mbox{\boldmath{$k$}}) =
        \left\{
        \begin{array}{ll}
        1, &\mbox{if $2|\kappa$ ,}\\
                  2, & \mbox{if $2\not | \kappa$},
        \end{array}
        \right.
    \end{equation*}
and for $p\ge 3$
    \begin{equation*}
    \nu_p(\mbox{\boldmath{$k$}}) =
        \left\{
        \begin{array}{ll}
        1, &\mbox{if $p|\kappa$,}\\
                  2, & \mbox{if $p|\Delta $, $p \not | \kappa$,} \\
              3, &\mbox{if $p\not | \Delta$}.
        \end{array}
        \right.
    \end{equation*}
We now see that $\gs(\mbox{\boldmath{$k$}})=0$ if $2\ndiv \kappa$ or $3\ndiv\Delta $, which by (\ref{2.50}) proves the lemma in these cases. Thus we now assume that $2|\kappa$ and $3|\Delta$, and  have
\begin{eqnarray*} \gs(\mbox{\boldmath{$k$}})&= &\frac{9}{\nu_3(\mbox{\boldmath{$k$}})}\prod_{p>3}\left(\frac{p}{p-1}\right)^2\left(\frac{p-\nu_p(\mbox{\boldmath{$k$}})}{p-1}\right) \\
&=& \frac{9}{\nu_3(\mbox{\boldmath{$k$}})}\prod_{p>3}\left(\frac{p}{p-1}\right)^2\left(\frac{p-3}{p-1}\right)\prod_{\stackrel{\scr p|\kappa}{\scr p>3}}\left(\frac{p-1}{p-3}\right)\prod_{\stackrel{\stackrel{\scr p|\Delta}{\scr p\ndiv \kappa}}{\scr p>3}}\left(\frac{p-2}{p-3}\right)\\
&=& \frac{9}{\nu_3(\mbox{\boldmath{$k$}})}\prod_{p>3}\left(1-\frac{1}{(p-1)^2}\right)\left(1- \frac{2}{(p-1)(p-2)}\right)\prod_{\stackrel{\scr p|\kappa}{\scr p>3}}\left(\frac{p-1}{p-2}\right)\prod_{\stackrel{\scr p|\Delta}{\scr p>3}}\left(\frac{p-2}{p-3}\right)\\
&=& 6C_2C_3H_2(\kappa)H_3(\Delta)\\
&=& \gs_2(\kappa)\gs_3(\Delta) ,
\end{eqnarray*}
which proves Lemma \ref{Lemma5}.

\section{Proof of Theorem \ref{Theorem2}}

We now prove that the mixed correlations can be reduced to the pure correlations through an application of the Bombieri-Vinogradov theorem. For $k\ge 2$
we consider the general sum, for $R\le P \le N$,
\begin{equation}\mathcal{S}_P(N,\mbox{\boldmath$j$}) =\sum_{n=1}^N \Lambda_R (n+j_1)\Lambda_R (n+j_2)\cdots \Lambda_R(n+j_{k-1})\Lambda_P(n) \label{4.10}\end{equation}
where the $j_i$'s are not necessarily distinct, but satisfy
\begin{equation} j_i\neq 0 , \quad  1\le i\le k-1. \label{4.20}
\end{equation}
Note that for $1<n\le N$, $\Lambda_N(n) = \Lambda(n)$. We  prove Theorem \ref{Theorem2} by proving the following theorem which shows that both $\mathcal{S}_R(N,\mbox{\boldmath$j$})$
and $\mathcal{S}_N(N,\mbox{\boldmath$j$})$ in certain ranges are asymptotic  to the same main term.  Let
\begin{equation} \label{4.30}
W_R(\mbox{\boldmath$j$})= \sum_{\substack{d_1,d_2,\ldots, d_{k-1} \le R\\ (d_i,j_i)=1,\, 1\le i\le k-1\\(d_r,d_s)| j_s-j_r, \, 1\le r<s\le k-1 }}\frac{1}{\phi(D_{k-1})}\prod_{i=1}^{k-1}\mu(d_i)\log\frac{R}{d_i},
\end{equation}
where
\begin{equation}
D_{k-1}=  [d_1,d_2,
\ldots , d_{k-1}]. \label{4.40}
\end{equation}

\begin{theorem}\label{Theorem4.1} We have, for $k\ge 2$, $N^\epsilon \le R \le N^{\frac{1}{k}}$,  and $\max_i|j_i|\ll R^{\frac{1}{k}}$,
\begin{equation} \mathcal{S}_R(N,\mbox{\boldmath$j$}) = N W_R(\mbox{\boldmath$j$}) + O(R^{k}) +O_k\Big(N(\log^{4^{k-1}+ k-2}N )e^{-c_1\sqrt{\log\frac{R}{\max|j_i|^{k-1}}}}\Big),\label{4.50}
\end{equation}
and, for $N^\epsilon \le R \le N^{\frac{\vartheta}{k-1} -\epsilon}$,
\begin{equation}\mathcal{S}_N(N,\mbox{\boldmath$j$}) = N W_R(\mbox{\boldmath$j$})  +O_k\big(\frac{N}{(\log N)^{\frac{A}{2} - 4^{k-3/2}+2^{k-1}-1}}\big). \label{4.60}
\end{equation}
\end{theorem}

We have
\begin{equation}
\mathcal{S}_P(N,\mbox{\boldmath$j$}) = \sum_{d_1,d_2,\ldots, d_{k-1} \le R}\Big(\prod_{i=1}^{k-1}\mu(d_i)\log\frac{R}{d_i}\Big)\Big(\sum_{\substack{n=1\\ d_i|n+j_i,\ 1\le i\le k-1}}^N\Lambda_P(n)\Big). \label{4.70}
\end{equation}
Let
\begin{equation}
T_P(N, \mbox{\boldmath$j$})= \sum_{\substack{n=1\\ d_i|n+j_i,\  1\le i\le k-1}}^N\Lambda_P(n). \label{4.80}
\end{equation}
The $k-1$ congruence relations $n\equiv -j_i (\mathrm{mod}\, d_i)$ will have no solutions unless $(d_r,d_s)| j_s-j_r$ for all $1\le r<s\le k-1$.
If these divisibility conditions hold, then by the Chinese remainder theorem there exists a unique solution to these congruences $n\equiv a (\mathrm{mod}\, [d_1,d_2,
\ldots , d_{k-1}])$ for some $a= a(\mbox{\boldmath{$d$}},\mbox{\boldmath{$j$}})$. Here $a$ satisfies the original congruences $a\equiv -j_i (\mathrm{mod}\, d_i)$ for $1\le i\le k-1$. Thus we have
\begin{eqnarray}
T_P(N, \mbox{\boldmath{$j$}})&=& [ (d_r,d_s)| j_s-j_r, \ 1\le r<s\le k-1 ]\Big(  \sum_{\substack{1\le n\le N\\ n\equiv a (\mathrm{mod}\, D_{k-1})}}\Lambda_P(n)\Big)\nonumber\\
&=& [ (d_r,d_s)| j_s-j_r, \ 1\le r<s\le k-1 ]\psi_P(N;D_{k-1},a). \label{4.90}
\end{eqnarray}
We next have that
\begin{equation*}
\psi_P(N;D_{k-1},a) = \sum_{d\le P} \mu(d)\log\frac{P}{d}\sum_{\substack{1\le n\le N\\ n\equiv a (\mathrm{mod}\, D_{k-1})\\ n\equiv 0(\mathrm{mod}\, d)}}1.
\end{equation*}
The two congruences in this sum are solvable provided $(d,D_{k-1})|a$, in which case we have that $n$ runs through a residue class modulo $[D_{k-1},d]$.
Hence
\begin{equation*}
\psi_P(N;D_{k-1},a) =N \sum_{\substack{d\le P\\ (d,D_{k-1})|a}} \frac{\mu(d)}{[D_{k-1},d]}\log\frac{P}{d} + O(P).
\end{equation*}
We now write $d= gd'$, where $g=(d,D_{k-1})$. Hence $[d,D_{k-1}] = D_{k-1}d'$, and by Lemma \ref{Lemma1}
\begin{eqnarray*}
\psi_P(N;D_{k-1},a) &=& \frac{N}{D_{k-1}}\sum_{\substack{g|D_{k-1},\ g|a \\ g\le P}}\mu(g)\sum_{\substack{d'\le \frac{P}{g}\\(d',D_{k-1})=1 }} \frac{\mu(d')}{d'}\log\frac{P}{gd'} + O(P)\\
&=&  \frac{N}{\phi(D_{k-1})}\sum_{\substack{g|D_{k-1},\ g|a \\ g\le P}}\mu(g)\\ & & \qquad + O\big(\frac{N}{D_{k-1}}\sum_{\substack{g|D_{k-1},\ g|a \\ g\le P}}\mu^2(g)e^{-c_1\sqrt{\log(P/g)}}\big) +O(P).
\end{eqnarray*}
Since $a\equiv -j_i (\mathrm{mod}\, d_i)$ for $1\le i\le k-1$ and $g|a$ we see that $(g,d_i)|j_i$ for $1\le i\le k-1$, and conversely since $g|D_{k-1}$ these divisibility conditions imply that $g|a$. Hence
\begin{eqnarray}
\psi_P(N;D_{k-1},a) &=&  \frac{N}{\phi(D_{k-1})}\sum_{\substack{g|D_{k-1}\\ g\le P\\ (g,d_i)|j_i, \, 1\le i\le k-1}}\mu(g)\nonumber \\ & & + O\big(\frac{N}{D_{k-1}}\sum_{\substack{g|D_{k-1} \\ g\le P\\ (g,d_i)|j_i, \, 1\le i\le k-1}}\mu^2(g)e^{-c_1\sqrt{\log(P/g)}}\big) +O(P). \label{4.100}
\end{eqnarray}
The truncated M\"obius function sum complicates the calculations of our pure correlations when one or more of the $j_i=0$, but when all the $j_i\neq 0$ the truncation problem disappears. Thus, we see in this sum that $g\le \prod_{i=1}^{k-1}| j_i|  \le \max|j_i|^{k-1}$, and hence provided
\begin{equation}
\max_{1\le i\le k-1}|j_i| \le P^{\frac{1}{k-1}} \label{4.110}
\end{equation}
we have
\begin{equation*}
\sum_{\substack{g|D_{k-1}\\ g\le P\\ (g,d_i)|j_i, \, 1\le i\le k-1}}\mu(g) = \sum_{\substack{g|D_{k-1}\\ (g,d_i)|j_i, \, 1\le i\le k-1}}\mu(g)=
        \left\{
        \begin{array}{ll}
        1, &\mbox{if $(d_i,j_i)=1$,\ $1\le i \le k-1$},\\
            0, &\mbox{otherwise}.
        \end{array}
        \right.
    \end{equation*}
We conclude that subject to \eqref{4.110},
\begin{equation}\begin{split}
\psi_P(N;D_{k-1},a) &=[(d_i,j_i)=1,\ 1\le i \le k-1]  \frac{N}{\phi(D_{k-1})} \\ &\qquad + O\big(\frac{N}{D_{k-1}}d(D_{k-1})e^{-c_1\sqrt{\log(P/\max|j_i|^{k-1})}}\big) +O(P). \label{4.120}
\end{split}\end{equation}
Hence by \eqref{4.30}, \eqref{4.70}, \eqref{4.80}, \eqref{4.90}, and \eqref{4.120}, we obtain
\begin{equation}\begin{split}
\mathcal{S}_P(N,\mbox{\boldmath$j$}) =& NW_R(\mbox{\boldmath$j$}) +O(PR^{k-1}) \\ & + O\Big(N\log^{k-1}N e^{-c_1\sqrt{\log\frac{P}{\max|j_i|^{k-1}}}}\sum_{d_1,d_2,\ldots , d_{k-1}\le R}\frac{d(D_{k-1})}{D_{k-1}}\Big). \label{4.130}
\end{split}
\end{equation}
We now estimate the sum in the second error term sufficiently well for our needs. Letting
\[ D_R(n)= \sum_{\substack{e|n \\ e\le R}}d(e),\]
then
\begin{eqnarray*}
\sum_{n\le N}D_R(n)^m &=& \sum_{n\le N} \sum_{\substack{e_1,e_2,\ldots , e_m \le R\\ e_i|n,\, 1\le i\le m }}d(e_1)d(e_2)\ldots d(e_m)\\
&=& \sum_{e_1,e_2,\ldots , e_m\le R}d(e_1)d(e_2)\ldots d(e_m)\sum_{\substack{n\le N\\ e_i|n, 1\le i\le m}}1\\
&=& N \sum_{e_1,e_2,\ldots , e_m\le R} \frac{d(e_1)d(e_2)\ldots d(e_m)}{[e_1,e_2,\ldots , e_m]}
+ O\Big(\big(\sum_{e\le R} d(e)\big)^m\Big).
\end{eqnarray*}
The last error term is $O(R^m\log^mR)$,
and $D_R(n)\le \sum_{e|n}d(n) = d^2(n)$. Hence, using the estimate \cite{Hua}
\begin{equation}
\sum_{m\le N}d(m)^k \ll_k N\log^{2^k-1} N,
\label{4.140}
\end{equation}
we have
\[   \sum_{n\le N} D_R(n)^m \le  \sum_{n\le N} d(n)^{2m} \ll_m  N\log^{4^{m}-1}N. \]
Thus
\begin{equation}\begin{split}
\sum_{d_1,d_2,\ldots , d_m\le R} \frac{d(D_{m})}{D_m}
&\le \sum_{e_1,e_2,\ldots , e_m\le R} \frac{d(e_1)d(e_2)\ldots d(e_m)}{[e_1,e_2,\ldots , e_m]}\nonumber \\ &\ll_m  \log^{4^{m}-1}N + \frac{R^m\log^mR}{N}.
\end{split}
\end{equation}
Hence we conclude
\begin{equation}\begin{split}
\mathcal{S}_P(N,\mbox{\boldmath$j$}) =& N W_R(\mbox{\boldmath$j$}) +O(PR^{k-1}) \\ & + O_k\bigg(\Big(N\log^{4^{k-1}+ k-2}N  + R^{k-1}\log^{2k-2}N\Big) e^{-c_1\sqrt{\log\frac{P}{\max|j_i|^{k-1}}}}\bigg). \label{4.150}
\end{split}
\end{equation}
Taking $P=R$ proves the first part of Theorem \ref{Theorem4.1}. Equation \eqref{4.150} may also be useful when $P$ is not too large but larger than $R$.

We next turn to the case $P=N$. In this case $\psi_P(N;q,a)=\psi(N;q,a)+ O(\log N)$, the error term coming from $\Lambda_N(1)$.
We apply \eqref{1.290} and have
\[ \psi(N;D_{k-1},a) =[(D_{k-1},a)=1]  \frac{N}{\phi(D_{k-1})}  + E(N;D_{k-1},a) .\]
The condition that $(D_{k-1},a)=1$ is equivalent to having $(d_i,a)=1$ for $1\le i\le k-1$, and since $a\equiv -j_i (\mathrm{mod}\, d_i)$ these conditions are equivalent to $(d_i,j_i)=1$ for $1\le i \le k-1$. We
conclude that
\begin{equation}\psi(N;D_{k-1},a) =[(d_i,j_i)=1, \, 1\le i\le k-1]  \frac{N}{\phi(D_{k-1})}  + E(N;D_{k-1},a).\label{4.160}
\end{equation}
By \eqref{4.90} we thus obtain in place of \eqref{4.130}
\begin{equation}\begin{split}
\mathcal{S}_N(N,\mbox{\boldmath$j$}) &= NW_R(\mbox{\boldmath$j$}) + O\Big(\log^{k-1}N \sum_{d_1,d_2,\ldots , d_{k-1}\le R}\mu^2(D_{k-1})|E(N;D_{k-1},a)|\Big)\\&
\qquad +O(\log N\log^{k-1}R\prod_{i=1}^{k-1}d(1+j_i)), \label{4.170} \end{split}
\end{equation}
the last error term coming from the $n=1$ term. This last error will be negligible since it is  $\ll (\log N (\max_i|j_i|)^\epsilon)^k \ll_k N^\epsilon $ since $\max_i|j_i| \ll R^{\frac{1}{k}}$.
For the sum in the error term, we have
\begin{equation*}\begin{split}
\sum_{d_1,d_2,\ldots , d_{k-1}\le R}&\mu^2(D_{k-1})|E(N;D_{k-1},a)|\\&=
\sum_{m\le R^{k-1}}\mu^2(m)\max_{a (\mathrm{mod}\, m)}|E(N;m,a)|\sum_{\substack{ m=D_{k-1}\\ d_1,d_2,\ldots , d_{k-1}\le R}}1\\
&\le  \sum_{m\le R^{k-1}}\mu^2(m)\max_{a (\mathrm{mod}\, m)}|E(N;m,a)|d_{l}(m), \end{split}
\end{equation*}
where $l=2^{k-1}-1$.  The factor of $d_l(m)$ arises since, given $m$, the number of solutions of $m=D_{r}$
is bounded by $d_{2^r-1}(m)$, since the least common multiple  of $r$ squarefree numbers can always be expressed uniquely as the product of up to $2^r-1$ numbers which are pairwise relatively prime, determined by exactly which of the original numbers $d_1,d_2,\cdots , d_r$ each factor divides. (We will use this decomposition in later sections.) Applying Cauchy's inequality we see the previous expression is
\begin{equation*}
\ll \sqrt{ \sum_{m\le R^{k-1}}\frac{{d_l}(m)^2}{m}}\sqrt{\sum_{m\le R^{k-1}}m \max_{a (\mathrm{mod}\, m)}|E(N;m,a)|^2}.
\end{equation*}
We now use the generalization of \eqref{4.140}
\begin{equation} \sum_{m\le N}d_r(m)^k \ll_k N\log^{r^k-1} N,
\label{4.180}
\end{equation}
and the trivial estimate $|E(N;m,a)|\ll \frac{N\log N}{m}$
to see the error term above is
\begin{equation*}
\ll_l \sqrt{(\log R)^{l^2}N\log N \sum_{m\le R^{k-1}} \max_{a (\mathrm{mod}\, m)}|E(N;m,a)|}.
\end{equation*}
We now apply \eqref{1.300} to conclude this error is, for $R^{k-1}\le N^{\vartheta -\epsilon}$,
\[ \ll_k \frac{N}{(\log N)^{\frac{A}{2} - 4^{k-3/2}+2^{k-1}-1}},\]
which finishes the proof of Theorem \ref{Theorem4.1}.

\section{ Pair Correlation of $\Lambda_R(n)$}

We first prove the case $k=1$ of Theorem \ref{Theorem1}. We have
\begin{eqnarray} \mathcal{ S}_1(N,j,(1)) &=& \sum_{n\le N}\Lambda_R(n+j)=
\sum_{d\le R}\mu(d)\log \frac{R}{ d}\sum_{\stackrel{\scr 1\le n\le N } { \scr d|n+j}} 1 \nonumber \\
& =& N \sum_{d\le R}\frac{\mu(d)}{ d}\log\frac{R}{ d} + O(R) \nonumber \\
&=& N + O( \frac{N}{ (\log R)^A}) + O(R)  \label{5.10}
\end{eqnarray}
by Lemma \ref{Lemma1}. This proves Theorem \ref{Theorem1} in this case.

We now examine the case $k=2$ of Theorem \ref{Theorem1}.
These results  have been proved before in \cite{GO2}. The proof we give here models the procedure we will use for higher correlations without any of the truncation complications which will arise there.
In view of the comment following Theorem \ref{Theorem1}, we need to consider
\begin{equation} \mathcal{ S}_2(k)= \sum_{n=1}^N \Lambda_R(n)\Lambda_R(n+k),  \label{5.20}
\end{equation}
In our earlier notation, we have $\mathcal{ S}_2(k)= \mathcal{ S}_2(N, (0,k),(1,1))$ if $k\neq 0$, and $\mathcal{ S}_2(0) = \mathcal{ S}_2(N,(0),(2))$.

\begin{theorem} \label{Theorem5.1} We have
\begin{equation} \mathcal{ S}_2(0)  = N\log R + O(N) + O(R^2), \label{5.30}
\end{equation}
and for $0<|k| \le R $ and any $A>0$,
we have
\begin{equation} \mathcal{ S}_2(k) =  \gs_2(k)N + O(\frac{k}{ \phi(k)}\frac{N}{ (\log 2 R/k)^A}) +O(R^2).\label{5.40}
\end{equation}
\end{theorem}

Graham \cite{Gr} has proved \eqref{5.30} for $1 \le R\le N$ with the error term $O(R^2)$ removed. By Theorem \ref{Theorem5.1} we see that Theorem \ref{Theorem1} is true for  $k=2$, with $\mathcal{ C}(2)=1$ and  $\mathcal{ C}(1,1)=1$.

\emph{ Proof of Theorem \ref{Theorem5.1}.}  Applying the definition of $\Lambda_R(n)$, we have
\[\mathcal{ S}_2(k) =\sum_{d,e\le R}\mu(d) \log(R/d)\mu(e)\log(R/e) \sum_{\stackrel{\scr n\le N} { \scr d|n , e|n+k}}1 .\]
In the inner sum the two divisibility conditions imply that $n$ will run through a residue class modulo $[d,e]$ provided $(d,e)|k$, and there will be no solution for $n$ otherwise. Therefore we have
\begin{equation} \begin{split} \mathcal{ S}_2(k) &= N\sum_{\stackrel{\scr d,e\le R } { \scr (d,e)|k}} \frac{\mu(d) \log(R/d)\mu(e)\log(R/e)}{ [d,e]} + O\Big(\sum_{d,e\le R} \log(R/d)\log(R/e)\Big) \\
& = NT_2(k) +O(R^2). \label{5.50}
\end{split}
\end{equation}
To evaluate $T_2(k)$ we break the sum into relatively prime summands in order to handle $[d,e]$. We let $d=a_1b_{12}$ and $e=a_2b_{12}$ where $b_{12} = (d,e)$ so that $a_1$, $a_2$, and $b_{12}$ are pairwise relatively prime.  For higher correlations this decomposition notation will be used as well.  Hence we have
\begin{eqnarray*} T_2(k)&=& \sumprime_{\stackrel{\stackrel{\scr a_1 b_{12} \le R }{\scr a_2b_{12}\le R }} { \scr b_{12}|k}}\frac{\mu^2(b_{12}) }{ b_{12}} \frac{\mu(a_1)}{ a_1}\frac{\mu(a_2)}{ a_2} \log \frac{R}{ a_1b_{12}}\log\frac{R}{ a_2b_{12}} \nonumber \\
& = &\sumprime_{\stackrel{\scr a_2 b_{12} \le R  } { \scr b_{12}|k}}\frac{\mu^2(b_{12})}{ b_{12}}\frac {\mu(a_2)}{ a_2} \log \frac{R}{ a_2b_{12}} \sum_{\stackrel{\scr a_1\le R/b_{12} } { \scr (a_1,a_2b_{12})=1}}\frac{\mu(a_1)}{ a_1}\log\frac{R/b_{12}}{ a_1} ,
\end{eqnarray*}
where the prime on the summation indicates that all the summands are relatively prime to each other. We now apply Lemma \ref{Lemma1} with $j=0$  to obtain for any $B>0$
\[ T_2(k)= \sumprime_{\stackrel{\scr a_2 b_{12} \le R  } { \scr b_{12}|k}}\frac{\mu^2(b_{12})}{\phi( b_{12})}\frac {\mu(a_2)}{ \phi(a_2)} \log \frac{R}{ a_2b_{12}} + O\Bigg(\sumprime_{\stackrel{\scr a_2 b_{12} \le R  } { \scr b_{12}|k}}\frac{\mu^2(b_{12})\mu^2(a_2) \log \frac{R/b_{12}}{ a_2}}{ b_{12}a_2\log^B(2R/b_{12})} \Bigg).\]
For the main term above we sum over $a_2$ and apply Lemma \ref{Lemma1} again with $j=1$ to see this term is equal to
\[ \sum_{\stackrel{\scr b_{12}\le R } { \scr b_{12}|k}}\frac{\mu^2(b_{12})}{ \phi(b_{12})}\gs_2(b_{12}) + O\Big( \sum_{\stackrel{\scr b_{12}\le R } { \scr b_{12}|k}}\frac{\mu^2(b_{12})}{ \phi(b_{12})\log^B(2R/b_{12})}\Big).\]
Summing over $a_2$ in the error term in the formula for $T_2(k)$ above, we conclude
\begin{equation} T_2(k) =  \sum_{\stackrel{\scr b_{12}\le R } { \scr b_{12}|k}}\frac{\mu^2(b_{12})}{ \phi(b_{12})}\gs_2(b_{12}) + O\Big( \sum_{\stackrel{\scr b_{12}\le R } { \scr b_{12}|k}}\frac{\mu^2(b_{12})}{ \phi(b_{12})\log^{B-2}(2R/b_{12})}\Big).\label{5.60}
\end{equation}
We now consider two cases. If $k=0$ the main term is
\[ \sum_{ b_{12}\le R }\frac{\mu^2(b_{12})}{ \phi(b_{12})}\gs_2(b_{12})= \log R + O(1)\]
by Lemma \ref{Lemma4} with $j=1$, and the error term in \eqref{5.60} is also $O(1)$. This proves \eqref{5.30}. If $0<|k| \le R$ the main term is
\[ \sum_{b_{12}|k}\frac{\mu^2(b_{12})}{ \phi(b_{12})}\gs_2(b_{12}) = \gs_2(k)\]
 by the second part of Lemma \ref{Lemma4}, and the error term is
\[ \ll \frac{1}{ (\log 2R/k)^{B-2}}\sum_{b_{12}|k}\frac{\mu^2(b_{12})}{ \phi(b_{12})} = \frac{k}{ \phi(k)(\log 2R/k)^{B-2}},\]
which proves \eqref{5.40}.

It is worth noting that we can also give a very short proof of \eqref{5.40} using Theorem \ref{Theorem4.1} from the last section. With  $\mbox{\boldmath$j$}=(0,k)$, we have
\begin{eqnarray} W_R(\mbox{\boldmath$j$}) &= &\sum_{\substack{d_1 \le R \\ (d_1,k)=1}} \frac{\mu(d_1)}{\phi(d_1)}\log(R/d_1)\nonumber  \\
&=& \gs_2(k) + O(e^{-c_1\sqrt{\log R}}) \label{5.70}
\end{eqnarray}
by Lemma \ref{Lemma1}. By Theorem \ref{Theorem4.1} this proves \eqref{5.40} and also evaluates the mixed second correlation as well.

\section{ Triple correlation for $\Lambda_R(n)$}

To prove Theorem \ref{Theorem1} when $k=3$ we need to evaluate the sums
\begin{equation} \mathcal{ S}_3(k) = \sum_{n=1}^N {\Lambda_R}^2(n)\Lambda_R(n+k),\label{6.10}
\end{equation}
and, for non-zero $k_1\neq k_2$,
\begin{equation} \mathcal{ S}_3(k_1,k_2) = \sum_{n=1}^N \Lambda_R(n)\Lambda_R(n+k_1)
\Lambda_R(n+k_2)\label{6.20}
\end{equation}
In the notation of Theorem \ref{Theorem1} we have $ \mathcal{ S}_3(0)=\mathcal{ S}_3(N, (0),(3))$, $ \mathcal{ S}_3(k)=\mathcal{ S}_3(N, (0,k), (2,1))$ if $k\neq 0$, and
$\mathcal{ S}_3(k_1,k_2) = \mathcal{ S}_3(N, (0,k_1,k_2), (1,1,1))$ for non-zero $k_1\neq k_2$. We will obtain the following results on these correlations.

\begin{theorem} We have
\begin{equation}\mathcal{ S}_3(0) = \frac{3}{ 4}N\log^2R + O(N\log N(\log\log R)^{18}) + O(R^3),\label{6.30}
\end{equation}
and, for $k\neq 0$,  $|k|\le R^{\frac{1}{2}-\epsilon}$,
\begin{equation} \mathcal{ S}_3(k) = \gs_2(k) N \log R + O(N(\log\log R)^{13}) + O(R^3), \label{6.40}
\end{equation}
and letting $\mbox{\boldmath$k$}=(0,k_1,k_2)$, $k_1\neq k_2 \neq 0$, if $(k^*)^2<R/2$, where $k^*=\max(|k_1|,|k_2|)$, then
\begin{equation}  \mathcal{ S}_3(k_1,k_2) = \gs(\mbox{\boldmath$k$}) N + O\left(N e^{-c_1\sqrt{\log \big(\frac{R}{2(k^*)^2}\big) }}\log^8R\right) + O(R^3).\label{6.50}
\end{equation}
\label{theorem5}
\end{theorem}

We consider the general situation and specialize later. Let
\begin{equation}\mathcal{ S}_3(k_1,k_2,k_3) = \sum_{n=1}^N \Lambda_R(n+k_1)\Lambda_R(n+k_2)\Lambda_R(n+k_3)
\label{6.60}
\end{equation}
Expanding,  we have
\[\mathcal{ S}_3(k_1,k_2,k_3) = \sum_{d_1,d_2,d_3\le R}
\mu(d_1)\log (R/d_1)\mu(d_2)\log (R/d_2)
\mu(d_3)\log (R/d_3)\sum_{\stackrel{\stackrel{\stackrel{\scr n\le N}
 {\scr d_1|n+k_1 }}  {\scr d_2|n+k_2 }} { \scr d_3|n+k_3}}1.\]
The sum over $n$ is zero unless $(d_1,d_2)|k_2-k_1$, $(d_1,d_3)|k_3-k_1$, and $(d_2,d_3)|k_3-k_2$, in which case the sum runs through a residue class modulo $[d_1,d_2,d_3]$, and we have
\[\sum_{\stackrel{\stackrel{\stackrel{\scr n\le N}
 {\scr d_1|n+k_1 }}  {\scr d_2|n+k_2 }} { \scr d_3|n+k_3}}1 = \frac{N}{[d_1,d_2,d_3]} + O(1).\]
We conclude
\begin{equation}
\begin{array}{lcl}\mathcal{ S}_3(k_1,k_2,k_3) &= & {\displaystyle N\sum_{\stackrel{\stackrel{\stackrel{\scr d_1,d_2,d_3\le R} {\scr (d_1,d_2)|k_2-k_1 }}  {\scr  (d_1,d_3)|k_3-k_1} }{ \scr (d_2,d_3)|k_3-k_2}}
\frac{\mu(d_1)\log (R/d_1)\mu(d_2)\log (R/d_2)
\mu(d_3)\log (R/d_3)}{ [d_1,d_2,d_3] } + O(R^3) }\\
&=& N T_3(k_1,k_2,k_3) +O(R^3).
\end{array}
 \label{6.70}
\end{equation}
We now decompose $d_1$, $d_2$, and $d_3$ into relatively prime factors
\begin{eqnarray*} d_1&=&a_1b_{12}b_{13}a_{123} \\
d_2&=&a_2b_{12}b_{23}a_{123} \\
d_3&=&a_3b_{13}b_{23}a_{123}
\end{eqnarray*}
where $a_\chi$ or $b_\chi$ is a divisor of the $d_i$'s where $i$ occurs in $\chi$. Since the $d_i$'s are squarefree,  these new variables are pairwise relatively prime. The letters $a$ and $b$ reflect the parity of the number of $d_i$'s that the new variable divides. We will let $\mathcal{ D}$ denote the set of $a_\chi$'s and $b_\chi$'s which satisfy the conditions
\begin{eqnarray}  a_1b_{12}b_{13}a_{123}&\le& R\nonumber \\
a_2b_{12}b_{23}a_{123}&\le& R \label{6.80}  \\
a_3b_{13}b_{23}a_{123}&\le& R \nonumber \\
b_{12}a_{123}|k_2-k_1, \quad b_{13}a_{123}|k_3-k_1,& & b_{23}a_{123}|k_3-k_2 .\label{6.90}
\end{eqnarray}
Letting
\[ L_i(R) = \log \frac{R}{d_i}, \]
we have
\begin{equation} \begin{split}
 T_3(k_1,k_2,k_3) &= \sumprime_{\mathcal{ D}}\frac{ \mu(a_1)\mu(a_2)\mu(a_3)\mu^2(b_{12})\mu^2(b_{13})\mu^2(b_{23})\mu(a_{123}) }{ a_1 a_2 a_3 b_{12}b_{13}b_{23}a_{123}} L_1(R)L_2(R)L_3(R) \\
&= \sumprime_{\mathcal{D}} f_R(d_1,d_2,d_3) .\label{6.100}
\end{split}
\end{equation}
We now will sum over $a_1$, $a_2$, and $a_3$  using Lemma \ref{Lemma1} and Lemma \ref{Lemma3}. In order to apply these lemmas we need each $a_i$ to range over a long enough interval, and therefore we need to restrict the ranges of some of the other variables. The excluded ranges will later be shown to make a lower order contribution. If $D$ is a product of some of the variables in $\mathcal{D}$, we let  $\mathcal{D}(D)$ denote the subset of $\mathcal{D}$ where the variables not occuring in $D$ are eliminated from the inequalities in (\ref{6.80}) and divisibility conditions in (\ref{6.90}). Thus, letting  $D_1= a_2a_3b_{12}b_{13}b_{23}a_{123}$, we have that $\mathcal{D}(D_1)$ no longer includes the variable $a_1$ and we take  $a_1 =1$ in (\ref{6.80}). We now obtain on summing over $a_1$ using Lemma \ref{Lemma1} and taking $R_1<R$,
\begin{eqnarray} T_3(k_1,k_2,k_3) & = &
 \sumprime_{\stackrel{ \scr \mathcal{ D}}{\scr b_{12}b_{13}a_{123} \le R_1}}f_R(d_1,d_2,d_3) +
 \sumprime_{\stackrel{ \scr \mathcal{ D}}{\scr R_1< b_{12}b_{13}a_{123} \le R}}f_R(d_1,d_2,d_3) \nonumber \\
&=& \sumprime_{\stackrel{ \scr \mathcal{ D}(D_1)}{\scr b_{12}b_{13}a_{123} \le R_1}} \frac{ \mu(a_2)\mu(a_3)\mu(a_{123})\mu^2(D_1) }{ \phi(D_1)} L_2(R)L_3(R) + E_1(R) \nonumber \\ \qquad &+&  \sumprime_{\stackrel{ \scr \mathcal{ D}}{\scr R_1< b_{12}b_{13}a_{123} \le R}}f_R(d_1,d_2,d_3) \nonumber \\
&=& U_3(k_1,k_2,k_3) + E_1(R) + \mathcal{E}_f(\mathcal{D}) , \label{6.110}
\end{eqnarray}
where
\begin{eqnarray}
E_1(R)&=& \sumprime_{\stackrel{\scr \mathcal{ D} }{\scr b_{12}b_{13}a_{123}\le R_1}}\frac{ \mu(a_2)\mu(a_3)\mu(a_{123})\mu^2(D_1) }{ D_1}L_2(R)L_3(R)r_0(\frac{R}{ b_{12}b_{13}a_{123}},D_1)\nonumber \\
&\ll & e^{-c_1\sqrt{\log(R/R_1)}} \log^8 R   .\nonumber
\end{eqnarray}
Hence
\begin{equation} T_3(k_1,k_2,k_3) = U_3(k_1,k_2,k_3) + O( e^{-c_1\sqrt{\log(R/R_1)}} \log^8 R ) + \mathcal{E}_f(\mathcal{D}). \label{6.120}
\end{equation}
 Denote the summand for $U_3(k_1,k_2,k_3)$ by $g_R(d_1,d_2,d_3)$, which does not depend on $a_1$. Because of the symmetry in our original variables in (\ref{6.100}), we could equally well have summed over $a_2$ or $a_3$ above and obtained the same expression for $g_R(d_1,d_2,d_3)$ with the appropriate change in variables and renumbering of the $k_i's$. We will later make use of this fact for some of our error terms, and will let the summation conditions $\mathcal{D}(D)$ determine which variables appear in $g_R$ and subsequent summands.  Returning to (\ref{6.110}), we obtain on summing over $a_2$ using Lemma \ref{Lemma1} that, with $D_2=a_3b_{12}b_{13}b_{23}a_{123}$ and $R_2<R$,
\begin{eqnarray} U_3(k_1,k_2,k_3)  & = &
 \sumprime_{\stackrel{\stackrel{ \scr \mathcal{ D}(D_1)}{\scr b_{12}b_{13}a_{123} \le R_1}}{\scr b_{12}b_{23}a_{123}\le R_2}}g_R(d_1,d_2,d_3) + \sumprime_{\stackrel{\stackrel{ \scr \mathcal{ D}(D_1)}{\scr b_{12}b_{13}a_{123} \le R_1}}{\scr R_2 < b_{12}b_{23}a_{123}\le R}}g_R(d_1,d_2,d_3)
 \nonumber \\
&=& \sumprime_{\stackrel{\stackrel{ \scr \mathcal{ D}(D_2)}{\scr b_{12}b_{13}a_{123} \le R_1}}{\scr b_{12}b_{23}a_{123}\le R_2}} \frac{ \mu(a_3)\mu(a_{123})\mu^2(D_2) }{ \phi(D_2)} L_3(R)\gs_2(D_2)\nonumber \\  &+&  O(e^{-c_1\sqrt{\log(R/R_2)}} \log^6 R) + \sumprime_{\stackrel{\stackrel{ \scr \mathcal{ D}(D_1)}{\scr b_{12}b_{13}a_{123} \le R_1}}{\scr R_2 < b_{12}b_{23}a_{123}\le R}}g_R(d_1,d_2,d_3)
\nonumber \\
&=& V_3(k_1,k_2,k_3) +  O(e^{-c_1\sqrt{\log(R/R_2)}} \log^6 R) +\mathcal{E}_g(\mathcal{D}(D_1)). \label{6.130}
\end{eqnarray}
Finally, we denote the summand in $V_3(k_1,k_2,k_3)$ by $h_R(d_1,d_2,d_3)$ and let $D_3 =b_{12}b_{13}b_{23}a_{123}$ and $R_3<R$. Then by Lemma \ref{Lemma3} with $j=1$ we obtain
\begin{eqnarray} V_3(k_1,k_2,k_3)  & = &
 \sumprime_{\stackrel{\stackrel{\stackrel{ \scr \mathcal{ D}(D_2)}{\scr b_{12}b_{13}a_{123} \le R_1}}{\scr b_{12}b_{23}a_{123}\le R_2}}{\scr b_{13}b_{23}a_{123}\le R_3}}h_R(d_1,d_2,d_3) + \sumprime_{\stackrel{\stackrel{\stackrel{ \scr \mathcal{ D}(D_2)}{\scr b_{12}b_{13}a_{123} \le R_1}}{\scr b_{12}b_{23}a_{123}\le R_2}}{\scr R_3< b_{13}b_{23}a_{123}\le R}}h_R(d_1,d_2,d_3)  \nonumber \\
&=& -\sumprime_{\stackrel{\stackrel{\stackrel{ \scr \mathcal{ D}(D_3)}{\scr b_{12}b_{13}a_{123} \le R_1}}{\scr b_{12}b_{23}a_{123}\le R_2}}{\scr b_{13}b_{23}a_{123}\le R_3}} \frac{ \mu^2(D_3)\mu(a_{123})\mu((D_3,2)) }{ \phi(D_3)}\gs_2(2D_3)\gs_3(2D_3) \nonumber \\ \qquad &+& O(e^{-c_1\sqrt{\log(R/R_3)}} \log^4R)+\sumprime_{\stackrel{\stackrel{\stackrel{ \scr \mathcal{ D}(D_2)}{\scr b_{12}b_{13}a_{123} \le R_1}}{\scr b_{12}b_{23}a_{123}<R_2}}{\scr R_3< b_{13}b_{23}a_{123}<R}}h_R(d_1,d_2,d_3)\nonumber \\
&=& W_3(k_1,k_2,k_3) +  O(e^{-c_1\sqrt{\log(R/R_3)}} \log^4 R) + \mathcal{E}_h(\mathcal{D}(D_2)). \label{6.140}
\end{eqnarray}

We now prove Theorem \ref{theorem5} by considering each case separately.
We first  prove (\ref{6.50}) which is the case where the error terms are the easiest to handle.

\section{Evaluation of $S_3(k_1,k_2)$}

We  consider $\mathcal{ S}_3(k_1,k_2)$ by taking $k_1\neq k_2\neq 0$ and $k_3=0$ in (\ref{6.90}) which therefore becomes
\begin{equation}
b_{12}a_{123}|k_2-k_1, \quad b_{13}a_{123}|k_1, \quad b_{23}a_{123}|k_2 .\label{7.10}
\end{equation}
These conditions imply, letting $k^*=\max(|k_1|,|k_2|)$, that
\begin{equation} b_{12}b_{13}a_{123} < 2(k^*)^2, \quad
b_{12}b_{23}a_{123} < 2(k^*)^2, \quad b_{13}b_{23}a_{123} <2(k^*)^2.\label{7.20}
\end{equation}
Hence, taking $R_1=R_2=R_3= 2(k^*)^2 <R$, we see that the error terms
$\mathcal{E}_f$, $\mathcal{E}_g$, and $\mathcal{E}_h$ in
\eqref{6.110}, \eqref{6.130}, and \eqref{6.140} are identically zero, and therefore
\begin{equation}
T_3(k_1,k_2) = W_3(k_1,k_2) + O(e^{-c_1\sqrt{\log\big(\frac{R}{2(k^*)^2}\big)}}\log^8 R).
\label{7.30}
\end{equation}
Now that the variables $a_1$, $a_2$, and $a_3$ have been eliminated, the bounds on the variables in $\mathcal{D}$ are automatically satisfied
from (\ref{7.20}), and, provided $(k^*)^2 < R/2$,  we have
\begin{equation}
W_3(k_1,k_2) = - \sumprime_{\stackrel{\stackrel{\scr b_{13}a_{123}|k_1 } { \scr b_{23}a_{123}|k_2}} { \scr b_{12}a_{123}|k_2-k_1}}\frac{ \mu^2(D_3)\mu(a_{123})\mu((D_3,2))}{ \phi(D_3)} \gs_2(2D_3)\gs_3(2D_3).\label{7.40}
\end{equation}
This sum is over square-free divisors, and therefore we let
\begin{equation}
k_1 = s_1K_1, \quad k_2 = s_2K_2,
\label{7.50}
\end{equation}
where $K_1$ and $K_2$ are the largest square-free divisors of $k_1$ and $k_2$, and let
\begin{equation}
\kappa = (k_1,k_2), \quad K_{12} = (K_1,K_2) .
\label{7.60}
\end{equation}
Then we may rewrite $W_3$ as
\begin{equation}
W_3(k_1,k_2) = - \sumprime_{\stackrel{\stackrel{\scr b_{13}a_{123}|K_1 } { \scr b_{23}a_{123}|K_2}} { \scr b_{12}a_{123}|k_2-k_1}}\frac{ \mu^2(D_3)\mu(a_{123})\mu((D_3,2))}{ \phi(D_3)} \gs_2(2D_3)\gs_3(2D_3).\label{7.70}
\end{equation}
The proof of \eqref{6.50} will follow from \eqref{6.70}, \eqref{7.30} and Lemma \ref{Lemma5} once we prove that
\begin{equation} W_3(k_1,k_2) =\gs_2(K_{12})\gs_3(\Delta) = \gs_2(\kappa)\gs_3(\Delta) = \gs(\mbox{\boldmath$k$}). \label{7.80}
\end{equation}
We now let
\begin{equation}
 \quad K_1 = K_{12} j_1, \quad K_2 = K_{12} j_2, \quad k_2-k_1 = K_{12}(s_2j_2-s_1j_1) \label{7.90}
\end{equation}
where $K_{12}$, $j_1$, and $j_2$ are square-free and pairwise relatively prime ($s_2j_2-s_1j_1$ may not be square-free or relatively prime with $K_{12}$.)
Next, let
\begin{equation}
b_{13}=c_{13}d_{13}, \quad b_{23}=c_{23}d_{23} \label{7.100}
\end{equation}
where $c_{13},c_{23}|K_{12}$, $d_{13}|j_1$, $d_{23}|j_2$, and $c_{13}$, $d_{13}$, $c_{23}$, $d_{23}$, $b_{12}$, $a_{123}$ are thus all pairwise relatively prime. We also see that $a_{123}|K_{12}$. Thus, with $D_3=b_{12}c_{13}c_{23}d_{13}d_{23}a_{123}$,
\begin{equation}
W_3(k_1,k_2) = - \sumprime_{c_{13}c_{23}a_{123}|K_{12}} \mu(a_{123})\sumprime_{\stackrel{\stackrel{\scr b_{12}|k_2-k_1 } { \scr d_{13}|j_1}} { \scr d_{23}|j_2}}\frac{ \mu^2(D_3)\mu((D_3,2))}{ \phi(D_3)} \gs_2(2D_3)\gs_3(2D_3).\label{7.110}
\end{equation}
We first sum over $d_{13}$ in the inner sum. To do this, we take $D_4=b_{12}c_{13}c_{23}d_{23}a_{123}$ and have
\begin{eqnarray*}
\lefteqn{\sum_{\stackrel{\scr d_{13}|j_1}{(d_{13},D_4)=1}}\frac{ \mu^2(D_3)\mu((D_3,2))}{ \phi(D_3)} \gs_2(2D_3)\gs_3(2D_3)  }\\ & & =2C_2\frac{\mu^2(D_4)}{\phi(D_4)}\sum_{\stackrel{\scr d_{13}|j_1}{(d_{13},D_4)=1}}\frac{\mu^2(d_{13})}{\phi(d_{13})}\mu((d_{13}D_4,2))H_2(2d_{13}D_4)\gs_3(2d_{13}D_4).
\end{eqnarray*}
We now break the sum on the right into two sums according to whether  $d_{13}$ is even or odd, in the former case we let $d_{13}=2d$, and on using Lemma \ref{Lemma4} we obtain that the right-hand side is
\begin{eqnarray*} \lefteqn{ = \frac{2C_2\mu^2(D_4)H_2(D_4)}{\phi(D_4)} \Bigg\{ -[(D_4,2)= 1] \sum_{\stackrel{\scr 2d|j_1}{(d,2D_4)=1}}\frac{\mu^2(d)}{\phi_2(d)}\gs_3(2dD_4) }\\ & & \qquad \qquad + \mu((D_4,2))\sum_{\stackrel{\scr d_{13}|j_1}{(d_{13},2D_4)=1}}\frac{\mu^2(d_{13})}{\phi_2(d_{13})}\gs_3(2d_{13}D_4)\Bigg\} \\
&=&  \frac{2C_2\mu^2(D_4)H_2(D_4)}{\phi(D_4)}\Big( -[(D_4,2)=1][2|j_1]+ \mu((D_4,2))\Big)\gs_3(2D_4j_1).
\end{eqnarray*}
We will denote
\begin{equation}
B(d)=  \frac{2C_2\mu^2(d)H_2(d)}{\phi(d)}.
\label{7.120}
\end{equation}
Now substituting into \eqref{7.110}, and  letting $D_5=b_{12}c_{13}c_{23}a_{123}$, the sum over $d_{23}$ is equal to
\begin{equation*}
B(D_5)\sum_{\stackrel{d_{23}|j_2}{(d_{23},D_5)=1}}\frac{\mu^2(d_{23}) H_2(d_{23})}{\phi(d_{23})}
 \Big( -[(d_{23}D_5,2)=1][2|j_1]+ \mu((d_{23}D_5,2))\Big)\gs_3(2d_{23}D_5j_1).
\end{equation*}
Since $(j_1,j_2)=1$, and $d_{23}|j_2$, we may replace the condition $(d_{23},D_5)=1$ in the sum by $(d_{23},j_1D_5)$ and divide the sum into two sums with even or odd terms as above to see that this expression is
\begin{eqnarray*}&=& -B(D_5)[(D_5,2)=1][2|j_2]\sum_{\stackrel{\scr 2d|j_2}{\scr (d,2j_1D_5)=1}} \frac{\mu^2(d)}{\phi_2(d)}\gs_3(2dD_5j_1) \\
& & \quad + B(D_5)\Big(-[(D_5,2)=1][2|j_1]+\mu((D_5,2))\Big)\sum_{\stackrel{\scr d_{23}|j_2}{\scr (d_{23},2j_1D_5)=1}} \frac{\mu^2(d_{23})}{\phi_2(d_{23})}\gs_3(2d_{23}D_5j_1) \\
&= & B(D_5)\Big( -[(D_5,2)=1]\big([2|j_1]+[2|j_2]\big)+\mu((D_5,2))\Big)\gs_3(2D_5j_1j_2).
\end{eqnarray*}
Now let
\begin{equation}
\Delta = k_1k_2(k_2-k_1)= s_1s_2{K_{12}}^2j_1j_2(k_2-k_1) = s_1s_2{K_{12}}^3j_1j_2(s_2j_2-s_1j_1). \label{7.130}
\end{equation}
We  substitute the last result into (\ref{7.110}) and sum over $b_{12}$.  Let  $D_6 = c_{13}c_{23}a_{123}$. We claim that the relatively prime condition  $(b_{12},D_6)=1$ may be replaced by $(b_{12}, D_6j_1j_2)=1$. To see this, note that $j_1$, $j_2$, and $K_{12}$ are pairwise relatively prime, and further $(s_1,j_2)=1$ and $(s_2,j_1)=1$. Thus
\[ (k_2 -k_1,j_1) = (K_{12}(s_2j_2-s_1j_1),j_1) = (s_2,j_1) =1 \]
and similarly $(k_2-k_1,j_2)=1$. Hence,  since $b_{12}|k_2-k_1$, we have $(b_{12},j_1j_2)=1$.  Now summing over $b_{12}$ our sum is
\begin{eqnarray*}
&=& B(D_6)\sum_{\stackrel{\scr b_{12}|k_2-k_1}{\scr (b_{12},D_6j_1j_2)=1}}\Bigg( -[(b_{12}D_6,2)=1]\big([2|j_1]+[2|j_2]\big)+\mu((b_{12}D_6,2))\Bigg) \\
& & \hskip 1.7in \times \  \frac{\mu^2(b_{12}) H_2(b_{12})}{\phi(b_{12})}\gs_3(2b_{12}D_6j_1j_2) \\
&=& B(D_6)\sum_{\stackrel{\scr 2b|k_2-k_1}{\scr (b,2D_6j_1j_2)=1}}\frac{\mu^2(b)}{\phi_2(b)}\Big( -[2|k_2-k_1][(D_6,2)=1]\Big)\gs_3(2bD_6j_1j_2) \\
& &  + B(D_6)\sum_{\stackrel{\scr b_{12}|k_2-k_1}{\scr (b_{12},2D_6j_1j_2)=1}}\frac{\mu^2(b_{12})}{\phi_2(b_{12})}
\Big( -[(D_6,2)=1]\big([2|j_1]+[2|j_2]\big)+\mu((D_6,2))\Big)\gs_3(2b_{12}D_6j_1j_2)\\
&=& B(D_6) \Big( -[(D_6,2)=1]\big([2|j_1]+[2|j_2]+[2|k_2-k_1]\big)+\mu((D_6,2))\Big)\gs_3(2D_6j_1j_2(k_2-k_1)).
\end{eqnarray*}

Now $D_6|K_{12}|\Delta$, and hence
\[ \gs_3(2D_6j_1j_2(k_2-k_1))=\gs_3(2D_6\Delta)= \gs_3(\Delta).  \]
We conclude that
\begin{equation}
W_3(k_1,k_2) = -\gs_3(\Delta)\sumprime_{D_6|K_{12}}\mu(a_{123})B(D_6)\Big( -[(D_6,2)=1]\big([2|j_1]+[2|j_2]+[2|k_2-k_1]\big)+\mu((D_6,2))\Big).
\label{7.140}
\end{equation}

If $K_{12}$ is odd, then exactly one of the variables $j_1$, $j_2$, and $k_2-k_1$ is even and the other two are odd.  If $K_{12}$ is even, then $j_1$ and $j_2$ are odd and $k_2-k_1$ is even. Hence in either case
\[ [2|j_1]+[2|j_2]+[2|k_2-k_1]=1 ,\]
and therefore
\begin{equation}
W_3(k_1,k_2) = -\gs_3(\Delta)\sumprime_{D_6|K_{12}}\mu(a_{123})B(D_6)\big( -[(D_6,2)=1]+\mu((D_6,2))\big).
\label{7.150}
\end{equation}
When $K_{12}$ is odd the expression in parentheses is zero, and hence
\begin{equation*}
W_3(k_1,k_2) = 0, \quad \textrm{ if $K_{12}$ is odd.}
\end{equation*}
If $K_{12}$ is even, then the expression in parentheses in \eqref{7.150} is zero if $D_6$ is odd, and is equal to $-1$ when $D_6$ is even. Hence, we conclude
\begin{equation}
W_3(k_1,k_2) = [2|K_{12}]\gs_3(\Delta)\sumprime_{\substack{D_6|K_{12}\\ 2|D_6}}\mu(a_{123})B(D_6).\label{7.160}
\end{equation}
We could now evaluate this sum as before by summing over each variable in turn, but there is an easier approach, based on the observation that if a square-free number is a product of some factors, then necessarily those factors must be relatively prime with each other. Let $z=vwy$, and $\mathcal{A}$ be a set of natural numbers. Then for any arithmetic function $a(z)$ we have
\begin{equation} \begin{split}
\sumprime_{z\in \mathcal{A}} \mu^2(z)a(z)\mu(y) &=
\sum_{z\in \mathcal{A}} \mu^2(z)a(z)\sum_{w|z}\sum_{y|
\frac{z}{w}}\mu(y) \\
&=  \sum_{z\in \mathcal{A}} \mu^2(z)a(z)\sum_{\substack{w|z\\ w=z}}1 \\
&= \sum_{z\in \mathcal{A}} \mu^2(z)a(z).\label{7.170}
\end{split} \end{equation}
The sum in \eqref{7.160} is of this form, and therefore we have
\begin{equation*} \begin{split}
W_3(k_1,k_2) &= [2|K_{12}]\gs_3(\Delta)\sum_{\substack{z|K_{12}\\ 2|z}}B(z) \\
& = 2C_2[2|K_{12}]\gs_3(\Delta) \sum_{\substack{2z'|K_{12} \\ (z',2)=1}}\frac{\mu^2(z')}{\phi_2(z')} \\
&= 2C_2[2|K_{12}]\gs_3(\Delta)H_2(K_{12}),
\end{split}\end{equation*}
where we used \eqref{2.30} in the last line.
We conclude by \eqref{2.50} that
\begin{equation}
W_3(k_1,k_2) = \gs_2(K_{12})\gs_3(\Delta) = \gs_2(\kappa)\gs_3(\Delta),
\label{7.180}
\end{equation}
which completes the proof of (\ref{7.80}).

\section{ Evaluation of $S_3(k),\ k\neq 0$}

We now take $k_1=k \neq 0 $, $k_2=k_3=0$ so that the divisibility conditions in (\ref{6.90}) become
\[
b_{12}a_{123}|k, \quad b_{13}a_{123}|k,
\]
 and since our variables are relatively prime
\begin{equation}
b_{12}b_{13}a_{123} |k. \label{8.10}
\end{equation}
Taking $R_1 =k<R$ we see that $\mathcal{E}_f$ in (\ref{6.120}) is identically zero.
We will take $R_2=R_3 \ge k^2$ to be chosen later as a function of $R$, and conclude from \eqref{6.120}, \eqref{6.130}, and  \eqref{6.140} that
\begin{equation}
T_3(k)= W_3(k) + \mathcal{E}_g(\mathcal{D}(D_1))+ \mathcal{E}_h(\mathcal{D}(D_2))+O(e^{-c_1\sqrt{\log(R/R_2)}}\log^8R). \label{8.20}
\end{equation}
Now, summing over $a_3$ using Lemma \ref{Lemma1}, we have
\begin{eqnarray}  \mathcal{E}_g(\mathcal{D}(D_1)) &=& \sumprime_{\stackrel{\stackrel{\stackrel{ \scr \mathcal{ D}(D_1)}{\scr b_{12}b_{13}a_{123} \le R_1}}{\scr R_2 < b_{12}b_{23}a_{123}\le R}}{\scr b_{13}b_{23}a_{123}\le R_3}}g_R(d_1,d_2,d_3)
+ \sumprime_{\stackrel{\stackrel{\stackrel{ \scr \mathcal{ D}(D_1)}{\scr b_{12}b_{13}a_{123} \le R_1}}{\scr R_2 < b_{12}b_{23}a_{123}\le R}}{\scr  R_3< b_{13}b_{23}a_{123}\le R}}g_R(d_1,d_2,d_3) \nonumber \\
&=& \sumprime_{\stackrel{\stackrel{\stackrel{ \scr \mathcal{ D}(a_2D_3)}{\scr b_{12}b_{13}a_{123} \le R_1}}{\scr R_2 < b_{12}b_{23}a_{123}\le R}}{\scr b_{13}b_{23}a_{123}\le R_3}}h_R(d_1,d_2,d_3) +O(e^{-c_1\sqrt{\log(R/R_3)}}\log^6R) \nonumber \\
&+& \sumprime_{\stackrel{\stackrel{\stackrel{ \scr \mathcal{ D}(D_1)}{\scr b_{12}b_{13}a_{123} \le R_1}}{\scr R_2 < b_{12}b_{23}a_{123}\le R}}{\scr  R_3< b_{13}b_{23}a_{123}\le R}}g_R(d_1,d_2,d_3) . \label{8.30}
\end{eqnarray}
The first sum above is the same as $\mathcal{E}_h(\mathcal{D}(D_2))$ with the appropriate relabeling of variables, and the estimate we now obtain applies to both expressions. We see in this sum that $a_2\le \frac{R}{b_{12}b_{23}a_{123}} \le \frac{R}{R_2}$, and hence  the sum is
\begin{eqnarray}
&\ll& \log(R/R_2)\sum_{ b_{12}b_{13}a_{123}|k } \sum_{ \frac{R_2}{b_{12}a_{123}} <  b_{23} \le \frac{R}{b_{12}a_{123}}} \frac{\mu^2(D_3)}{\phi(D_3)}\sum_{\substack{a_2\le \frac{R}{R_2}\\ (a_2,D_3)=1}}\frac{\mu^2(a_2)}{\phi(a_2)}\gs_2(a_2D_3)
\nonumber \\
&\ll& \log^2(R/R_2)\log\log R\sum_{ b_{12}b_{13}a_{123}|k } \sum_{ \frac{R_2}{b_{12}a_{123}} <  b_{23} \le \frac{R}{b_{12}a_{123}}} \frac{\mu^2(D_3)}{\phi(D_3)}
\nonumber \\
&\ll & \log^3(R/R_2)\log\log R\left(\sum_{b|k}\frac{\mu^2(b)}{\phi(b)}\right)^3 \nonumber
\\
&\ll & \log^3(R/R_2)\log\log R(\frac{k}{\phi(k)})^3 ,\nonumber
\end{eqnarray}
by Lemma \ref{Lemma4}, and hence, since $\frac{k}{\phi(k)} \ll \log \log 3k$,
\begin{equation}
\mathcal{E}_h(\mathcal{D}(D_2)) \ll \log^3(R/R_2)(\log \log R)^4. \label{8.40}
\end{equation}
Similarly, for the second sum in (\ref{8.30}) both $a_2 \le \frac{R}{R_2}$ and $a_3 \le \frac{R}{R_2}$, so that the sum is
\[
\ll  \log^5(R/R_2)(\log\log 3k)^3\]
and hence
\begin{equation}
\mathcal{E}_g(\mathcal{D}(D_1)) \ll \log^5(R/R_2)(\log \log R)^3+e^{-c_1\sqrt{\log(R/R_2)}}\log^6 R. \label{8.50}
\end{equation}
We take
\begin{equation}
R_2 = R_3 =Re^{- c_2(\log\log R)^2} \label{8.60}
\end{equation}
where $c_2$ is a sufficiently large constant, and conclude by \eqref{8.20}, \eqref{8.40}, and \eqref{8.50} that, for $|k|\le R^{\frac{1}{2}-\epsilon}$,
\begin{equation}
T_3(k) = W_3(k) + O\big((\log\log R)^{13}\big), \label{8.70}
\end{equation}
where
\begin{equation}
 W_3(k) = - \sumprime_{\stackrel{\stackrel{\scr b_{12} b_{23}a_{123}\le R_2 } { \scr b_{13}b_{23}a_{123}\le R_2}} {
\scr b_{12}b_{13} a_{123}|k}} \frac{ \mu^2(D_3)\mu(a_{123})\mu(( D_3,2))}{ \phi(D_3)} \gs_2(2D_3)\gs_3(2D_3),\label{8.80}
\end{equation}
and as before $D_3=b_{12}b_{13}b_{23}a_{123}$.
Since $b_{23}$ is the only variable not constrained
by the divisibility condition, we have
on letting
\begin{equation} R_4 = \frac{R_2}{ a_{123}} \min(\frac{1}{ b_{12}}, \frac{1}{ b_{13}}),\label{8.90}
\end{equation}
and $E= b_{12}b_{13}a_{123}$, that
\[W_3(k) = - \sumprime_{E|k}\frac{\mu^2(E)}{ \phi(E)}\mu(a_{123})\sum_{\stackrel{\scr b_{23}\le R_4 } { \scr (b_{23}, E)=1}}\frac{\mu^2(b_{23})}{ \phi(b_{23})}\mu((b_{23}E,2)) \gs_2(2b_{23}E)\gs_3(2 b_{23} E) .\]
We break the inner sum into two subsums according to whether $b_{23}$ is even or odd, the former case forcing $(E,2)=1$. We thus obtain, on taking $b_{23} = 2b$ in the first subsum and applying (\ref{2.40})
\begin{eqnarray}  W_3(k) &=&  \sumprime_{\stackrel{\scr E|k} { \scr (E,2)=1}}\frac{\mu^2(E)\gs_2(2E)}{ \phi(E)}\mu(a_{123})\sum_{\stackrel{\scr b\le \frac{R_4}{ 2} } { \scr (b, 2E)=1}}\frac{\mu^2(b)}{ \phi_2(b)} \gs_3(2 b E) \nonumber \\
& &\qquad - \sumprime_{E|k}\frac{\mu^2(E)\gs_2(2E)}{ \phi(E)}\mu((E,2))\mu(a_{123})\sum_{\stackrel{\scr b_{23}\le R_4 } { \scr (b_{23}, 2E)=1}}\frac{\mu^2(b_{23})}{ \phi_2(b_{23})}\gs_3(2 b_{23} E) .\nonumber
\end{eqnarray}
The inner sums in both sums above are by Lemma \ref{Lemma4} and the estimate \eqref{2.190}
\[
= \log R_4 +O(\log\log 3k) + O(\frac{\exp\left(\frac{c\sqrt{\log k}}{\log\log 3k}\right)}{\sqrt{R_4}}).\]
 Substituting we now see that the main term from the first sum is canceled by the part of the main term in the second sum when $(E,2)=1$, and hence we obtain
\begin{equation*}\begin{split}  & W_3(k)\\ &= \sumprime_{ E|k }\frac{\mu^2(E)\gs_2(2E)}{ \phi(E)}\Big( \mu(a_{123})[2|E]\log R_4 + +O(\log\log 3k) +  O(\frac{\exp\left(\frac{c\sqrt{\log k}}{\log\log 3k}\right)}{\sqrt{R_4}})\Big)
 \\
 & = \log R_2 \sumprime_{ \substack{E|k \\ 2|E}}\frac{\mu^2(E)\gs_2(2E)}{ \phi(E)} \mu(a_{123})   \\ & \quad + O\left(\sumprime_{ E|k }\frac{\mu^2(E)\gs_2(2E)}{ \phi(E)}\Big( \log E + \log\log 3k + \frac{\exp\left(\frac{c\sqrt{\log k}}{\log\log 3k}\right)}{\sqrt{\frac{R_2}{E}}} \right)  \\
 &=Y_3(k)\log R_2  + O\left( \big(\sum_{b|k}\frac{\mu^2(b)\log b}{\phi(b)}\big)^3\log \log 3k\right)  \\
&  \quad +O\left(\frac{\exp\left(\frac{c\sqrt{\log k}}{\log\log 3k}\right)}{\sqrt{R_2}}\big( \sum_{b|k}\frac{\mu^2(b)\sqrt{b}}{\phi(b)}\big)^3\right) ,\end{split}
\end{equation*}
where we used the estimate $\gs_2(k) \ll \log\log 3k$ in the last two error terms. By Lemma 5 of \cite{GY}, the sum in the first error term is
\[ = \frac{kh_1(k)}{\phi(k)} \ll (\log \log 3k)^2 \]
where $h_1(k)$ is given by \eqref{2.170}. The sum in the second error term is
\[ \ll \exp\left(c\sum_{n<\log k}\frac{\mu^2(n)\sqrt{n}}{\phi(n)} \right) \ll e^{c\sqrt{\log k}} \]
and hence we conclude that
\begin{equation} W_3(k) = Y_3(k)\log R_2 +O((\log\log R)^7). \label{8.100}
\end{equation}
We complete the evaluation of $W_3(k)$ by proving
\begin{equation}
Y_3(k) = \gs_2(k). \label{8.110}
\end{equation}
We note first that if $2\ndiv k$ then the sum defining $Y_3(k)$ is empty and therefore $Y_3(k)=0$, in agreement with $\gs_2(k)$.
Therefore by \eqref{7.170}
\begin{equation*}\begin{split}
Y_3(k) &= [2|k] \sum_{ \substack{E|k \\ 2|E}}\frac{\mu^2(E)\gs_2(2E)}{ \phi(E)} \\
&= 2C_2[2|k]\sum_{\substack{2E'|k \\ (E',2)=1}}\frac{\mu^2(E')}{ \phi_2(E')} \\
&= 2C_2[2|k]H_2(k) =\gs_2(k),
\end{split}\end{equation*}
which proves \eqref{8.110}.

\section{ Evaluation of $S_3(0)$: Main Term}
We now consider $k_1=k_2=k_3=0$.   On applying equations (\ref{6.120}), (\ref{6.130}), and (\ref{6.140}) and taking
\[ S= R_1=R_2=R_3 \]
we have that
\begin{equation} T_3(0) = W_3(0) + O(e^{-c_1\sqrt{\log(R/S)}} \log^8 R) +\mathcal{E}_f(\mathcal{D})+ \mathcal{E}_g(\mathcal{D}(D_1))+ \mathcal{E}_h(\mathcal{D}(D_2)),\label{9.10}
\end{equation}
where we relabel variables for simplicity and have
\begin{equation}
 W_3(0) = - \sumprime_{\stackrel{\stackrel{\scr uvy \le S } { \scr uwy\le S}}
{ \scr vwy\le S}} \frac {\mu(y)\mu^2(D)\mu((D, 2))}{ \phi(D)} \gs_2(2D)\gs_3(2D),\label{9.20}
\end{equation}
and $D=uvwy$. We will prove in this section that
\begin{equation}
W_3(0) = \frac{3}{4}\log^2S + O(\log S). \label{9.30}
\end{equation}
In the next section we will prove that $\mathcal{E}_f(\mathcal{D})$, $\mathcal{E}_g(\mathcal{D}(D_1))$, and $\mathcal{E}_h(\mathcal{D}(D_2))$
all satisfy the bound
\begin{equation}
\ll  \log^9(\frac{R}{S})\log R \label{9.40}
\end{equation}
from which we conclude from (\ref{9.10}) and (\ref{9.30})
on taking
\[ S = Re^{-c_2(\log\log R)^2} \]
that
\begin{equation}
T_3(0) = \log^2R +O(\log R (\log \log R)^{18})
\label{9.50}
\end{equation}
and therefore by (\ref{6.70}) we obtain (\ref{6.30}).

The first step in evaluating $W_3(0)$ is to sum over $u$ in (\ref{9.20}) by separating into sums according to whether $u$ is even or odd, and apply Lemma \ref{Lemma4}. We let $z=vwy$, and have
\begin{equation*}\begin{split}\sum_{\stackrel{\scr u \le \min(\frac{S}{vy},\frac {S}{wy})}{\scr (u,z)=1}}\frac{\mu^2(u)}{\phi(u)}& \mu((uz,2))\gs_2(2uz)\gs_3(2uz) \\ &=
2C_2H_2(2z) \sum_{\stackrel{\scr u \le \frac{S}{z}\min(v,w)}{\scr (u,z)=1}}\frac{\mu^2(u)}{\phi(u)}\mu((uz,2))H_2(2u)\gs_3(2uz) \\
&= 2C_2H_2(2z)\Bigg(-[(z,2)=1] \sum_{\stackrel{\scr 2u' \le \frac{S}{z}\min(v,w)}{\scr (u',2z)=1}}\frac{\mu^2(u')}{\phi_2(u')}\gs_3(2u'z) \\
& \quad +\mu((z,2))\sum_{\stackrel{\scr u \le \frac{S}{z}\min(v,w)}{\scr (u,2z)=1}}\frac{\mu^2(u)}{\phi_2(u)}\gs_3(2uz).\Bigg)\\
&= 2C_2H_2(2z)\Bigg(-[(z,2)=2]\bigg( \log\Big(\frac{S(z,3)}{3z}\min(v,w)\Big) +D_3 \\
& \qquad \qquad +h_3(6z)\bigg) + [(z,2)=1]\log 2 + O\Big(\frac{m(6z)}{\sqrt{\frac{S}{z}\min(v,w)}}\Big)\Bigg).
\end{split}
\end{equation*}
Substituting this result into (\ref{9.20}) we obtain
\begin{equation}\begin{split}
W_3(0) = &2C_2 \sumprime_{\stackrel{\scr z\le S}{\scr 2|z}}\frac{\mu^2(z)\mu(y)}{\phi(z)}H_2(2z)
\Big( \log\bigg(\frac{S(z,3)}{3z}\min(v,w)\bigg) +D_3 +h_3(6z)\Big)\\
&  - 2C_2 \log 2\sumprime_{\stackrel{\scr z\le S}{\scr (z,2)=1}}\frac{\mu^2(z)\mu(y)}{\phi(z)}H_2(2z) +
O\Big(   \sumprime_{ z\le S} \frac{\mu^2(z)m(6z)}{\phi(z)\sqrt{\frac{S}{z}\min(v,w)}}H_2(2z)\Big).
\label{9.60}
\end{split}
\end{equation}
Except for the factor of $\log \min(v,w)$, the sums above are of the form
considered in \eqref{7.170}.
Hence we have
\begin{eqnarray}
2C_2 \sumprime_{\stackrel{\scr z\le S}{\scr 2|z}}\frac{\mu^2(z)\mu(y)}{\phi(z)}H_2(2z)\log \frac{S}{z} &=& 2C_2 \sum_{\stackrel{\scr z\le S}{\scr 2|z}}\frac{\mu^2(z)}{\phi(z)}H_2(2z)\log\frac{S}{z} \nonumber \\
&=& 2C_2 \sum_{\stackrel{\scr 2z'\le S}{\scr (z',2)=1}}\frac{\mu^2(z')}{\phi_2(z')} \log \frac{S}{2z'}\nonumber \\
&=& \frac{1}{2} \log^2S + O(\log S), \label{9.70}
\end{eqnarray}
where we used \eqref{2.210}.
Similarly we have, by Lemma \ref{Lemma2}
\begin{equation}
2C_2 \sumprime_{\stackrel{\scr z\le S}{\scr 2|z}}\frac{\mu^2(z)\mu(y)}{\phi(z)}H_2(2z) = \log S + O(1).\label{9.80}
\end{equation}
We also see that the condition $2|z$ in \eqref{9.80} may be replaced by $(z,2)=1$ and \eqref{9.80} will still hold. Similarly it is clear that
\begin{equation}
2C_2 \sumprime_{\stackrel{\scr z\le S}{\scr 2|z}}\frac{\mu^2(z)\mu(y)}{\phi(z)}H_2(2z)\log((z,3)) \ll \log S, \label{9.90}
\end{equation}
and also that
\begin{equation}\begin{split}
2C_2 &\sumprime_{\stackrel{\scr z\le S}{\scr 2|z}}\frac{\mu^2(z)\mu(y)}{\phi(z)}H_2(2z)h_3(6z) \\
& \quad \ll \sum_{\stackrel{\scr z\le S}{\scr 2|z}}\frac{\mu^2(z)}{\phi(z)}H_2(2z)h_3(6z)\\
& \quad \ll \sum_{p\le S}\frac{\log p}{p}\sum_{\stackrel{\scr z\le S}{\scr p|z}}\frac{\mu^2(z)}{\phi(z)}H_2(2z) \\
& \quad \ll \sum_{3\le p\le S}\frac{\log p}{p\phi_2(p)}\sum_{\stackrel{\scr m\le\frac{ S}{p}}{\scr (m,p)=1}}\frac{\mu^2(m)}{\phi(m)}H_2(2m)\\
&\quad \ll \log S. \label{9.100}
\end{split}
\end{equation}
We conclude that
\begin{equation}\begin{split}
W_3(0) = &  \frac{1}{2}\log^2 S + O(\log S) + 2C_2 \sumprime_{\stackrel{\scr z\le S}{\scr 2|z}}\frac{\mu^2(z)\mu(y)}{\phi(z)}H_2(2z) \log(\min(v,w)) \\
&   +
O\Big(   \sumprime_{ z\le S} \frac{\mu^2(z)m(6z)}{\phi(z)\sqrt{\frac{S}{z}\min(v,w)}}H_2(2z)\Big).
\label{9.110}
\end{split}
\end{equation}
We now evaluate the first sum in \eqref{9.110}, and show that
\begin{equation}
2C_2 \sumprime_{\stackrel{\scr z\le S}{\scr 2|z}}\frac{\mu^2(z)\mu(y)}{\phi(z)}H_2(2z) \log(\min(v,w)) = \frac{1}{4}\log^2 S +O(\log S) . \label{9.120}
\end{equation}
By the symmetry in the variables $v$ and $w$ we have
that the sum above is
\begin{equation} = 4C_2\sumprime_{\stackrel{\stackrel{\scr vwy\le  S}{\scr 2|vwy}}{\scr v<w}}\frac{\mu^2(v)\mu^2(w)\mu(y)}{\phi(vwy)}H_2(2vwy) \log v,\label{9.130} \end{equation}
since the condition $(v,w)=1$ implies that the only term with $v=w$ is when $v=w=1$, and this term is zero.
We will sum over $w$ and apply Lemma \ref{Lemma2} according to the parity of $w$, which gives
\begin{eqnarray*} \sum_{\stackrel{\stackrel{\scr v< w \le  \frac{S} {vy}}{\scr 2|vwy}}{\scr (w,vy)=1}}\frac{\mu^2(w)}{\phi(w)}H_2(w) &=&
[(vy,2)=1]\sum_{\stackrel{\scr v< 2w' \le \frac{ S}{ vy}}{\scr (w',2vy)=1}}\frac{\mu^2(w')}{\phi_2(w')}+[2|vy]\sum_{\stackrel{\scr v< w \le \frac{ S}{ vy}}{\scr (w,2vy)=1}}\frac{\mu^2(w)}{\phi_2(w)} \\
&=& \Big([(vy,2)=1]+ [2|vy]\Big)\left(\frac{1}{\gs_2(2vy)}\log\frac{S}{v^2y}\right) +O(\frac{m(2vy)}{\sqrt{v}})\\
&=& \frac{1}{\gs_2(2vy)} \log\frac{S}{v^2y} +O(\frac{m(2vy)}{\sqrt{v}}).
\end{eqnarray*}
Substituting this result into
\eqref{9.130} we find that this sum is
\begin{equation*}\begin{split}
=& 4C_2\sum_{v\le \sqrt{S}}\frac{\mu^2(v)}{\phi(v)}H_2(v) \log v
\sum_{\stackrel{\scr y\le \frac{S}{v^2}}{\scr (y,v)=1}}\frac{\mu(y)}{\phi(y)}H_2(y)\left(\frac{1}{\gs_2(2vy)}\log\frac{S}{v^2y} +O(\frac{m(2vy)}{\sqrt{v}})\right) \\
=& 2\sum_{v\le \sqrt{S}}\frac{\mu^2(v)}{\phi(v)} \log v
\sum_{\stackrel{\scr y\le \frac{S}{v^2}}{\scr (y,v)=1}}\frac{\mu(y)}{\phi(y)}\log\frac{S}{v^2y} \\ & \qquad +O\Big(
\sum_{v\le \sqrt{S}}\frac{\mu^2(v)H_2(v)}{\sqrt{v}\phi(v)}m(v) \log v
\sum_{ y\le \frac{S}{v^2}}\frac{\mu^2(y)H_2(y)}{\phi(y)}m(y)\Big).
\end{split}
\end{equation*}
The inner sum in the error term here is
\begin{eqnarray*}
&\ll& \sum_{y\le S}\frac{\mu^2(y)H_2(y)}{\phi(y)}\sum_{d|y}\frac{\mu^2(d)}{\sqrt{d}} \\
&=& \sum_{d\le S}\frac{\mu^2(d)H_2(d)}{\sqrt{d}\phi(d)}\sum_{\substack{m \le S/d \\ (m,d)=1}}\frac{\mu^2(m)H_2(m)}{\phi(m)}\\
&\ll& \log S,
\end{eqnarray*}
and by \eqref{2.190} the sum over $v$ in the error term converges;  hence the error term is $O(\log S)$. Thus by Lemma \ref{Lemma1} and partial summation in Lemma \ref{Lemma4} the above expression is
\begin{eqnarray*} & =& 2\sum_{v\le \sqrt{S}}\frac{\mu^2(v)}{\phi(v)}\gs_2(v) \log v
       + O(\sum_{v\le \sqrt{S}}\left(\frac{\mu^2(v)}{\phi(v)}\log v\right) e^{-c_1\sqrt{\log(\frac{S}{v^2})}}) +O(\log S) \\ &=& \log^2\sqrt{S} + O(\log S) ,
\end{eqnarray*}
 which proves \eqref{9.120}.

 It remains to deal with the error term in \eqref{9.110}. By symmetry we may assume that $v\le w$, and thus the error term is
\begin{eqnarray*}& \ll& \frac{1}{\sqrt{S}}\sumprime_{\stackrel{\scr vwy\le S}{\scr (vwy,2)=1}}\frac{\mu^2(vwy)\sqrt{wy}}{\phi_2(vwy)}m(vwy) \\
&\ll &   \frac{1}{\sqrt{S}}\sum_{w\le S}\frac{\mu^2(w)\sqrt{w}}{\phi_2(w)}m(w)\sumprime_{\stackrel{\scr v\le \frac{S}{w}}{\scr (v,2)=1}}\frac{\mu^2(v)}{\phi_2(v)}m(v) \sumprime_{\stackrel{\scr y\le \frac{S}{vw}}{\scr (y,2)=1}}\frac{\mu^2(y)\sqrt{y}}{\phi_2(y)}m(y).
\end{eqnarray*}
The inner sum over $y$ is
\begin{eqnarray*} &\ll & \sum_{\stackrel{\scr y\le \frac{S}{vw}}{\scr (y,2)=1}}\frac{\mu^2(y)\sqrt{y}}{\phi_2(y)}\sum_{d|y}\frac{\mu^2(d)}{\sqrt{d}} \\
&\ll & \sum_{\stackrel{\scr d \le \frac{S}{vw}}{(d,2)=1}}\frac{\mu^2(d)}{\phi_2(d)}
\sum_{\stackrel{m\le \frac{S}{dvw}}{(m,2)=1}}\frac{\mu^2(m)\sqrt{m}}{\phi_2(m)} \\
&\ll &  \sqrt{\frac{S}{vw}}\sum_{\stackrel{\scr d \le \frac{S}{vw}}{(d,2)=1}}\frac{\mu^2(d)}{\phi_2(d)\sqrt{d}} \\
&\ll & \sqrt{\frac{S}{vw}}.
\end{eqnarray*}
Substituting we see the sum over $v$ now converges, and on summing over $w$ and treating $m(w)$ as in the previous estimate we see the error is $O(\log S)$. This completes the proof of \eqref{9.30}.

\section{ Evaluation of $S_3(0)$: Error Terms}

We now treat the error terms $\mathcal{E}_f(\mathcal{D})$, $\mathcal{E}_g(\mathcal{D}(D_1))$, and $\mathcal{E}_h(\mathcal{D}(D_2))$. We proceed as we did before in \eqref{8.30}; in $\mathcal{E}_f(\mathcal{D})$ we break the sum into two sums according to whether $b_{12}b_{23}a_{123}\le R_1$ or
$R_1<b_{12}b_{23}a_{123}\le R$, in the former sum we sum over $a_2$ using Lemma \ref{Lemma1} and obtain a sum of the same form as $\mathcal{E}_g(\mathcal{D}(D_1))$ and an error term
\begin{equation} \ll e^{-c_1\sqrt{\log(R/R_1)}}\log^8 R.\label{10.10}
\end{equation}
In  the second sum where $R_1 < b_{12}b_{23}a_{123}$ we have $a_2\le \frac{R}{R_1}$ and we do not sum over $a_2$. We continue this process with regard to $a_3$, and likewise deal with the error term $\mathcal{E}_g(\mathcal{D}(D_1))$. The result of this process is that we are left with  errors bounded by \eqref{10.10} and three types of sums of the forms
\begin{equation}
 \mathcal{E}_1 = \sum_{a\le R/R_1} \frac{\mu(a)}{\phi(a)}\sumprime_{\stackrel{\stackrel{\scr uvy \le R_1 } { \scr uwy\le R_1}}
{ R_1<\scr vwy\le R/a}} \frac {\mu(y)\mu^2(D)}{ \phi(D)} \gs_2(aD)\log\frac{R}{avwy},\label{10.20}
\end{equation}
\begin{equation}
 \mathcal{E}_2 = \sumprime_{a,b\le R/R_1} \frac{\mu(ab)}{\phi(ab)}\sumprime_{\stackrel{\stackrel{\scr uvy \le R_1 } {\scr R_1< uwy\le R/a}}
{\scr R_1< vwy\le R/b}} \frac {\mu(y)\mu^2(D)}{ \phi(D)} \log\frac{R}{auwy}\log\frac{R}{bvwy},\label{10.30}
\end{equation}
and
\begin{equation}
 \mathcal{E}_3 = \sumprime_{a,b,c\le R/R_1} \frac{\mu(abc)}{abc}\sumprime_{\stackrel{\stackrel{\scr R_1< uvy \le R/a } { \scr R_1< uwy\le R/b}}
{ R_1<\scr vwy\le R/c}} \frac {\mu(y)\mu^2(D)}{D} \log\frac{R}{auvy}\log\frac{R}{buwy}\log\frac{R}{cvwy},\label{10.40}
\end{equation}
where $D=uvwy$.

We can handle $\mathcal{E}_3$ immediately. Estimating trivially, we
have
\[  \mathcal{E}_3 \ll (\log R/R_1)^6 \sum_{\stackrel{\stackrel{\scr R_1< uvy \le R} {\scr R_1 <uwy\le R}}
{ R_1<\scr vwy\le R}} \frac{1}{uvwy}. \]
The top two inequalities in the summation conditions imply
\[ \left(\frac{R_1}{u}\right)^2 < vwy^2 \le \left(\frac{R}{u}\right)^2 \]
and hence
\[ \left(\frac{R_1}{u}\right)^2\frac{1}{vwy} < y \le \left(\frac{R}{u}\right)^2\frac{1}{vwy} \]
Thus the bottom inequality in the summation conditions implies
\[ \left(\frac{R_1}{u}\right)^2\frac{1}{R} < y \le \left(\frac{R}{u}\right)^2\frac{1}{R_1} ,\]
and therefore the sum above is
\begin{eqnarray*} &\ll & \sum_{u\le R}\frac{1}{u}\sum_{\left(\frac{R_1}{u}\right)^2\frac{1}{R} < y \le \left(\frac{R}{u}\right)^2\frac{1}{R_1}}\frac{1}{y} \sum_{\frac{R_1}{uy}<w\le \frac{R}{uy}}\frac{1}{w} \sum_{\frac{R_1}{uy}<v\le \frac{R}{uy} }\frac{1}{v} \\
&\ll & \log R (\log R/R_1)^3 .
\end{eqnarray*}
Thus
\begin{equation}
\mathcal{E}_3 \ll \log R(\log R/R_1)^9.
\label{10.50}
\end{equation}

Consider next $\mathcal{E}_2$. The trivial estimate used for $\mathcal{E}_3$ would give the bound $\ll \log^2 R(\log R/R_1)^6$, and therefore we need to save a factor of $\log R$, which will occur when we sum over $y$.
We first note that the conditions on the summation variables for the sum in \eqref{10.30} imply that
$u,v \le w$.  Next, we extend the summation range $uvy\le R_1$ to $uvy \le R$, which may be done with an error  $\ll \log R (\log R/R_1)^7$ in the same way that \eqref{10.50} was obtained. Finally, the terms with $R_1 < wy$ also contribute this same error, since this condition implies with the other summation conditions that $u,v \le \frac{R}{R_1}$ and $\frac{R_1}{uy} <w < \frac{R}{uy}$, so that only $y$ has a full summation range. Hence we have
\begin{equation} \mathcal{E}_2 = \sumprime_{a,b\le R/R_1} \frac{\mu(ab)}{\phi(ab)}\sumprime_{\substack{ u,v\le w \\ wy \le R_1  \\ R_1< uwy\le R/a \\  R_1< vwy\le R/b }} \frac {\mu(y)\mu^2(D)}{ \phi(D)} \log\frac{R}{auwy}\log\frac{R}{bvwy}+O\big(\log R(\log (R/R_1))^7\big).\label{10.60}
\end{equation}
We now sum over $u$, which satisfies
\[ \frac{R_1}{wy} < u \le \min(\frac{R}{awy},w) .\]
If $\min(\frac{R}{awy},w) =w$ then $w^2\le \frac{R}{ay}$ and there will only be terms when $\frac{R_1}{y} < w^2$. We conclude in this case that
\[  \sqrt{\frac{R_1}{y}} < w \le \sqrt{\frac{R}{ay}} .\]
Hence, as in the estimate to obtain \eqref{10.50}, these terms contribute at most $\ll \log R(\log(R/R_1))^7$ since only the variable $y$ runs through a full summation range.
We conclude, with $z=vwy$,
\begin{equation}\begin{split} \mathcal{E}_2 = \sumprime_{a,b\le R/R_1} \frac{\mu(ab)}{\phi(ab)}\sumprime_{\substack{ v\le w \\ wy\le R_1\\ R_1<z\le R/b}} \frac {\mu(y)\mu^2(z)}{ \phi(z)}\log\frac{R}{bz}&\sum_{\substack{\frac{R_1}{wy} < u \le \frac{R}{awy}\\ (u,abz)=1}} \frac{\mu^2(u)}{\phi(u)}\log\frac{R}{auwy}\\ &+O\big(\log R(\log (R/R_1))^7\big).\label{10.70}\end{split}
\end{equation}
To evaluate the inner sum, we use the relation, for $1\le S\le R$, $p(j)|k$,
\begin{equation}
\sum_{\substack{\frac{R}{S} < d \le R\\ (d,k)=1}} \frac{\mu^2(d)}{\phi_j(d)}\log \frac{R}{d} = \frac{1}{2\gs_j(k)} \log^2 S + O(\frac{m(k)}{\sqrt{\frac{R}{S}}}\log S).
\label{10.80} \end{equation}
This result follows immediately on writing the sum on the left-hand side above as
\begin{equation*}
= \sum_{\substack{d \le R\\ (d,k)=1}} \frac{\mu^2(d)}{\phi_j(d)}\log \frac{R}{d} -\sum_{\substack{d\le \frac{R}{S}\\ (d,k)=1}}\frac{\mu^2(d)}{\phi_j(d)}\big( \log \frac{R}{Sd} + \log S\big)
\end{equation*}
and applying \eqref{2.210} and Lemma \ref{Lemma2}.
Thus we have
\begin{equation}
\sum_{\substack{\frac{R_1}{wy} < u \le \frac{R}{awy}\\ (u,abz)=1}} \frac{\mu^2(u)}{\phi(u)}\log\frac{R}{auwy}
 = \frac{\phi(abz)}{2abz}\log^2\frac{R}{aR_1} + O(\frac{m(abz)}{\sqrt{\frac{R_1}{wy}}}\log\frac{R}{aR_1}). \label{10.90}
\end{equation}
Substituting this expression into \eqref{10.70} we obtain
\begin{equation} \mathcal{E}_2 = \frac{1}{2}\sumprime_{a,b\le R/R_1} \frac{\mu(ab)}{ab}\log^2\frac{R}{aR_1}\sumprime_{\substack{ v\le w \\ wy\le R_1\\ R_1< z\le R/b}} \frac {\mu(y)\mu^2(z)}{ z}\log\frac{R}{bz} +O\big(\log R(\log (R/R_1))^7\big),\label{10.100}
\end{equation}
once we show that the contribution
\[ \ll \frac{\log^2\frac{R}{R_1}}{\sqrt{R_1}}\sumprime_{a,b\le R/R_1} \frac{\mu^2(ab)m(ab)}{\phi(ab)}\sumprime_{\substack{ v\le w \\ wy\le R_1\\ R_1<z\le R/b}} \frac {\mu^2(z)\sqrt{wy}m(z)}{ \phi(z)}\]
from the error term in \eqref{10.90} is covered by the error term in \eqref{10.100}. This expression is of the same form as the error term in \eqref{9.110} estimated at the end of the last section, except it is over a more restricted summation range. The factors $m(ab)$ and $m(z)$ are handled as in that argument, and make no contribution, so we may ignore them. Hence the expression above is, by Lemma \ref{Lemma2},
\[ \begin{split}&\ll  \frac{\log^4\frac{R}{R_1}}{\sqrt{R_1}}\sum_{\substack{ v\le w \\ wy\le R_1\\ R_1<vwy\le R}} \frac {\mu^2(vwy)\sqrt{wy}}{ \phi(vwy)}\\
&\ll \frac{\log^4\frac{R}{R_1}}{\sqrt{R_1}} \sum_{w\le R_1}\frac{\sqrt{w}}{\phi(w)}\sum_{y\le \frac{R_1}{w}}\frac{\sqrt{y}}{\phi(y)}\sum_{\frac{R_1}{wy}< v \le \frac{R}{wy}}\frac{\mu^2(v)}{\phi(v)} \\
& \ll \log R_1 \log^5\frac{R}{R_1}, \end{split}\]
which is acceptable. Thus we have established \eqref{10.100}.

We now treat the sum in \eqref{10.100} and show it is also bounded by the error term, from which we conclude that
\begin{equation}\mathcal{E}_2 \ll  \log R (\log  R/R_1)^7. \label{10.110} \end{equation}
To see this, consider
the sum over $y$  in equation \eqref{10.100}
\[ \sum_{\substack{ \frac{R_1}{vw}<y\le \min(\frac{R_1}{w},\frac{R}{bvw})\\ (y,abvw)=1}}\frac{\mu(y)}{y}\log\frac{R}{bvwy}. \]
For the terms with $\min\big( \frac{R_1}{w}, \frac{R}{bvw}\big) = \frac{R_1}{w}$, we have $v\le \frac{R}{bR_1}$, and the sum in \eqref{10.100} is
\[ \begin{split} &\ll \log^5\frac{R}{R_1}\sum_{ v\le \frac{R}{R_1} } \frac{\mu^2(v)}{v}\sum_{ w\le R_1}\frac{\mu^2(w)}{w}\sum_{\frac{R_1}{vw}< y\le \frac{R_1}{w}} \frac {\mu^2(y)}{ y} \\ & \ll \log^5\frac{R}{R_1}\sum_{ v\le \frac{R}{R_1} } \frac{\mu^2(v)}{v}\sum_{ w\le R_1}\frac{\mu^2(w)}{w}\left([vw\le R_1]\log v +[vw>R_1]\log \frac{R_1}{w}\right) \\ &\ll \log R\log^5\frac{R}{R_1}\left(\sum_{ v\le \frac{R}{R_1} } \frac{\mu^2(v)\log v}{v} + \sum_{ v\le \frac{R}{R_1} } \frac{\mu^2(v)}{v}\sum_{ \frac{R_1}{v}< w\le R_1}\frac{\mu^2(w)}{w}\right) \\
&\ll \log R \log^7\frac{R}{R_1}, \end{split} \]
which is acceptable.
For the remaining terms when $\min\big( \frac{R_1}{w}, \frac{R}{bvw}\big) = \frac{R}{bvw}$,
we have $ \frac{R}{bR_1}< v$, and by Lemma \ref{Lemma1} with the error term estimate \eqref{2.140} the sum over $y$ is
\[ \begin{split} &= \sum_{\substack{ y\le \frac{R}{bvw}\\ (y,abvw)=1}}\frac{\mu(y)}{y}\log\frac{R}{bvwy}- [vw\le R_1]\sum_{\substack{ y\le \frac{R_1}{vw}\\ (y,abvw)=1}}\frac{\mu(y)}{y}\Big(\log\frac{R_1}{vwy}+\log\frac{R}{bR_1}\Big)\\
& = (1-[vw\le R_1] )\gs_1(abvw) + O(m(abvw)e^{-c_1\sqrt{\log(\frac{R}{bvw})}}) \\& \quad +[vw\le R_1]O(m(abvw) e^{-c_1\sqrt{\log(\frac{R_1}{vw})}}\log \frac{R}{R_1}). \end{split} \]
Hence the sum in \eqref{10.100} is in this case
\[ \begin{split} &\ll  \log^5(R/R_1) \sum_{vw\le R_1}\frac{m(vw)}{vw}e^{-c_1\sqrt{\log \frac{R_1}{vw}}}\\ & \quad + \log^4(R/R_1) \sum_{R_1<vw\le R}\frac{1}{vw}(\gs_1(vw) + m(vw))\\
&\ll\log R (\log(R/R_1))^5, \end{split}\]
where as before the factors of $m(vw)$ and $\gs(vw)$ make no contribution to the error when they are summed.
This finishes the proof of \eqref{10.110}.

Finally consider $\mathcal{E}_1$. The inner sum in \eqref{10.20} is, with $E=uwy$,
\begin{equation*}\sumprime_{\substack{ uvy \le R_1 \\ uwy\le R_1\\
 R_1< vwy\le R/a \\ (D,a)=1}}\frac {\mu(y)\mu^2(D)}{ \phi(D)} \gs_2(aD)\log\frac{R}{avwy}  = 2C_2\sumprime_{\substack{ E \le R_1 \\ (E,a)=1}} \frac {\mu(y)\mu^2(E)}{ \phi(E)} H_2(aE)S_R(E),
\end{equation*}
where
\[ S_R(E) =  \sum_{\substack{\frac{R_1}{wy}< v\le \min( \frac{R}{awy}, \frac{R_1}{uy}) \\ (v,aE)=1}}\chi(aE,v)\frac{\mu^2(v)H_2(v)}{\phi(v)}\log\frac{R}{avwy} \]
and
\[ \chi(aE,v) = [2|v][(aE,2)=1] +[(v,2)=1][2|aE]. \]
 We apply \eqref{10.80} and Lemma \ref{Lemma2} to evaluate $S_R(E)$. When  $\min( \frac{R}{awy}, \frac{R_1}{uy})= \frac{R}{awy}$
then in \eqref{10.80}   $S= \frac{R}{aR_1}$, and
\begin{equation*}
S_R(E) = \frac{1}{2\gs_2(2aE)}\log^2\frac{R}{aR_1}  +O(\frac{m(2aE)\log \frac{R}{aR_1}}{\sqrt{\frac{R_1}{wy}}});
\end{equation*}
while if $\min( \frac{R}{awy}, \frac{R_1}{uy})= \frac{R_1}{uy}$ then $1<S=\frac{w}{u}\le \frac{R}{aR_1}$, whence $u< w\le \frac{R}{aR_1}u$,
and thus we obtain
\begin{equation*}
S_R(E)= \frac{1}{2\gs_2(2aE)}\log\frac{w}{u} \log \left(\left(\frac{R}{aR_1}\right)^2\frac{u}{w}\right)  +O(\frac{m(2aE)\log \frac{R}{aR_1}}{\sqrt{\frac{R_1}{wy}}}).
\end{equation*}

On substituting, the inner sum in \eqref{10.20} becomes
\begin{equation*}\begin{split}
=&\frac{1}{2}\log^2\frac{R}{aR_1}\sumprime_{\substack{ uw \le R_1\\ \frac{R}{aR_1} \le \frac{w}{u} \\
  (uw,a)=1}} \frac {\mu^2(uw)}{ \phi(uw)}  \sum_{\substack{ y \le \frac{R_1}{uw} \\ (y,auw)=1}}\frac{\mu(y)}{\phi(y)} \\
&+ \frac{1}{2}\sumprime_{\substack{ uw \le R_1\\ 1< \frac{w}{u} <\frac{R}{aR_1} \\
  (uw,a)=1}} \frac {\mu^2(uw)}{ \phi(uw)}\log\frac{w}{u} \log \bigg(\Big(\frac{R}{aR_1}\Big)^2\frac{u}{w}\bigg)  \sum_{\substack{ y \le \frac{R_1}{uw} \\ (y,auw)=1}}\frac{\mu(y)}{\phi(y)}\\
&+ O\Bigg(\sumprime_{\substack{ uw \le R_1 \\
  (uw,a)=1}} \frac {\mu^2(uw)}{ \phi(uw)} \Big( \sum_{\substack{ y \le \frac{R_1}{uw} \\ (y,auw)=1}}\frac{\mu^2(y)H_2(2auwy)}{\phi^(y)} \frac{m(2auwy)\log \frac{R}{aR_1}}{\sqrt{\frac{R_1}{wy}}}\Big)\Bigg) .\end{split}
\end{equation*}
As before, in estimating the contribution of the error term above to \eqref{10.20} the factor $m(2auwy)$ may be ignored, and therefore this contribution is
\[ \begin{split} &\ll \frac{\log^2\frac{R}{R_1}}{\sqrt{R_1}}\sum_{\substack{uwy\le R_1\\ (uwy,2)=1}}\frac{\mu^2(uwy)\sqrt{wy}}{\phi_2(uwy)} \\
&\ll \frac{\log^2\frac{R}{R_1}}{\sqrt{R_1}}\sum_{\substack{u\le R_1\\ (u,2)=1}}\frac{\mu^2(u)}{\phi_2(u)} \sum_{\substack{w\le \frac{R_1}{u}\\ (w,2)=1}}\frac{\mu^2(w)\sqrt{w}}{\phi_2(w)}\sum_{\substack{y\le \frac{R_1}{uw}\\ (y,2)=1}}\frac{\mu^2(y)\sqrt{y}}{\phi_2(y)}\\
&\ll \log^2\frac{R}{R_1}\sum_{\substack{u\le R_1\\ (u,2)=1}}\frac{\mu^2(u)}{\sqrt{u}\phi_2(u)} \sum_{\substack{w\le \frac{R_1}{u}\\ (w,2)=1}}\frac{\mu^2(w)}{\phi_2(w)}\\
&\ll \log R_1 \log^2\frac{R}{R_1}. \end{split} \]
Next, for the main terms above the sum over $y$ is $\ll m(auw) e^{-c_2\sqrt{\log\frac{R_1}{uw}}}$ by Lemma \ref{Lemma1} with the error term in \eqref{2.140}, and hence both sums contribute to  \eqref{10.20}
\begin{eqnarray*}
 &\ll & \log^3\frac{R}{R_1} \sum_{u \le R_1}\frac{\mu^2(u)m(u)}{\phi(u)} \sum_{w\le R_1/u} \frac{\mu^2(w)m(w)}{\phi(w)}e^{-c_1\sqrt{\log \frac{R_1}{uw}}} \\
&\ll & \log R\log^3\frac{R}{R_1}.
\end{eqnarray*}
We conclude
\begin{equation}
\mathcal{E}_1 \ll \log R (\log R/R_1)^3 .
\label{10.120}
\end{equation}
By \eqref{10.10},\eqref{10.50},\eqref{10.110}, and \eqref{10.120} we have proved \eqref{9.40} and thus completed the proof of \eqref{5.30}.

\section{ Mixed Triple Correlations}
The case $k=2$ of Theorem \ref{Theorem3} has already been handled by \eqref{1.50} and \eqref{5.70}.
In this section we evaluate $W_R(\mbox{\boldmath$k$})$ from Section 4 in the case $k=3$. This will not only prove Theorem \ref{Theorem3}  but will also give an alternative and simpler proof of the second two parts of Theorem \ref{theorem5}.

We let $ \mbox{\boldmath$k$} =(k_1,k_2,0)$ and
\begin{equation}
W_R(\mbox{\boldmath$k$})= \sum_{\substack{d_1,d_2\le R\\ (d_1,k_1)=1,\, (d_2,k_2)=1\\(d_1,d_2)| k_2-k_1}}\frac{\mu(d_1)\mu(d_2)}{\phi([d_1,d_2])}\log\frac{R}{d_1}\log\frac{R}{d_2}.\label{11.10}
\end{equation}
The results we obtain are contained in the following theorem. Recall $H_j(m)$ is defined in \eqref{2.30}, $h_j(m)$ is defined in \eqref{2.170}, and $H_j(m),h_j(m)  \ll_j \log \log 3m$.
\begin{theorem} If $\mbox{\boldmath$k$} = (k,k,0)$, $k\neq 0$, and $\log|k|\ll \log R$,  we have
\begin{equation}
W_R(\mbox{\boldmath$k$})= \gs(\mbox{\boldmath$k$}) \log R +O(H_2(k)h_3(6k)), \label{11.20}
\end{equation}
and if  $ \mbox{\boldmath$k$} =(k_1,k_2,0)$ and $k_1\neq k_2\neq 0$, letting $\Delta = k_1k_2(k_2-k_1)$, and assuming $\Delta <R/2$, we have
\begin{equation}
W_R(\mbox{\boldmath$k$})= \gs(\mbox{\boldmath$k$}) + O(\frac{k_1}{\phi(k_1)}H_2(\Delta)e^{-c_1\sqrt{\log R/2\Delta}}).
\label{11.30}
\end{equation}
\label{Theorem11.1}
\end{theorem}

We decompose into relatively prime variables by letting $d_1 =a_1b_{12}$ and $d_2 = a_2b_{12}$ where $(d_1,d_2)=b_{12}$ and thus $a_1$, $a_2$, and $b_{12}$ are pairwise relatively prime. Then we have
\begin{equation}
W_R(\mbox{\boldmath$k$})= \sumprime_{\substack{a_1b_{12}\le R\\ a_2b_{12}\le R \\ (a_1b_{12},k_1)=1,\, (a_2b_{12},k_2)=1\\b_{12}| k_2-k_1}}\frac{\mu(a_1)\mu(a_2)\mu^2(b_{12})}{\phi(a_1a_2b_{12})}\log\frac{R}{a_1b_{12}}\log\frac{R}{a_2b_{12}}.\label{11.40}
\end{equation}
We first sum over $a_1$ and apply Lemma \ref{Lemma1} to see that
\begin{equation} \sum_{\substack{a_1\le R/b_{12}\\ (a_1,a_2b_{12}k_1)=1}}\frac{\mu(a_1)}{\phi(a_1)}\log\frac{R}{a_1b_{12}} = \gs_2(a_2b_{12}k_1) + r_1(\frac{R}{b_{12}}, a_2b_{12}k_1),\label{11.50}
\end{equation}
and hence we have
\begin{equation}\begin{split}
W_R(\mbox{\boldmath$k$})&= \sumprime_{\substack{ a_2b_{12}\le R \\ (b_{12},k_1)=1,\, (a_2b_{12},k_2)=1\\b_{12}| k_2-k_1}}\frac{\mu(a_2)\mu^2(b_{12})\gs_2(a_2b_{12}k_1)}{\phi(a_2b_{12})}\log\frac{R}{a_2b_{12}}+ E_1(R) \\
&=V_R(\mbox{\boldmath$k$})+E_1(R),\label{11.60}
\end{split}
\end{equation}
where
\begin{eqnarray}
E_1(R) &=& \sumprime_{\substack{ a_2b_{12}\le R \\ (b_{12},k_1)=1,\, (a_2b_{12},k_2)=1\\b_{12}| k_2-k_1}} \frac{\mu(a_2)\mu^2(b_{12}) r_1(\frac{R}{b_{12}}, a_2b_{12}k_1)}{\phi(a_2b_{12})}\log\frac{R}{a_2b_{12}} \nonumber \\
& \ll& \sum_{\substack{b_{12}\le R\\ b_{12}|k_2-k_1}} \frac{\mu^2(b_{12})}{\phi(b_{12})}e^{-c_1\sqrt{\log(R/b_{12})}}\sum_{a_2\le R/b_{12}}\frac{\mu^2(a_2)}{\phi(a_2)}\log \frac{R}{a_2b_{12}} \nonumber \\
&\ll& \sum_{\substack{b_{12}\le R\\ b_{12}|k_2-k_1}} \frac{\mu^2(b_{12})}{\phi(b_{12})}e^{-c_2\sqrt{\log(R/b_{12})}}. \label{11.70} \end{eqnarray}

We first consider the case when $k_1=k_2=k$. Then from \eqref{11.70} we see that
\begin{equation} E_1(R)\ll 1 .\label{11.80} \end{equation}
By Lemma \ref{Lemma3} we have
\begin{eqnarray}
V_R(\mbox{\boldmath$k$})&=& \sumprime_{\substack{a_2b_{12}\le R \\  (a_2b_{12},k)=1}}\frac{\mu(a_2)\mu^2(b_{12})\gs_2(a_2b_{12}k)}{\phi(a_2b_{12})}\log\frac{R}{a_2b_{12}}\nonumber  \\
&=& \sum_{\substack{b_{12}\le R  \\  (b_{12},k)=1}}\frac{\mu^2(b_{12})}{\phi(b_{12})}\sum_{\substack{ a_2\le R/b_{12}\\ (a_2,b_{12}k)=1}}\frac{\mu(a_2)\gs_2(a_2b_{12}k)}{\phi(a_2)}\log\frac{R}{a_2b_{12}} \nonumber \\
&=& -\sum_{\substack{b_{12}\le R  \\  (b_{12},k)=1}}\frac{\mu^2(b_{12})\mu((b_{12}k,2))}{\phi(b_{12})}\gs_2(2b_{12}k)\gs_3(2b_{12}k)\nonumber \\ & & \qquad  + \ O\left( \sum_{b_{12}\le R}\frac{\mu^2(b_{12})H_2(b_{12})}{\phi(b_{12})\log^A(2R/b_{12})}\right) \nonumber \\
&=& Y_R(k) +O(1). \label{11.90}
\end{eqnarray}
To evaluate $Y_R(k)$ we divide into even terms and odd terms and apply Lemma \ref{Lemma4}. Thus
\begin{eqnarray*}
V_R(k) &=& [(k,2)=1]\sum_{\substack{2b\le R \\  (b,2k)=1}}\frac{\mu^2(b)}{\phi(b)}\gs_2(2bk)\gs_3(2bk) \\
& & \qquad -\mu((k,2)) \sum_{\substack{b_{12}\le R \\  (b_{12},2k)=1}}\frac{\mu^2(b_{12})}{\phi(b_{12})}\gs_2(2b_{12}k)\gs_3(2b_{12}k)  \\
&=& 2C_2H_2(k)\Bigg([(k,2)=1]\sum_{\substack{b\le R/2 \\  (b,2k)=1}}\frac{\mu^2(b)}{\phi_2(b)}\gs_3(2bk) \\ & & \hskip 1.5in  -\mu((k,2))\sum_{\substack{b_{12}\le R \\  (b_{12},2k)=1}}\frac{\mu^2(b_{12})}{\phi_2(b_{12})}\gs_3(2b_{12}k)
\Bigg)\\
&=& 2C_2H_2(k)\bigg(\Big( [(k,2)=1] - \mu((k,2))\Big)\log R +O(h_3(6k)) +O\big(\frac{\gs_3(6k)m(k)}{\sqrt{R}}\big) )\bigg)\\
&=&  2C_2[2|k]H_2(k)\log R +O(H_2(k)h_3(6k))\\
&=& \gs_2(k) \log R +O(H_2(k)h_3(6k)).
\end{eqnarray*}
On combining these results we have proved the first part of Theorem \ref{Theorem11.1}.

We now turn to the case that $k_2\neq k_1$. As in Lemma \ref{Lemma5} we let $\kappa = (k_1,k_2)$ and $\Delta = k_1k_2(k_2-k_1)$. As before let $k^* = \max (|k_1|,|k_2|)$. In this case we see from \eqref{11.70} that
\begin{equation} E_1(R) \ll  H_1(k_2-k_1) e^{-c_2\sqrt{\log(R/2k^*)}}. \label{11.100}
\end{equation}
We now let $k_1=s_1K_1$, $k_2=s_2K_2$, where $K_1$ and $K_2$ are the largest squarefree divisors of $k_1$ and $k_2$ respectively. Let $K_{12}=(K_1,K_2)$. Then we have
\begin{equation*} V_R(\mbox{\boldmath$k$}) =\sumprime_{\substack{ a_2b_{12}\le R \\ (b_{12},K_1)=1,\, (a_2b_{12},K_2)=1\\b_{12}| k_2-k_1}}\frac{\mu(a_2)\mu^2(b_{12})\gs_2(a_2b_{12}K_1)}{\phi(a_2b_{12})}\log\frac{R}{a_2b_{12}}
\end{equation*}
Next let $a_2=c_2d_2$  where $c_2|K_1$ and  $(d_2,K_1)=1$ so that $(a_2,K_1)=c_2$, from which we see
\begin{equation*} V_R(\mbox{\boldmath$k$}) =\sumprime_{\substack{ c_2d_2b_{12}\le R \\ (d_2b_{12},K_1)=1,\, (c_2d_2b_{12},K_2)=1\\b_{12}| k_2-k_1, \, c_2|K_1}}\frac{\mu(c_2)\mu(d_2)\mu^2(b_{12})\gs_2(d_2b_{12}K_1)}{\phi(c_2d_2b_{12})}\log\frac{R}{c_2d_2b_{12}}
\end{equation*}
On summing over $d_2$ by dividing the sum according to whether $d_2$ is even or odd, we see  on applying Lemma \ref{Lemma1} that
\begin{equation*}\begin{split}
\sum_{\substack{ d_2\le R/c_2b_{12} \\ (d_2,b_{12}K_1K_2)=1}}&\frac{\mu(d_2)\gs_2(d_2b_{12}K_1)}{\phi(d_2)}\log\frac{R}{c_2d_2b_{12}}\\
&=  2C_2H_2(b_{12}K_1)\sum_{\substack{ d_2\le R/c_2b_{12} \\ (d_2,b_{12}K_1K_2)=1\\ 2|d_2b_{12}K_1}}\frac{\mu(d_2)}{\phi(d_2)}H_2(d_2)\log\frac{R}{c_2d_2b_{12}}\\
&= 2C_2H_2(b_{12}K_1)\Bigg(-[(K_1K_2b_{12},2)=1] \sum_{\substack{ 2d\le R/c_2b_{12} \\ (d,2b_{12}K_1K_2)=1}}\frac{\mu(d)}{\phi_2(d)}\log\frac{R}{2c_2db_{12}}\\
& \qquad \qquad + [2|b_{12}K_1]\sum_{\substack{ d_2\le R/c_2b_{12} \\ (d_2,2b_{12}K_1K_2)=1}}\frac{\mu(d_2)}{\phi_2(d_2)}\log\frac{R}{c_2d_2b_{12}}\Bigg)\\
&=2C_2H_2(b_{12}K_1)\chi(b_{12})\gs_3(2b_{12}K_1K_2)+O(H_2(b_{12}K_1)e^{-c_1\sqrt{\log(R/2c_2b_{12})}}),
\end{split}
\end{equation*}
where
\[  \chi(b_{12}) = -[(K_1K_2b_{12},2)=1] +[2|b_{12}K_1]. \]
Since $c_2b_{12}|K_1(k_2-k_1)$, we see that the condition $|\Delta| < R/2$ implies $c_2b_{12}\le R/2$ is automatically satisfied. Therefore we  have
\begin{equation} \begin{split}
V_R(\mbox{\boldmath$k$}) &=2C_2\sum_{\substack{ b_{12}| k_2-k_1, \, c_2|K_1 \\ (b_{12},K_1K_2)=1,\, (c_2,K_2)=1}}\frac{\mu(c_2)\mu^2(b_{12})H_2(b_{12}K_1)\chi(b_{12})\gs_3(2b_{12}K_1K_2)}{\phi(c_2)\phi(b_{12})} + E_2(R)  \\
&= Y_R(\mbox{\boldmath$k$})    +E_2(R), \label{11.110}
\end{split}
\end{equation}
where
\begin{equation*}  E_2(R) \ll \frac{k_1}{\phi(k_1)}H_2(\Delta)e^{-c_1\sqrt{\log R/2\Delta}}.
\end{equation*}
Since
\[  \sum_{\substack{c_2|K_1\\ (c_2,K_2)=1}}\frac{\mu(c_2)}{\phi(c_2)}=\prod_{p|\frac{K_1}{K_{12}}}\left( 1 - \frac{1}{p-1}\right)= [(K_1/K_{12},2)=1]\frac{H_2(K_{12})}{H_2(K_1)},\]
we have
\begin{equation}
Y_R(\mbox{\boldmath$k$}) =2C_2[(K_1/K_{12},2)=1]\frac{H_2(K_{12})}{H_2(K_1)}\sum_{\substack{ b_{12}| k_2-k_1\\ (b_{12},K_1K_2)=1 }}\frac{\mu^2(b_{12})H_2(b_{12}K_1)\chi(b_{12})\gs_3(2b_{12}K_1K_2)}{\phi(b_{12})}. \label{11.120}
\end{equation}
We divide the sum in the equation above according to the parity of $b_{12}$ and apply Lemma \ref{Lemma4} to see that the sum is
\begin{equation*} \begin{split}
& =[(K_1K_2,2)=1]\sum_{\substack{ 2b| k_2-k_1\\ (b,2K_1K_2)=1 }}\frac{\mu^2(b)H_2(bK_1)\chi(2b)\gs_3(2bK_1K_2)}{\phi(b)} \\
& \qquad  + \sum_{\substack{ b_{12}| k_2-k_1\\ (b_{12},2K_1K_2)=1 }}\frac{\mu^2(b_{12})H_2(b_{12}K_1)\chi(b_{12})\gs_3(2b_{12}K_1K_2)}{\phi(b_{12})} \\
& =\eta(K_1,K_2)H_2(K_1)\sum_{\substack{ b| k_2-k_1\\ (b,2K_1K_2)=1 }}\frac{\mu^2(b)\gs_3(2bK_1K_2)}{\phi_2(b)} \\
& = \eta(K_1,K_2)H_2(K_1)
\gs_3(2K_1K_2(k_2-k_1)),
\end{split}
\end{equation*}
where
\[ \eta(K_1,K_2) = [(K_1K_2,2)=1][2|k_2-k_1] -[(K_1K_2,2)=1] +[2|K_1] \]
On substituting we have that
\begin{equation*}
Y_R(\mbox{\boldmath$k$}) =2C_2[(K_1/K_{12},2)=1]\eta(K_1,K_2)H_2(K_{12})
\gs_3(2K_1K_2(k_2-k_1)).
\end{equation*}
Now $[(K_1/K_{12},2)=1] =0$ if  $K_1$ is even and  $K_2$ is odd, and $\eta(K_1,K_2) =0$ if both $K_1$ and $K_2$ are odd or $K_1$ is odd and $K_2$ is even. Hence $Y_R(\mbox{\boldmath$k$})$ is zero unless $K_{12}$ is even, and therefore
\begin{equation*}
Y_R(\mbox{\boldmath$k$}) =[2|K_{12}]2C_2H_2(K_{12})
\gs_3\big(2K_1K_2(k_2-k_1)\big) = \gs_2(\kappa)\gs_3(\Delta)=\gs(\mbox{\boldmath$k$}),
\end{equation*}
by Lemma \ref{Lemma5}.
This completes the proof of Theorem \ref{Theorem11.1}

\section{Application to primes}

The use of correlations of short divisor sums to study primes goes back at least to Selberg's work on the sieve. Our mixed correlation result that,
for $ R \le N^{\frac{1}{4} -\epsilon}$,
\begin{equation}
\tilde{\mathcal{S}}_3(N,k) =  \sum_{n=1}^N {\Lambda_R}^2(n)\Lambda(n+k)
= \gs_2(k) N \log R +o(N\log N) \label{12.10}
\end{equation}
provides the upper bound for prime pairs  in \eqref{1.390} with $\mathcal{B}=4$, since for $n\ge R$
\[ \mu^2(n) \Lambda(n) \le \frac{\log n}{\log^2 R}\Lambda_R^2(n) ,\]
and the prime powers make a contribution $\ll N^{1/2 + \epsilon}$.
The Selberg sieve provides the same information, and while the optimal majorant obtained with the Selberg sieve is different from $\Lambda_R(n)$ (and also $\lambda_R(n)$ in \eqref{1.20}), there is nothing lost asymptotically in the use of $\Lambda_R(n)$.

To study primes in short intervals, we consider the modified moments
\begin{equation}
M_k'(N, h,\psi_R,C) =  \sum_{n=N+1}^{2N}
\big(\psi_R(n+h)-\psi_R(n)- C\log N\big)^k ,\label{12.20}
\end{equation}
and
\begin{equation}
\tilde{M}_k'(N, h,\psi_R,C) =  \sum_{n=N+1}^{2N}\big(\psi(n+h)-\psi(n)\big)
\big(\psi_R(n+h)-\psi_R(n)- C\log N\big)^{k-1} ,\label{12.30}
\end{equation}
where $C$ is a function of $h$ and $R$ that will be chosen to optimize our applications.  If we take $C=0$ these moments reduce to the moments considered in Section 1. We will assume in this section that $h = \lambda \log N$,  $\lambda \ll 1$, and thus
\begin{equation} h \ll \log N ,\label{12.40} \end{equation}
which we will make free use of in our estimates.
We now consider, for $\rho\ge 0$,
\begin{equation} \begin{split} \mathcal{M}(h,\rho) &= \tilde{M}_3'(N, h,\psi_R,C) - \rho\log N M_2'(N,h,\psi_R, C)  \\ &=  \sum_{n=N+1}^{2N}\big(\psi(n+h)-\psi(n) -\rho\log N \big)
\big(\psi_R(n+h)-\psi_R(n)- C\log N\big)^{2} .
\label{12.50}
\end{split}\end{equation}
To evaluate $\mathcal{M}(h,\rho)$ we see first that
\[
\tilde{M}_3'(N, h,\psi_R,C) = \tilde{M}_3'(N, h,\psi_R)-2C\log N\tilde{M}_2'(N, h,\psi_R) +C^2\log^2N\tilde{M}_1'(N, h,\psi_R). \]
We apply  Corollary \ref{Corollary2} (which as mentioned in Section 1 applies immediately to $\tilde{M}_k'$ as well as $\tilde{M}_k$), with $R=N^\theta$ and $0< \theta <\frac{\vartheta}{2}$,
\[ \tilde{M}_3'(N, h,\psi_R,C) \sim \Big( (\theta^2\lambda +3\theta \lambda^2 + \lambda^3) -2C(\theta\lambda + \lambda^2)+ C^2\lambda \Big) N\log^3N \]
and by Corollary \ref{Corollary1}, for $0< \theta \le \frac{1}{2}$,
\[ M_2'(N, h,\psi_R,C) \sim \Big( \theta\lambda +\lambda^2 -2C\lambda +C^2\Big) N\log^2N . \]
We therefore see that $\mathcal{M}(h,\rho)$ is quadratic in $C$ when $\lambda \neq \rho$, and therefore on completing the square we find that, for $\lambda \neq \rho$, and $0<\theta < \frac{\vartheta}{2}$,
\[ \mathcal{M}(h,\rho)\sim \Big((\lambda -\rho)\Big( C - \frac{ \lambda(\lambda - \rho +\theta)}{\lambda -\rho}\Big)^2 + \frac{\lambda \theta}{\lambda -\rho}\big( (\lambda-\rho)^2-\theta \rho\big)\Big)N\log^3N .\]
By choosing
\begin{equation}C = \frac{ \lambda(\lambda - \rho +\theta)}{\lambda -\rho}= \lambda\Big( 1 + \frac{\theta}{\lambda-\rho}\Big)
\label{12.60} \end{equation}
we maximize $\mathcal{M}(h,\rho)$ if $\lambda <\rho$ and minimize it if $\lambda >\rho$. We conclude that with this choice of $C$, and $0<\theta < \frac{\vartheta}{2}$,
\begin{equation}\mathcal{M}(h,\rho)
\sim  \frac{\lambda \theta}{\lambda -\rho}\big( (\lambda-\rho)^2-\theta \rho\big)N\log^3N.
\label{12.70} \end{equation}
We see that $\mathcal{M}(h,\rho)$ is positive (and $\gg N\log^3 N$ ) when  $\lambda$ is a fixed number in  the range $ \rho - \sqrt{\theta \rho} < \lambda < \rho$ but is negative when $\rho <\lambda < \rho + \sqrt{\theta \rho}$.

We now let $P_r(N,h)$ denote the number of integers $N<n\le 2N$ for which the interval $(n, n+h]$ contains exactly $r$ primes. Thus
\begin{equation}
P_r(N,h) = \sum_{\substack{n=N+1 \\ \pi(n+h)-\pi(n) = r}}^{2N}1 .
\label{12.80}
\end{equation}
The Poisson model for primes in short intervals (see \cite{GA}) is equivalent to the conjecture that
\begin{equation}
P_r(N,h) \sim \frac{\lambda^r e^{-\lambda}}{r!} N .
\label{12.90}
\end{equation}
We  let
\begin{equation} Q_r^-(N,h) = \sum_{m=0}^r P_m(N,h) = \sum_{\substack{n=N+1\\ \pi(n+h)-\pi(n)\le r}}^{2N} 1
\label{12.100}
\end{equation}
and
\begin{equation} Q_r^+(N,h) = \sum_{m=r+1}^\infty P_m(N,h) = \sum_{\substack{n=N+1\\ \pi(n+h)-\pi(n)> r}}^{2N} 1.
\label{12.110}
\end{equation}
Thus we have
\begin{equation}
Q_{r}^-(N,h) + Q_{r}^+(N,h) = N.
\label{12.120}
\end{equation}
We let $p_{j_0}=p_{j_0}(N)$ and $p_{j_1}=p_{j_1}(N)$ denote respectively the smallest and the largest primes in the interval $[N+1,2N]$.
For smaller than average gaps between primes, we use the relation, for $r\ge 1$,
\[Q_r^+(N,h) =\sum_{\substack{N+1\le n<p_{j_0}\\ \pi(n+h)-\pi(n)>r}}1 + \sum_{j=j_0}^{j_1} \sum_{\substack{p_j\le n < p_{j+1}\\p_{j+r+1}\le n+h}}1 - \sum_{\substack{2N< n <p_{j_1+1}\\ \pi(n+h)-\pi(n)>r}}1.\]
The first and third sums are $O(Ne^{-c_1\sqrt{\log N}})$ by the prime number theorem with error term, and hence
\begin{eqnarray}
Q_r^+(N,h) &=&\sum_{N+1\le p_j \le 2N} \sum_{\substack{p_j\le n < p_{j+1}\\p_{j+r+1}\le n+h}}1 + O(Ne^{-c_1\sqrt{\log N}})\nonumber \\
&=& \sum_{\substack{N+1\le p_j \le 2N \\ p_{j+r+1}-p_{j+1}< h}}\big( p_{j+1}-\max(p_j,p_{j+r+1}- \lfloor h\rfloor) \big)+ O(Ne^{-c_1\sqrt{\log N}} )\nonumber \\
&\le & h \sum_{\substack{N+1\le p_j \le 2N \\ p_{j+r+1}-p_{j+1}< h}}1 +O(Ne^{-c_1\sqrt{\log N}}).
\label{12.130}
\end{eqnarray}
For larger than average gaps between primes a similar argument shows, for $r\ge 0$,
\begin{eqnarray}
Q_{r}^-(N,h) &=&\sum_{N+1\le p_j \le 2N} \sum_{\substack{p_j\le n < p_{j+1}\\p_{j+r+1}> n+h}}1 +O(Ne^{-c_1\sqrt{\log N}}) \nonumber \\
&=& \sum_{\substack{N+1\le p_j \le 2N \\ p_{j+r+1}-p_j > h}}\big( \min(p_{j+1},p_{j+r+1}- \lfloor h\rfloor) - p_j \big)+O(Ne^{-c_1\sqrt{\log N}}) \nonumber \\
&\le &  \sum_{\substack{N+1\le p_j \le 2N \\ p_{j+r+1}-p_{j}> h}}(p_{j+r+1}-p_j)+O(Ne^{-c_1\sqrt{\log N}}).
\label{12.140}
\end{eqnarray}
Next, we have
\begin{equation}
Q_r^+(N,h) = (1 +o(1)) \sum_{\substack{n=N+1\\ \psi(n+h)-\psi(n) \ge \rho \log N}}^{2N} 1
\label{12.150}
\end{equation}
where $\rho$ can be taken to be any number in the range $r<\rho < r+1$, since the prime powers may be discarded with an error $\ll N^{\frac{1}{2}}$. Also
 \begin{equation}
Q_r^-(N,h) = (1 +o(1)) \sum_{\substack{n=N+1\\ \psi(n+h)-\psi(n) \le \rho \log N}}^{2N} 1
\label{12.160}
\end{equation}
where again $\rho$ can be taken to be any number in the range $r<\rho < r+1$. Returning to  \eqref{12.50}, we see on applying Cauchy's inequality twice and using \eqref{12.150} that
\begin{equation}\begin{split}
\mathcal{M}(h,\rho) &\le \sum_{\substack{n=N+1\\ \psi(n+h)-\psi(n) \ge \rho\log N }}^{2N}\big(\psi(n+h)-\psi(n)\big)
\big(\psi_R(n+h)-\psi_R(n)- C\log N\big)^{2}  \\
&\le \left(\sum_{\substack{n=N+1\\ \psi(n+h)-\psi(n) \ge \rho \log N}}^{2N} 1
\right)^{\frac{1}{4}}\left(\sum_{n=N+1}^{2N}(\psi(n+h)-\psi(n))^4\right)^{\frac{1}{4}}  \\  & \qquad \left(\sum_{n=N+1}^{2N}(\psi_R(n+h)-\psi_R(n)- C\log N)^4\right)^{\frac{1}{2}}  \\
&=  (1 + o(1))Q_r^+(N,h)^{\frac{1}{4}}M_4'(N,h,\psi)^{\frac{1}{4}}M_4'(N,h,\psi_R,C)^{\frac{1}{2}}.
\label{12.170}
\end{split}
\end{equation}
The same argument also shows that
\begin{equation}
-\mathcal{M}(h,\rho) \le (1 + o(1))\rho \log N Q_r^-(N,h)^{\frac{1}{2}}M_4'(N,h,\psi_R,C)^{\frac{1}{2}},
\label{12.180}
\end{equation}
and therefore we conclude that for any $r<\rho <r+1$
\begin{equation}\begin{split}
- (1 + o(1))Q_r^-(N,h)^{\frac{1}{2}}& \le \frac{\mathcal{M}(h,\rho)}
{\rho \log N M_4'(N,h,\psi_R,C)^{\frac{1}{2}}}\\ &\qquad  \le
(1 + o(1))Q_r^+(N,h)^{\frac{1}{4}}\left(\frac {M_4'(N,h,\psi)^{\frac{1}{4}}}{\rho \log N}\right) \end{split}
\label{12.190}
\end{equation}
To prove the first part of Theorem \ref{Theorem5} we estimate the moments
$M_4'(N,h,\psi)$ and  $M_4'(N,h,\psi_R,C)$
trivially when  $h \gg 1$ using the inequality
\[|abcd| \le \frac{1}{4}(a^4+b^4+c^4+d^4)\]
 and the equation above \eqref{1.100} to see that
\[ \begin{split} M_4'(N,h,\psi) &= \sum_{\substack{1\le m_i\le h\\ 1\le i\le 4}} \sum_{n=N+1}^{2N} \Lambda(n+m_1)\Lambda(n+m_2)\Lambda(n+m_3)\Lambda(n+m_4) \\
&\ll h^4 \sum_{n\le 3N}\Lambda(n)^4 \\
&\ll h^4N\log^3N, \end{split}\]
and similarly
\begin{eqnarray*} M_4'(N,h,\psi_R,C) &\ll & M_4'(N,h,\psi_R) + N\log^4N \\
&\ll & h^4 \sum_{n\le 3N}\Lambda_R(n)^4
+ N\log^4 N \\  &\ll & h^4 \log^4 R\sum_{n\le 3N}d(n)^4  \\
&\ll & h^4 N \log^{19} N
\end{eqnarray*}
by \eqref{4.140}.  Hence, subject to \eqref{12.40}, we see by \eqref{12.70}, \eqref{12.130}, and (12.19) that, for $r\ge 1$, and some positive constant $C$,
\begin{equation}
\sum_{\substack{N+1\le p_j \le 2N \\ p_{j+r+1}-p_{j+1}< h}}1 \gg
\frac{N}{\log^C N}
\label{12.200}
\end{equation}
provided $\rho - \sqrt{\theta\rho}< \lambda < \rho $, $r<\rho <r+1$, and
$0<\theta < \frac{\vartheta}{2}$. Since $\rho$ can be taken as close to $r$ as we wish, we conclude that
\[ \Xi_r \le r - \sqrt{\theta r}, \]
where unconditionally we may take any $0< \theta < 1/4$.  This proves the first part of Theorem \ref{Theorem5}. If we assume $\vartheta =1$ we can take $0<\theta < \frac{1}{2}$. The corresponding result for larger than average gaps between primes is proved in the same way.

In order to obtain positive proportion results, we need to use the generalization of the sieve upper bound \eqref{1.390} for prime $k$-tuples.
This result states that for the function $\psi_{\mbox{\boldmath$j$}}(N)$ defined in \eqref{1.80} where ${\mbox{\boldmath$j$}}=(j_1,j_2,  \ldots , j_r)$ with the $j_i$'s distinct and $\gs(\mbox{\boldmath$j$})\neq 0$
\begin{equation}
\psi_{\mbox{\boldmath$j$}}(N) \le (2^rr!+\epsilon)\gs(\mbox{\boldmath$j$})N ,
\label{12.210}
\end{equation}
 see Theorem 5.7 of \cite{HR}.  On applying this bound to the formulas leading to \eqref{1.150} we see that, subject to \eqref{12.40},
\begin{equation}
M_k(N,h,\psi) \ll N(\log N)^k,
\label{12.220}
\end{equation}
which implies the same estimate holds for $M_k'(N,h,\psi)$.
Next, as above
\[ M_4'(N,h,\psi_R,C) \ll M_4'(N,h,\psi_R) + N\log^4 N, \]
and therefore assuming \eqref{1.470} we have, for $0< \theta < \frac{1}{4}$,
\begin{equation}
M_4'(N,h,\psi_R,C) \ll  N\log^4 N.
\label{12.230}
\end{equation}
Using these estimates in \eqref{12.190} we obtain
\begin{equation}
\sum_{\substack{N+1\le p_j \le 2N \\ p_{j+r+1}-p_{j+1}< h}}1 \gg
\frac{N}{\log N} \label{12.240}
\end{equation}
under the same conditions as \eqref{12.200} and $0<\theta <\frac{1}{4}$.
This proves the remaining part of Theorem \ref{Theorem5}.

In an identical fashion we see that if $\mathcal{M}(h,\rho)<0$
then assuming \eqref{1.470} we have
\begin{equation}
\sum_{\substack{N+1\le p_j \le 2N \\ p_{j+r+1}-p_{j}> h}}(p_{j+r+1}-p_j) \gg
 N,
\label{12.250}
\end{equation}
where $r\ge 0$,  $r<\rho<r+1$,  $ \rho < \lambda <\rho +\sqrt{\theta\rho}$, and $0<\theta < \frac{1}{4}$. Since $\rho$ can be taken as close to $r+1$ and $\theta$ as close to $\frac{1}{4}$ as we wish, this completes the proof of Theorem \ref{Theorem6}.

\section*{ Acknowledgement}

A preliminary version of this paper was presented at the First Workshop on $L$-functions and Random Matrices at the American Institute of Mathematics in May, 2001. Thanks to a suggestion of Peter Sarnak during the talk, and encouragement of John Friedlander after the talk, the authors found an alternative method for proving Theorem \ref{Theorem1} which generalizes to the case of $k$-correlations. This new method has its own complications which make it very easy for mistakes to creep into the calculations. Following the conference, the first-named author visited AIM where these problems were solved jointly with Brian Conrey and David Farmer. Farmer has written a program in Mathematica to compute the constants $\mathcal{ C}_k(\mbox{\boldmath$a$})$ and has found that $\mathcal{C}_4(4) = \frac{3}{4}$, $\mathcal{C}_5(5) = \frac{11065}{2^{14}} = .67535\ldots$, and $\mathcal{C}_6(6)= \frac{11460578803}{2^{34}}= .66709\ldots$. The authors would like to thank these individuals and the American Institute of Mathematics.